\documentclass[10pt,reqno]{amsart}

\usepackage{amsmath,amstext,amsfonts,%
enumitem,%
extpfeil,graphicx,mathrsfs,%
mathtools,setspace,tikz,xparse}
\usepackage[alphabetic]{amsrefs}
\usepackage[unicode]{hyperref}
\usepackage[utf8]{inputenc}
\usepackage[T1]{fontenc}
\usepackage[normalem]{ulem}
\usepackage[cmtip,all]{xy}
\usepackage[light]{sourcesanspro}
\usepackage{lmodern}
\hypersetup{colorlinks=false,linkbordercolor={0.93 0.71 0.13},citebordercolor={1 0 .5},urlbordercolor={0.27 0.39 0.68}}
\usetikzlibrary{calc,matrix}

\definecolor{MyColor1}{HTML}{fff3b0}
\definecolor{MyColor2}{HTML}{e09f3e}
\definecolor{MyColor3}{HTML}{335c67}
\definecolor{MyColor4}{HTML}{393e41}

\usepackage{filecontents}

\DeclareMathOperator{\End}{End}

\DeclareMathOperator{\Ker}{ker}

\DeclareMathOperator{\Supp}{supp}
\DeclareMathOperator{\GL}{GL}
\DeclareMathOperator{\Sym}{Sym}
\DeclareMathOperator{\Img}{im}
\DeclareMathOperator{\Ind}{Ind}
\DeclareMathOperator{\ind}{ind}

\newcommand{\m}[1]{\(\smash{#1}\)}

\newcommand{\bb}[1]{\mathbb{#1}}
\newcommand{\QQ}{\bb{Q}}
\newcommand{\FF}{\bb{F}}
\newcommand{\ZZ}{\bb{Z}}
\newcommand{\CC}{\bb{C}}
\newcommand{\HECKEQ}[1]{T_{\mathrm{q},{#1}}}
\newcommand{\HECKES}[1]{T_{\mathrm{s},{#1}}}
\newcommand{\FFp}{\overline{\bb{F}}_p}
\newcommand{\QQp}{\smash{\overline{\bb{Q}}_p}}
\newcommand{\QQptimes}{\smash{\overline{\bb{Q}}_p^\times}}
\newcommand{\FFptimes}{\smash{\overline{\bb{F}}_p^\times}}
\newcommand{\ZZp}{\smash{\overline{\bb{Z}}_p}}
\newcommand{\asterisk}{\star}
\newcommand{\defeq}{\vcentcolon=}
\newcommand{\SCR}{\mathscr}
\newcommand{\IE}{i.e.\ }

\newcommand{\DARKCIRC}{\bullet}

\newcommand{\HECKET}{\SCR T}

\newcommand{\MU}{\mu}

\newcommand{\SIGMA}{\Sigma}

\newcommand{\LAMBDA}{\lambda}

\newcommand{\LESS}[1]{{\uline{#1}}}
\newcommand{\MORE}[1]{{\overline{#1}}}
\newcommand{\EXTEND}[1]{{#1}^\circ}
\newcommand{\BIGO}[1]{\smash{\mathsf{O}\!\left({#1}\right)}}

\newcommand{\SIZEEE}{4.25pt}

\newcommand{\TEMPORARYA}{\xi}

\newcommand{\TEMPORARYD}{\vartheta}%

\newcommand{\TEMPORARYI}{\nu}
\newcommand{\TEMPORARYM}{\mathrm{S}}
\newcommand{\VERIFYIF}[1]{[{#1}]}
\newcommand{\absdet}{|\!\det\!|}

\newcommand{\PADICVALUATION}{v_p}

\newcommand{\FAMILY}[1]{({#1})}

\newcommand{\taga}{\m{\ast}1}
\newcommand{\tagb}{\m{\ast}2}
\newcommand{\tagc}{\m{\ast}3}
\newcommand{\tagd}{\m{\ast}4}
\newcommand{\tage}{\m{\ast}5}

\newcommand{\sfirst}{s_{1}}
\newcommand{\ssecond}{s_{2}}
\newcommand{\TILDE}{\widetilde}

\newcommand{\GEQSLANT}{\geqslant}
\newcommand{\LEQSLANT}{\leqslant}
\newcommand{\subm}{sub}
\newcommand{\UUU}{U}

\newcommand{\maxideal}{\mathfrak m}
\newcommand{\repnn}{\pi}

\newcommand{\SUBMODULE}{\mathrm{sub}}
\newcommand{\QUOTIENTI}{\mathrm{quot}}

\newcommand{\FILT}{\widehat}
\newcommand{\FORMX}{X}
\newcommand{\FORMY}{Y}
\newcommand{\DELTAA}{\Delta}
\newcommand{\TEMPORARYJ}{\epsilon}

\newcommand{\COMPACTIND}[2]{\ind^{#2}}
\newcommand{\sqbrackets}[4]{{#1}\bullet_{%
\,#3}{#4}}

\renewcommand{\HECKET}{T}

\newcommand{\polyr}{z}
\newcommand{\polyH}{\Phi}

\newcommand{\jj}{t}

\newtheorem{theorem}{Theorem}
\newtheorem{lemma}[theorem]{Lemma}
\newtheorem{proposition}[theorem]{Proposition}

\newtheorem{conjecture}{Conjecture}

\newenvironment{Smatrix}%
    {\tiny (\begin{smallmatrix}}%
    {\end{smallmatrix})}
\newenvironment{Lmatrix}%
    {\footnotesize (\begin{smallmatrix}}%
    {\end{smallmatrix})}
\newcommand{\BLACKSQUARE}{\hspace{\stretch{1}}\m{\mathbin{\vcenter{%
    \hbox{\rule{.9ex}{.9ex}}}}}}
\renewenvironment{proof}[1]%
    {\emph{{#1}}}%
    {\BLACKSQUARE}

\makeatletter
\newcommand{\xhookdoubleheadrightarrow}[2][]{%
  \lhook\joinrel
  \ext@arrow 0359\rightarrowfill@ {#1}{#2}%
  \mathrel{\mspace{-15mu}}\rightarrow}
\renewcommand{\xtwoheadrightarrow}[2][]{%
  \ext@arrow 0359\rightarrowfill@ {#1}{#2}%
  \mathrel{\mspace{-15mu}}\rightarrow}
\def\@settitle{\begin{center}%
  \baselineskip14\p@\relax\normalfont\LARGE
    \uppercasenonmath\@title
  \@title\end{center}%
}
\makeatother

\title%
[On the reductions of crystalline representations]%
{On the reductions of certain
two-dimensional crystalline representations}

\author{Bodan Arsovski}

\begin{document}
\maketitle

\begin{abstract}
The question of computing the reductions modulo~\m{p} of two-dimensional crystalline \m{p}-adic Galois representations has been studied extensively, and partial progress has been made for representations that have small weights, very small slopes, or very large slopes. It was conjectured by Breuil, Buzzard, and Emerton that these reductions are irreducible if they have even weight and non-integer slope. We prove some instances of this conjecture for slopes up to \m{\frac{p-1}{2}}.
\end{abstract}

\renewcommand{\HECKET}{T}

\section{Introduction}
\label{secInt}

Throughout this article we assume that \m{p} is an odd prime number. The main objects we study are two-dimensional crystalline representations of the absolute Galois group of \m{\QQ_p}. These are up to a twist parametrized by an integer \m{k\GEQSLANT2} and an element \m{a \in \ZZp} such that \m{\PADICVALUATION(a)>0}, and we denote by  \m{V_{k,a}} the representation corresponding to the parameters \m{(k,a)}\footnote{Thus all two-dimensional crystalline representations of \m{G_{\QQ_p}} are of the form \m{V_{k,a} \otimes \eta} for some character \m{\eta} of \m{G_{\QQ_p}} that is the product of an unramified character and a power of the cyclotomic character.}.
 In  section~\ref{secAss} we outline a construction of \m{V_{k,a}} and reference proofs of some of its basic properties. Thus in the remainder of this article we assume that  \m{k\in \ZZ} and \m{a \in \ZZp}  are parameters satisfying the properties  \m{k\GEQSLANT2} and \m{\PADICVALUATION(a)>0}. We also define the modulo \m{p} representation \m{\overline{V}_{k,a}} as the semi-simplification of the reduction modulo the maximal ideal \m{\maxideal} of \m{\ZZp} of a Galois stable \m{\ZZp}-lattice in \m{{V_{k,a}}}.
 The most interesting property of \m{V_{k,a}} is that
   the representation associated with a non-ordinary finite slope classical eigenform
  \[\textstyle f = \sum_{n=1}^\infty a_n q^n\]
 ---which is  normalized so that \m{a_1=1} and has
 weight \m{k\GEQSLANT 2}, level \m{\Gamma_0(N)} such that \m{(N,p)=1}, and  character \m{\chi}---is crystalline and equal to \m{V_{k,a_p\sqrt{\chi}(p)} \otimes \sqrt{\chi}}, where \m{\sqrt{\chi}} is an unramified character of \m{G_{\QQ_p}} whose square is \m{\chi}.

  \textbf{The conjecture.} The motivation behind the conjecture considered in this article comes from an observation made by Buzzard and Gouv\^ea in \cite{Buz05} and~\cite{Gou01} (based on computational data) that, if \m{p} belongs to a certain class of ``\m{\Gamma_0(N)}-regular''-primes, the slope \m{\PADICVALUATION(a_p)} of \m{f} appears to always be an integer. A definition of \m{\Gamma_0(N)}-regularity is given in~\cite{Buz05} (see also~\cite{BG16}), and its connection with Galois representations is that a prime \m{p>2} is \m{\Gamma_0(N)}-regular if and only if, for all weights \m{k\GEQSLANT 2} and all \m{f \in S_k(\Gamma_0(N))}, the modulo \m{p} representation associated with \m{f} is reducible. Note that, since the level of \m{f} is \m{\Gamma_0(N)}, its weight \m{k} must be even. This naturally leads to the general conjecture, made by Breuil, Buzzard, and Emerton, that if \m{k} is even and \m{\overline{V}_{k,a}} is reducible then \m{\PADICVALUATION(a)} is an integer. We consider this conjecture in the following slightly rephrased 
 form (see conjecture~4.1.1 in~\cite{BG16}).
\begin{conjecture}\label{MAINCONJ}
{If \m{k} is even and \m{\PADICVALUATION(a) \not\in \ZZ} then \m{\overline{V}_{k,a}} is irreducible.}
\end{conjecture}

\textbf{Known results.}
 The question of computing %
  \m{\overline{V}_{k,a}} has been studied extensively. A classification of \m{\overline{V}_{k,a}} for \m{k\LEQSLANT 2p+1} follows from the work of Berger and Breuil in~\cite{Ber10}, \cite{Bre03a}, and \cite{Bre03b}. A classification of \m{\overline{V}_{k,a}} when the slope \m{\PADICVALUATION(a)} is greater than \m{\lfloor \frac{k-2}{p-1}\rfloor} follows from the work of Berger, Li, and Zhu, in~\cite{BLZ04}. With some exceptions, \m{\overline{V}_{k,a}} has been computed when \m{\PADICVALUATION(a)<2} in the work of Buzzard, Gee, Ganguli, Ghate, Bhattacharya, Rozensztajn, and Rai in \cite{BG15}, \cite{BGR18}, \cite{BG09}, \cite{BG13}, \cite{GG15}, and~\cite{GG19a}.
     Finally, Rozensztajn's work \cite{Roz18} gives an algorithm that computes \m{\overline{V}_{k,a}} for given \m{p,k,a}, and her work~\cite{Roz1} gives an algorithm that finds the locus of all \m{a} such that \m{\overline{V}_{k,a}=\overline{\rho}} for given \m{p,k}, and a modulo \m{p} representation \m{\overline\rho}.

\textbf{The main theorem in this article.} 
We prove the following theorem.
\begin{theorem}\label{MAINTHM1}
{If \m{\PADICVALUATION(a) \not\in \ZZ} and \m{\TEMPORARYI=\lceil\PADICVALUATION(a)\rceil} is such that
 \[\textstyle k \not\equiv 3,4,\ldots,2\TEMPORARYI,2\TEMPORARYI+1 \bmod p-1,\] then \m{\overline{V}_{k,a}} is irreducible.
In particular, conjecture~\ref{MAINCONJ} is true when
\[\textstyle k \not\equiv 4,\ldots,2\TEMPORARYI \bmod p-1.\]}
\end{theorem}

Note that \m{\TEMPORARYI \in\{1,2,3,\ldots\}} and the  theorem is vacuous when  \m{\TEMPORARYI > \frac{p-1}{2}}. The proof of theorem~\ref{MAINTHM1} is based on the method  (developed by Breuil, Buzzard, and Gee) of using the local Langlands correspondence and its compatibility with reduction modulo \m{p} to compute the representations of \m{G_{\QQ_p}} by computing certain corresponding representations of \m{\GL_2(\QQ_p)}.
We give detailed outlines of the proof in sections~\ref{secLoc} and~\ref{secProof} after setting up the necessary framework.

\section*{Acknowledgements}\label{ackref}
I would like to thank professor Kevin Buzzard for his helpful remarks, suggestions, and support. This work was supported by Imperial College London and its President's PhD Scholarship.

\section{Crystalline representations}
\label{secAss}

Recall that we assume that \m{p} is an odd prime number and \m{k \in \ZZ} and \m{a \in \ZZp} are such that \m{k\GEQSLANT2} and \m{\PADICVALUATION(a)>0}. The  representation \m{{V}_{k,a}} is defined in subsection~3.1 of~\cite{Bre03b} as the two-dimensional representation of \m{G_{\QQ_p}} that is crystalline on \m{\QQ_p} and such that 
  \[\textstyle D_{cris}(V_{k,a}^\ast) = D_{k,a}\]
  for the weakly admissible filtered \m{\varphi}-module
\m{\textstyle D_{k,a} = \overline{\QQ}_p e_1 \oplus \overline{\QQ}_p e_2}
which has Hodge--Tate weights \m{(0,k-1)}, a filtration  \[\textstyle \textstyle\mathrm{Fil}^{i} ( D_{k,a}) = \left\{\begin{array}{rl} D_{k,a} & \text{if } i \LEQSLANT 0, \\[2pt] \QQp e_1 & \text{if } 0 < i < k,\\[2pt]
0 & \text{if } i \GEQSLANT k,  \end{array}\right.\] and a Frobenius map \m{\varphi} such that  
\[\textstyle \left\{\begin{array}{ll} \varphi (e_1)  = p^{k-1} e_2, &  \\[2pt] \varphi (e_2)  = -e_1+a e_2. &  \end{array}\right.\]
 All  two-dimensional crystalline representations of \m{G_{\QQ_p}} are of the form \m{V_{k,a}\otimes \eta} for some character \m{\eta} of \m{G_{\QQ_p}} that is the product of an unramified character and a power of the cyclotomic character (see proposition~3.1 in~\cite{Bre03b}).
 We define \m{\overline{V}_{k,a}} to be the semi-simplification of the reduction of a Galois stable \m{\ZZp}-lattice in \m{V_{k,a}} modulo the maximal ideal \m{\maxideal} of \m{\ZZp} (the
resulting representation is independent of the choice of lattice).
 
  In the remainder of this article we assume that \m{\PADICVALUATION(a) \not\in \ZZ} and we denote \m{\TEMPORARYI \defeq{} \lceil \PADICVALUATION(a) \rceil}, so that \m{\TEMPORARYI\in\{1,2,3,\ldots\}} and \m{\TEMPORARYI-1<\PADICVALUATION(a)<\TEMPORARYI}. Since theorem~\ref{MAINTHM1} is vacuous when \m{\PADICVALUATION(a)> \frac{p-1}{2}}, we assume that \m{\TEMPORARYI\in \{1,\ldots,\frac{p-1}{2}\}}. Moreover, since \m{\overline V_{k,a}} has been completely classified  when \m{k\LEQSLANT 2p+1}  (see for instance theorem~3.2.1 in~\cite{Ber10}), we assume that \m{a^2\ne 4p^{k-1}}  and \m{a\ne \pm(1+p^{-1})p^{k/2}}.
 
     \section{The \emph{p}-adic and modulo~\emph{p} local Langlands correspondences}
\label{secLoc}

 Let \m{B} be the Borel subgroup of  \m{G=\GL_2(\QQ_p)}  consisting of the upper triangular elements, 
 let \m{K=\GL_2(\ZZ_p) \subset G}, and let \m{Z} be the center of \m{G}. Let \m{W} be a finite-dimensional representation of the closed subgroup \m{KZ} of \m{G}. By a locally algebraic map \m{KZ\to W} we mean a map which on an open subgroup of \m{KZ} is the restriction of a rational map on (the algebraic group) \m{KZ}, and we say that \m{W} is locally algebraic if the map \m{h \in KZ \mapsto h w \in W} is locally algebraic for all \m{w \in W}. For a %
  locally algebraic finite-dimensional representation \m{W} of \m{KZ} we define the compact induction of \m{W} by
    \begin{align*}\textstyle \ind^G W \defeq %
     \textstyle \big\{\text{locally algebraic } G \to W \,\,\big|\,\, & f (h g) = h f(g) \text{ for all } h \in KZ \\ &\textstyle{} \text{\& }\Supp f\text{ is compact in }KZ\backslash G \}.\end{align*}
 Suppose that the representation \m{W} is over the field \m{\FF \in \{\QQp,\FFp\}}. For %
  elements \m{g\in G} and \m{w \in W} let \m{\sqbrackets{g}{H}{\FF}{w}} be the unique element of \m{\ind^G W} that is supported on \m{KZg^{-1}} and maps \m{g^{-1}} to \m{w}. Since \m{KZ\backslash G} is discrete, every element of \m{\ind^G W} can be written as a finite linear combination of functions of the type \m{\sqbrackets{g}{H}{\FF}{w}}. It is easy to check that \[\textstyle g_1(\sqbrackets{g_2}{H}{\FF}{(h w)}) = \sqbrackets{(g_1g_2h)}{H}{\FF}{w}.\] %
 We define \[\textstyle{}\sqbrackets{\xi}{H}{\FF}{w} = \xi(\sqbrackets{\mathrm{id}}{H}{\FF}{w})\] for \m{\xi \in \FF}.
Let \m{\mu_x} be the unramified character of the Weil group that sends the geometric Frobenius to \m{x}, and let \m{|\ | : \QQ_p^\times \to \QQ_p^\times \hookrightarrow\QQptimes} be the \m{p}-adic norm, so that \m{|x| = p^{-\PADICVALUATION(x)}}. We can view the symmetric power \m{{\Sym}^{k-2}(\overline{\QQ}_p^2)} as the \m{G}-module of homogeneous polynomials in \m{x} and \m{y} of total degree~\m{k-2} with coefficients in \m{\QQp}, with \m{G} acting by \[\textstyle \begin{Lmatrix}g_1 & g_2 \\ g_3 & g_4
\end{Lmatrix} \cdot v(x,y) = v(g_1x+g_3y,g_2x+g_4y)\] for \m{\begin{Smatrix}
g_1 & g_2 \\ g_3 & g_4
\end{Smatrix} \in G} and \m{v \in {{\Sym}^{k-2}(\overline{\QQ}_p^2)}}.
Let
\[\textstyle\TILDE\SIGMA_{k-2} = \textstyle{\uline{\Sym}^{k-2}}({\overline{\QQ}_p^2)} \defeq {\Sym^{k-2}}({\overline{\QQ}_p^2)} \otimes \absdet^{\frac{k-2}{2}}.\]
In particular, \m{\begin{Smatrix}
    p & 0 \\ 0 & p
    \end{Smatrix}} acts on \m{{\Sym^{k-2}}({\overline{\QQ}_p^2)}} as multiplication by \m{p^{k-2}} and it acts on \m{\TILDE\SIGMA_{k-2}} trivially. If \m{R \in \{\ZZp,\FFp\}}, we can view the \m{KZ}-module \m{\textstyle{}\uline{\Sym}^{k-2}(R^2)} %
     of homogeneous polynomials in \m{x} and \m{y} of total degree~\m{k-2} with coefficients in \m{R}, with \m{\begin{Smatrix}
    p & 0 \\ 0 & p
    \end{Smatrix}} acting trivially and \m{K} acting by \[\textstyle \begin{Lmatrix}g_1 & g_2 \\ g_3 & g_4
\end{Lmatrix} \cdot v(x,y) = v(g_1x+g_3y,g_2x+g_4y)\] for \m{\begin{Smatrix}
g_1 & g_2 \\ g_3 & g_4
\end{Smatrix} \in K} and \m{v \in {\uline{\Sym}^{k-2}}({R^2)}}. 
 Let \m{\SIGMA_{k-2}} be the reduction  of \m{\textstyle{}\uline{\Sym}^{k-2}({\overline{\ZZ}_p^2)}} modulo \m{\maxideal}.  Write \m{\lambda} for one of the roots of \m{X^2-aX+p^{k-1}}, so that the other root is \m{\lambda^{-1}p^{k-1}}. Let us consider the  \m{\QQp}-representation \[\textstyle \pi = \textstyle\Ind_B^G ({\mu_{\lambda p^{1-k}}}|\ |^{1/2} \times {\mu_{\lambda^{-1}}}|\ |^{-1/2}).\]
    Here the Borel subgroup \m{B} is seen as a parabolic subgroup of \m{G} and \m{\Ind_B^G} denotes parabolic induction (as opposed to \m{\ind^G} which denotes compact induction).
Let us fix embeddings \m{\overline{\QQ}\hookrightarrow \QQp} and \m{\overline{\QQ}\hookrightarrow \CC}. %
In section~2.1 of~\cite{Bre03b}, Breuil defines the Hecke operator \m{T \in \End_G(\COMPACTIND{KZ}{G} \TILDE\SIGMA_{k-2})} corresponding to the double coset of \m{\begin{Smatrix}
p & 0 \\ 0 & 1
\end{Smatrix}} and proves the explicit formula  
    \begin{align*} & \textstyle \HECKET ({\gamma}\sqbrackets{}{KZ}{\QQp}{v}) %
     \textstyle = {\sum_{\MU\in\FF_p} \gamma\begin{Lmatrix}
    p &[\MU] \\ 0 & 1
    \end{Lmatrix} }\sqbrackets{}{KZ}{\QQp}{\left(\begin{Lmatrix}
    1 & -[\MU] \\ 0 & p
    \end{Lmatrix} \cdot v\right)} + \gamma\begin{Lmatrix}
    1 & 0 \\ 0 & p
    \end{Lmatrix} \sqbrackets{}{KZ}{\QQp}{\left(\begin{Lmatrix}
    p & 0 \\ 0 & 1
    \end{Lmatrix} \cdot v\right)},\end{align*} where \m{[\TEMPORARYA]} is the Teichm\"uller lift of \m{\TEMPORARYA \in \FF_p} to \m{\ZZ_p}. 
 In section~3.2 of the same article he defines the \m{\GL_2(\QQ_p)}-representation
\[\textstyle \Pi_{k,a} = \COMPACTIND{KZ}{G} \TILDE\SIGMA_{k-2}/(T-a),\] proves that
\begin{align}\label{compat}\textstyle \Pi_{k,a} \cong {\Sym^{k-2}}({\overline{\QQ}_p^2)} \otimes \absdet^{-\frac12} \otimes \pi,\end{align} %
and defines
\[\textstyle \Theta_{k,a} = \Img\left(\COMPACTIND{KZ}{G}(\uline{\Sym}^{k-2}(\overline{\ZZ}_p^2)) \xrightarrow{\hspace{\SIZEEE} \hphantom{\cong} \hspace{\SIZEEE}}\Pi_{k,a}\right)\]
and \[\textstyle{}\overline\Theta_{k,a} = \Theta_{k,a} \otimes_{\ZZp} \FFp.\]

    For \m{h\in\ZZ} let us write \m{\sigma_h \defeq{} \uline{\Sym}^h(\FFp^2)}. For \m{t\in\{0,\ldots,p-1\}}, \m{\lambda \in\FFp}, and a character \m{\psi:\QQ_p^\times \to \overline\FF_p^\times}, let %
    \[\textstyle \repnn(t,\lambda,\psi) \defeq{} ({\COMPACTIND{KZ}{G} \sigma_t}/{(T_{\sigma}-\lambda)})\otimes\psi,\] where \m{T_\sigma\in \End_G(\COMPACTIND{KZ}{G} \sigma_t)} is the Hecke operator corresponding to the double coset of \m{\begin{Smatrix} p&0 \\ 0 &1\end{Smatrix}}. Let \m{\omega} be the modulo~\m{p} reduction of the cyclotomic character. Let \m{\ind(\omega_2^{t+1})} be the unique irreducible representation whose determinant is \m{{\omega^{t+1}}} and that is equal to \m{{\omega_2^{t+1}\oplus\omega_2^{p(t+1)}}} on inertia.  Let \m{\MORE{h}\in\{1,\ldots,p-1\}}  and \m{\LESS{h}\in\{0,\ldots,p-2\}} be the numbers in the corresponding sets that are congruent to \m{h} modulo \m{p-1}.
    
    In~\cite{Ber10}, Berger proves the following correspondence between \m{\overline V_{k,a}} and \m{\overline \Theta_{k,a}}, which was conjectured by Breuil (conjecture~3.3.5 in~\cite{Bre03b}).\footnote{In fact, what Berger proves is a more general correspondence for so-called ``trianguline'' representations.}
    
          \begin{theorem}\label{BERGERX1}
          There are \m{t\in\{0,\ldots,p-1\}} and \m{\psi:\QQ_p^\times \to \overline\FF_p^\times} such that either 
          \[\textstyle \overline\Theta_{k,a} \cong \repnn(t,0,\psi) %
          \]
          or
          \[\textstyle \overline\Theta_{k,a}^{\rm ss} \cong \left(\repnn(t,\lambda,\psi) \oplus\repnn(\LESS{p-3-t},\lambda^{-1},\omega^{t+1}\psi) \right)^{\rm ss}\]
          for some \m{\lambda \in\FFp}.
          In the former case we have
          \[\textstyle \overline V_{k,a} \cong \ind(\omega_2^{t+1})\otimes\psi,\]
          and in the latter case we have
          \[\textstyle \overline V_{k,a} \cong \left(\mu_\lambda \omega^{t+1} \oplus \mu_{\lambda^{-1}}\right)\otimes\psi.\]
   \end{theorem}
    
    Equation~\eqref{compat} implies that Berger's modulo \m{p} correspondence is compatible with the \m{p}-adic local Langlands correspondence, which associates with \m{{V}_{k,a}} a principal series representation. %
     As \m{\overline\Theta_{k,a}} has finite length as a \m{\FFp[G]}-module, we have \m{\overline\Theta_{k,a}^{\rm ss} \cong \overline\Pi_{k,a}^{\rm ss}} by lemma~3.3.4 in~\cite{Bre03b}, so we may replace \m{\Theta_{k,a}} with \m{\Pi_{k,a}} in theorem~\ref{BERGERX1}. %
     Theorem~\ref{BERGERX1} implies that theorem~\ref{MAINTHM1} can be rewritten in the following equivalent form. %
    
    \begin{theorem}\label{MAINTHM1aa}
{If \m{\PADICVALUATION(a) \not\in \ZZ} and \m{\TEMPORARYI=\lceil\PADICVALUATION(a)\rceil} is such that \[\textstyle k \not\equiv 3,4,\ldots,2\TEMPORARYI,2\TEMPORARYI+1 \bmod p-1,\] then \m{\overline{\Theta}_{k,a}} is irreducible.}
\end{theorem}
    
    Thus the goal of the sections that follow is to prove theorem~\ref{MAINTHM1aa}.

\section{Notation}
\label{secNot}

 For \m{\alpha \in \ZZp}, let \m{\BIGO{\alpha}} denote the  sub-\m{\ZZp}-module \[\alpha \COMPACTIND{KZ}{G}\left(\uline{\Sym}^{k-2}(\overline{\ZZ}_p^2)\right) \subseteq \COMPACTIND{KZ}{G}\left(\uline{\Sym}^{k-2}(\overline{\ZZ}_p^2)\right).\] We  abuse this notation and write \m{\BIGO{\alpha}} to represent a term \m{f \in \BIGO{\alpha}}. Let 
\[\textstyle \SCR I_{a} = \ker\left(\COMPACTIND{KZ}{G} \SIGMA_{k-2} \xtwoheadrightarrow{\hspace{\SIZEEE} \hphantom{\cong} \hspace{\SIZEEE}} \overline\Theta_{k,a}\right).\]
We use the shorthand ``\m{\Img(T-a)}'' for the image of the map \[\textstyle T-a \in \End_G(\COMPACTIND{KZ}{G} \TILDE\SIGMA_{k-2}).\] This image is a \m{G}-submodule of the \m{G}-module \m{\COMPACTIND{KZ}{G} \TILDE\SIGMA_{k-2}}. A crucial property of \m{\Img(T-a)} is that if an element of \m{\COMPACTIND{KZ}{G} \SIGMA_{k-2}} is the reduction modulo \m{\maxideal} of an integral element of \m{\Img(T-a)} then it is also in the ``kernel'' \m{\SCR I_a}.

If \m{V} is a \m{\FFp [KZ]}-module, let {\m{V(m)}} denote {the twist \m{V \otimes \det^m}}. Let {\m{I_h}} denote the %
 \m{\FFp [K Z]}-module of degree \m{h} homogeneous functions \m{{\FF_p^2} \to \FFp } that vanish at the origin, with \m{\begin{Smatrix}
    p & 0 \\ 0 & p
    \end{Smatrix}} acting trivially  and \m{K} acting by  \[\textstyle{}\left(\begin{Lmatrix}
g_1 & g_2 \\ g_3 & g_4
\end{Lmatrix}f\right)(x,y)=f(g_1x+g_3y,g_2x+g_4y)\]  for \m{\begin{Smatrix}
g_1 & g_2 \\ g_3 & g_4
\end{Smatrix} \in K} and \m{f(x,y)\in I_h}.
Note that if \m{h_1 \equiv h_2 \bmod {p-1}} then \m{I_{h_1} \cong I_{h_2}}.
Recall that, for \m{h\in\ZZ}, %
 \[\textstyle{}\sigma_h=\uline{\Sym}^h(\FFp^2).\] Then \m{\sigma_{\MORE{h}} \subset I_h}, since any element of  \m{\sigma_{\MORE{h}}} is also a function \m{{\FF_p^2 \to \FFp}} which is homogeneous  of degree \m{h} and vanishes at the origin, and the actions of \m{KZ} on \m{\sigma_{\MORE{h}}} and \m{I_h} match. %
  Due to lemma~3.2 in~\cite{AS86}, there is a map
    \begin{align}\label{EQUATION1}
    \textstyle f \in I_h \xmapsto{\hspace{\SIZEEE} \hphantom{\cong} \hspace{\SIZEEE}} \sum_{u,v} f(u,v) (v  X-u Y)^{\LESS{-h}} \in \sigma_{\LESS{-h}}(h),%
    \end{align} which gives an isomorphism  \[\textstyle{}I_h/\sigma_{\MORE{h}} \cong \sigma_{\LESS{-h}}(h),
    \]and therefore the only two factors of \m{I_h} are \m{\sigma_{\MORE{h}}} (``the submodule'') and \m{\sigma_{\LESS{-h}}(h)} (``the quotient''). If \m{\MORE{h} \ne p-1}, then \m{\sigma_{\LESS{-h}}(h)} is not a submodule of \m{I_h}. This is because the actions of \m{\begin{Smatrix}
    \lambda & 0 \\ 0 & \lambda
    \end{Smatrix}} on \m{\sigma_{\LESS{-h}}(h)} and \m{I_h} do not match. Hence in this case \m{\sigma_{\MORE{h}}} is the only submodule of \m{I_h}.
    If, on the other hand, \m{\MORE{h}=p-1}, then \m{\sigma_{\LESS{-h}}(h) = \sigma_0} is also a submodule of \m{I_h = I_0}: the submodule of those functions that are equal to a  constant everywhere except (possibly) at the origin.
    
    If \m{A\subseteq B\subseteq\COMPACTIND{KZ}{G}\SIGMA_r}, we say that an element \m{x \in \COMPACTIND{KZ}{G}\SIGMA_r} represents an element \m{y \in B/A} if \m{x \in B} and \m{x+A = y \in B/A}. Let \[\textstyle{}\theta\defeq{}xy^p-x^py,\qquad r\defeq{}k-2 \geqslant 0, \qquad s\defeq{}\MORE{r}, %
    \]  and for  a polynomial \m{f} with coefficients in \m{\ZZp} let \m{\overline{f}} denote its reduction modulo \m{\maxideal}.
     Suppose that \m{r \GEQSLANT \TEMPORARYI(p+1)}. 
 For \m{\alpha\in\{0,\ldots,\TEMPORARYI-1\}} let \[\textstyle N_\alpha = \overline{\theta}^{\alpha}\Sigma_{r-\alpha(p+1)}/\overline{\theta}^{\alpha+1}\Sigma_{r-(\alpha+1)(p+1)} \cong I_{r-2\alpha}(\alpha).\]
These are the first \m{\TEMPORARYI} subquotients of the filtration
\[\textstyle \Sigma_{r} \supset \overline{\theta}\Sigma_{r-p-1}\supset \cdots \supset \overline{\theta}^{\alpha}\Sigma_{r-\alpha(p+1)} \supset\cdots\]
This filtration corresponds to a filtration
\[\textstyle \COMPACTIND{KZ}{G}\Sigma_{r} \supset \COMPACTIND{KZ}{G}(\overline{\theta}\Sigma_{r-p-1})\supset \cdots \supset \COMPACTIND{KZ}{G}(\overline{\theta}^{\alpha}\Sigma_{r-\alpha(p+1)}) \supset\cdots\]
of \m{\COMPACTIND{KZ}{G}\SIGMA_r}, which has corresponding subquotients \m{\FAMILY{\FILT N_\alpha}_{0\LEQSLANT\alpha<\TEMPORARYI}}. Since \[\textstyle \overline\Theta_{k,a}\cong \COMPACTIND{KZ}{G}\SIGMA_r/\SCR I_a,\] there is also a filtration
\[\textstyle \overline\Theta_{k,a} = \overline\Theta_0 \supset \overline\Theta_1 \supset \cdots \supset \overline \Theta_\alpha \supset \cdots\]
whose first \m{\TEMPORARYI} subquotients are quotients of \m{\FAMILY{\FILT N_\alpha}_{0\LEQSLANT\alpha<\TEMPORARYI}}, \IE if \m{\alpha\in\{0,\ldots,\TEMPORARYI-1\}} then there is a surjection \[\textstyle\textstyle \FILT N_\alpha \xtwoheadrightarrow{\hspace{\SIZEEE} \hphantom{\cong} \hspace{\SIZEEE}} \overline\Theta_\alpha/\overline\Theta_{\alpha+1},\] and the kernel of this surjection consists of those elements of the subquotient \m{\FILT N_\alpha} of \m{\COMPACTIND{KZ}{G}\SIGMA_r} that are represented by an element of \m{\SCR I_a \subset \COMPACTIND{KZ}{G}\SIGMA_r}. 
In the proof of theorem~\ref{MAINTHM1} we compute \m{\overline\Theta_{k,a}} by %
 finding elements in the kernels of these surjections and consequently obtaining data about the quotients \m{
 \FAMILY{\overline\Theta_\alpha/\overline\Theta_{\alpha+1}}_{0\LEQSLANT\alpha<\TEMPORARYI}} of \m{\FAMILY{\FILT N_\alpha}_{0\LEQSLANT\alpha<\TEMPORARYI}}.
 Therefore, the strategy of the proof is to find elements
 of \m{\SCR I_a} that represent non-trivial elements of the subquotients \m{\FAMILY{\FILT N_\alpha}_{0\LEQSLANT\alpha<\TEMPORARYI}}.

For a property \m{P} let us define \m{[P]=1} if \m{P} is true and \m{[P]=0} if \m{P} is false.
Lemmas~4.1 and~4.3 and remark~4.4 in~\cite{BG09} imply that the ideal \m{\SCR I_a} contains \m{1\sqbrackets{}{KZ}{\QQp}{x^jy^{r-j}}}
 for all \m{\textstyle 0\LEQSLANT j < \TEMPORARYI} and  \m{1\sqbrackets{}{KZ}{\QQp}{\overline{\theta}^{\TEMPORARYI}h}} for all  \m{h\in\Sigma_{r-\TEMPORARYI(p+1)}}, and thus \m{\overline{\Theta}_{k,a}} is a subquotient of \[\textstyle \COMPACTIND{KZ}{G}(\Sigma_r/\langle y^r,\ldots,x^{\TEMPORARYI-1}y^{r-\TEMPORARYI+1},\overline{\theta}^{\TEMPORARYI}\Sigma_{r-\TEMPORARYI(p+1)}\rangle),\]
 a module which has a series whose factors are subquotients of \m{\FILT N_0,\ldots,\FILT N_{\TEMPORARYI-1}}.
In particular, \m{\overline{\Theta}_\TEMPORARYI=0}. For \m{\alpha\in\{0,\ldots,\TEMPORARYI-1\}} let us denote
\begin{align*}
    \SUBMODULE(\alpha) & \textstyle = \sigma_{\MORE{r-2\alpha}}(\alpha) \subset N_\alpha,\\
    \QUOTIENTI(\alpha) & = N_\alpha/\sigma_{\MORE{r-2\alpha}}(\alpha) \cong \sigma_{\LESS{2\alpha-r}}(r-\alpha).
\end{align*}
We denote by \m{T,\HECKEQ{\alpha},\HECKES{\alpha}} the Hecke operators corresponding to the double coset of \m{\begin{Smatrix} p&0 \\ 0 &1\end{Smatrix}} on the modules \m{\COMPACTIND{KZ}{G}\TILDE \SIGMA_r, \COMPACTIND{KZ}{G}\QUOTIENTI(\alpha),\COMPACTIND{KZ}{G}\SUBMODULE(\alpha)}, respectively.

\section{Local constancy in the parameter \emph{k}}
\label{secAss2}

The main theorem is true for \m{p=3} (this follows from theorem~1.6 in~\cite{BG09}), so we may assume that \m{p\GEQSLANT 5}. In particular, since we assume that \m{\PADICVALUATION(a) <\frac{p-1}{2}} and \m{k\GEQSLANT 2p+2}, and since \m{(4p+2)(p^2-3p+1) \GEQSLANT 3(p-1)^3} for \m{p\GEQSLANT5}, we can conclude that
\[\textstyle{}k > 3 \PADICVALUATION(a) + \frac{(k-1)p}{(p-1)^2} + 1.\] Therefore theorem~B in~\cite{Ber12} implies that there is a constant \m{m = m(k,a)} such that \[\textstyle \textstyle \overline{V}_{k,a} \cong \overline{V}_{k^{\prime},a}\] whenever \m{k^{\prime}\GEQSLANT k} and \m{k^{\prime}\equiv k \bmod (p-1)p^m}. By using this isomorphism we can conclude that if the main theorem is true for \m{k>p^{100}} then it is true for all \m{k}.  Therefore in the remainder of this article we assume that \m{k>p^{100}}.

\section{Combinatorial definitions}
\label{secAss3}

For a formal variable \m{\FORMX} and \m{n\in \ZZ} let us define \[\textstyle\binom{\FORMX}{n} \defeq \frac{\FORMX_n}{n!} \defeq \left\{\begin{array}{rl}\frac{1}{n!} \prod_{j=0}^{n-1}(X-j) \in \QQ_p[X] & \text{if } n \GEQSLANT 0, \\[2pt] 0\in \QQ_p[X] & \text{if } n<0.  \end{array}\right.
\] Therefore \[\textstyle\FORMX_n=\prod_{j=0}^{n-1}(X-j) \in \ZZ_p[X]\] denotes the falling factorial when \m{n\GEQSLANT 0}.
 For \m{m \in \QQ_p} and \m{n\in \ZZ} let us define \m{\binom{m}{n} \in \QQ_p} as the evaluation of \m{\binom{\FORMX}{n}} at \m{\FORMX = m \in \QQ_p}. In particular, binomial coefficients with negative denominators are always zero. Let
\[\textstyle  \binom{\FORMX}{n}^\partial \defeq \frac{\partial}{\partial \FORMX}\binom{\FORMX}{n} = \binom{X}{n}\sum_{j=0}^{n-1} \frac{1}{X-j} \in \QQ_p[X].\]
For \m{m \in \QQ_p} and \m{n\in \ZZ} let us define \m{\binom{m}{n}^\partial \in \QQ_p} as the evaluation of \m{\binom{\FORMX}{n}^\partial} at \m{\FORMX = m \in \QQ_p}.
 If \m{\vartheta(X) \in \ZZ_p[X]} then there is the Taylor series expansion
\[\vartheta(b+\epsilon) = \vartheta(b) + \epsilon \vartheta^{\prime}(b) + \BIGO{\epsilon^2}\]
for \m{b,\epsilon\in \ZZ_p}. Indeed, this follows from the facts that
\[\textstyle (b+\epsilon)^{m} = b^m + \epsilon mb^{m-1} + \epsilon^2 \sum_{j=2}^m \binom{m}{j}\epsilon^{j-2}b^{m-j}\]
and
\[\textstyle \sum_{j=2}^m \binom{m}{j}\epsilon^{j-2}b^{m-j} \in \ZZ_p\]
for \m{b,\epsilon\in \ZZ_p} and \m{m \in \ZZ_{\GEQSLANT 0}}. Since \m{p^{\PADICVALUATION(n!)}\binom{\FORMX}{n} \in \ZZ_p[\FORMX]} for \m{n\in\ZZ_{\GEQSLANT 0}},
\begin{align}\label{eqTaylor}\textstyle \binom{b+\epsilon}{n} = \binom{b}{n} + \epsilon \binom{b}{n}^\partial + \BIGO{\epsilon^2p^{-\PADICVALUATION(n!)}}\end{align}
for \m{b,\epsilon\in \ZZ_p} and \m{n \in \ZZ_{\GEQSLANT0}}.
For \m{n,k\in\ZZ_{\GEQSLANT 0}} let us define the Stirling number of the first kind \m{\sfirst(n,k)} as %
 the coefficient of \m{\FORMX^k} in \m{\FORMX_n=\FORMX\cdots(\FORMX-n+1)}. Therefore, %
 \[\textstyle \left(\frac{\sfirst(i,j)}{i!}\right)_{0\LEQSLANT i,j\LEQSLANT m} \cdot \left(1,\ldots,\FORMX^m\right)^T = \left(\binom{\FORMX}{0},\ldots,\binom{\FORMX}{m}\right)^T.\]
Let us  also define the Stirling number of the second kind \m{\ssecond(n,k)} as the coefficient of \m{\FORMX_k} in \m{\FORMX^n}, in the sense that
\[\textstyle \FORMX^n = \sum_{k=0}^n \ssecond(n,k) \FORMX_k = \sum_{k=0}^n \ssecond(n,k) \prod_{j=0}^{k-1}(\FORMX-j).\]
In particular, \m{\left(\ssecond(i,j)\right)_{0\LEQSLANT i,j\LEQSLANT m}}
is the inverse of \m{\left(\sfirst(i,j)\right)_{0\LEQSLANT i,j\LEQSLANT m}} and
\[\textstyle \left(j! \ssecond(i,j)\right)_{0\LEQSLANT i,j\LEQSLANT m}  \cdot \left(\binom{\FORMX}{0},\ldots,\binom{\FORMX}{m}\right)^T = \left(1,\ldots,\FORMX^m\right)^T.\]
For %
 a family \m{\FAMILY{D_i}_{i\in\ZZ}} of elements of \m{\ZZ_p} supported on a finite set of indices, and for \m{w\in\ZZ}, let us define \[\textstyle \TEMPORARYD_{w}(D_\DARKCIRC) = \TEMPORARYD_{w}(\FAMILY{D_i}_{i\in \ZZ}) = \sum_{i\in\ZZ} D_i \binom{i(p-1)}{w}.\]
 Let us also define
    \[\textstyle\TEMPORARYM_{u,n} = %
    \sum_{i\in\ZZ} \binom{u}{i(p-1)+n}\]
for \m{u\in\ZZ_{\GEQSLANT 0}} and \m{n \in \ZZ}, and let \m{\TEMPORARYM_u = \TEMPORARYM_{u,0}}.
We use the convention that when the range of summation is
not specified, it is assumed to be over all of
\m{\ZZ}, so that ``\m{\sum_{m_1,\ldots,m_k}}'' means ``\m{\sum_{m_1 \in \ZZ} \cdots \sum_{m_k \in \ZZ}}''.

\section{Table of assumptions and definitions}
\label{secList}
In this section we collate the assumptions, definitions, and notation we have introduced in  a convenient table.

\newpage

\newcommand{\NEWLINEB}{\\[7pt]}
\null
\vspace{\stretch{1}}

\begin{tabular}{r|l}
    \m{p} & \m{p>2} is the prime. \NEWLINEB
    \m{G, K, Z} & \m{G=\GL_2(\QQ_p)}, \m{K=\GL_2(\ZZ_p)}, and \m{Z} is the center of \m{G}. \NEWLINEB
         \m{\MORE{h}} & the number in \m{\{1,\ldots,p-1\}}  congruent to \m{h} modulo \m{p-1}. \NEWLINEB
\m{\LESS{h}} & the number in \m{\{0,\ldots,p-2\}}  congruent to \m{h} modulo \m{p-1}. \NEWLINEB
    \m{k,r,s} & \m{k} is the weight and we assume \m{k>p^{100}}, \m{r = k-2}, and \m{s=\MORE{r}}. %
     \NEWLINEB
        \m{a,\TEMPORARYI} & \m{a} is the eigenvalue and \m{\TEMPORARYI = \lceil \PADICVALUATION(a)\rceil \in \{1,\ldots,\frac{p-1}{2}\}}. \NEWLINEB
     \m{\sqbrackets{g}{H}{\FF}{w}} & is in \m{\COMPACTIND{KZ}{G} W}, is supported on \m{KZg^{-1}}, and maps \m{g^{-1}\mapsto w}. \NEWLINEB
     \m{I_t} & \m{\FFp[KZ]}-module of degree \m{t} maps \m{\FF_p^2\to\FFp} that vanish at \m{(0,0)}. \NEWLINEB
          \m{\sigma_t} & \m{\uline{\Sym}^t(\FFp^2)}. \NEWLINEB
              \m{V(m)} & the twist \m{V \otimes \det^m}.
              \NEWLINEB
                  \m{\textstyle\TILDE\SIGMA_{k-2}} &  \m{\uline{\Sym}^{k-2}(\overline{\QQ}_p^2) = {{\Sym}^{k-2}}({\overline{\QQ}_p^2)}\otimes \absdet^{\frac{k-2}{2}}}. \NEWLINEB
        \m{\SIGMA_{k-2}} & the reduction  of \m{\textstyle \uline{\Sym}^{k-2}(\overline{\ZZ}_p^2)} modulo \m{\maxideal}. \NEWLINEB
        \m{\HECKET} & Hecke operator  for the double coset of \m{\begin{Smatrix}
p & 0 \\ 0 & 1
\end{Smatrix}} on \m{\COMPACTIND{KZ}{G}\TILDE \SIGMA_{k-2}}. \NEWLINEB
\m{\HECKEQ{\alpha}} & Hecke operator  for the double coset of \m{\begin{Smatrix}
p & 0 \\ 0 & 1
\end{Smatrix}} on \m{ \COMPACTIND{KZ}{G}\QUOTIENTI(\alpha)}. \NEWLINEB
\m{\HECKES{\alpha}} & Hecke operator  for the double coset of \m{\begin{Smatrix}
p & 0 \\ 0 & 1
\end{Smatrix}} on \m{\COMPACTIND{KZ}{G}\SUBMODULE(\alpha)}. \NEWLINEB
\m{\BIGO{\alpha}} & sub-\m{\ZZp}-module of multiples of \m{\alpha}; also used for \m{f \in \BIGO{\alpha}}. \NEWLINEB
\m{\theta} & \m{\theta = xy^p - x^py}. \NEWLINEB
\m{N_\alpha} & \m{N_\alpha  = \overline{\theta}^{\alpha}\Sigma_{r-\alpha(p+1)}/\overline{\theta}^{\alpha+1}\Sigma_{r-(\alpha+1)(p+1)} \cong I_{r-2\alpha}(\alpha)}.
\NEWLINEB
\m{\binom{X}{n}} & \m{\binom{X}{n} = \frac{X\cdots(X-n+1)}{n!} \in \QQ_p[X]} for \m{n\GEQSLANT 0} and \m{\binom{X}{n} = 0} otherwise. \NEWLINEB
\m{X_n} & \m{X_n = n!\binom{X}{n}\in \QQ_p[X]} for \m{n\GEQSLANT 0} and \m{{X_n} = 0} otherwise. \NEWLINEB
\m{\binom{\FORMX}{n}^\partial} & \m{\binom{\FORMX}{n}^\partial= \frac{\partial}{\partial \FORMX}\binom{\FORMX}{n} = \binom{X}{n}\sum_{j=0}^{n-1} \frac{1}{X-j} \in \QQ_p[X]}. %
 \NEWLINEB
\m{\TEMPORARYM_{u,n}} & \m{\TEMPORARYM_{u,n} = %
    \sum_i \binom{u}{i(p-1)+n}} 
for \m{u\GEQSLANT 0} and \m{n \in \ZZ}, and \m{\TEMPORARYM_u = \TEMPORARYM_{u,0}}.\NEWLINEB
\m{\sfirst(n,k)} &
 the coefficient of \m{\FORMX^k} in \m{\FORMX_n=\FORMX\cdots(\FORMX-n+1)}. \NEWLINEB  \m{\ssecond(n,k)} & the coefficient of \m{\FORMX_k} in \m{\FORMX^n}. \NEWLINEB
 \m{\TEMPORARYD_{w}(D_\DARKCIRC)}  &
 \m{\TEMPORARYD_{w}(D_\DARKCIRC) = \TEMPORARYD_{w}(\FAMILY{D_i}_{i\in \ZZ}) = \sum_{i\in\ZZ} D_i \binom{i(p-1)}{w}}.  \NEWLINEB
\m{\SCR I_a} & the kernel of the quotient map 
\m{\textstyle \COMPACTIND{KZ}{G} \SIGMA_{k-2} \xtwoheadrightarrow{} \overline\Theta_{k,a}}.  \NEWLINEB
\m{[P]} & \m{[P]=1} if \m{P} is true, \m{[P]=0} if \m{P} is false.  \NEWLINEB
\end{tabular}

\vspace{\stretch{1}}
\null

\newpage

\section{Combinatorial identities}\label{section2}

This section comprises nine lemmas about combinatorial identities involving binomial sums.  The proofs of these nine lemmas are all standard. Therefore the reader may wish to skip this section on initial reading and refer back to it as required.

Recall that we use the convention that \m{\binom{X}{n} \in \QQ_p[X]} is a polynomial of degree \m{n} if \m{n\GEQSLANT 0}, and  \m{\binom{X}{n}=0\in \QQ_p[X]} if \m{n<0}. Moreover, \m{\binom{m}{n} \in \QQ_p} is the evaluation of  \m{\binom{X}{n}} at \m{X=m \in \QQ_p}. In particular, we use the convention that binomial coefficients with negative denominators are always zero, which is in contrast with some of the literature. We also use the convention that when the range of summation is not specified, it is assumed to be over all of \m{\ZZ}, \IE we define
\[\textstyle \sum_{m_1,\ldots,m_k} F(m_1,\ldots,m_k)  =  \sum_{m_1 \in \ZZ} \cdots \sum_{m_k \in \ZZ} F(m_1,\ldots,m_k)  .\]
This never gives rise to convergence issues as  \m{F(m_1,\ldots,m_k)} is  supported on a finite subset of \m{\ZZ^k} whenever such a sum appears in this article.

\textbf{A formula for \m{\boldsymbol{\binom{m}{n}^\partial}} when \m{\boldsymbol{0 \LEQSLANT m < n}}.} In general, the derivative of the binomial coefficient is difficult to compute as its formula involves the harmonic numbers. There is a closed formula in one special case, however, given by the following lemma.

\begin{lemma}\label{partHARMONIC}
Suppose that \m{m,n\in\ZZ_{\GEQSLANT0}} are such that \m{m<n}. Then
\[\textstyle{} \binom{m}{n}^\partial = \frac{(-1)^{m+n+1}}{n \binom{n-1}{m}}= \frac{(-1)^{m+n+1}}{(m+1) \binom{n}{m+1}}.\]
\end{lemma}

\begin{proof}{Proof. }
For \m{j\in\{0,\ldots,n-1\}}, let \m{f_{n,j}\in \ZZ_p[X]} denote the polynomial \[\textstyle{} f_{n,j}(X) = \prod_{i \in \{0,\ldots,n-1\}\backslash \{j\}} (X-i),\]
so that
\[\textstyle{}\binom{m}{n}^\partial = \frac1{n!}\sum_{j=0}^{n-1}f_{n,j}(m).\]
If \m{j\ne m}, then \m{X-m\mid f_{n,j}} and therefore \m{f_{n,j}(m)=0}, so in fact
\[\textstyle{}\binom{m}{n}^\partial = \frac{f_{n,m}(m)}{n!} %
= \frac{(-1)^{n-m-1} m! (n-m-1)!}{n!} = \frac{(-1)^{m+n+1}}{n \binom{n-1}{m}}= \frac{(-1)^{m+n+1}}{(m+1) \binom{n}{m+1}}.\]
\end{proof}

\textbf{Combinatorial identities involving binomial sums.}
In this subsection we prove four lemmas about generic identities between binomial sums.

\begin{lemma}[Properties of \m{\TEMPORARYM_{u,n}}]\label{part1X5}
Suppose throughout this lemma that 
\begin{align*}
    n,t \in \ZZ, \qquad
    b,d,k,l,w  \in \ZZ_{\GEQSLANT 0}, \qquad
    m,u,v  \in \ZZ_{\GEQSLANT 1}.
\end{align*}
\begin{enumerate}
    \item If \m{u \equiv v \bmod (p-1)p^{m-1}} then
\begin{align}
\label{ceqn1}\tag{\emph{c-a}}
\textstyle \TEMPORARYM_{u,n} \equiv \TEMPORARYM_{v,n} \bmod p^m.
\end{align}
\item Suppose that \m{u = t_u(p-1)+s_u} with \m{s_u = \MORE{u}}, so that \m{s_u \in \{1,\ldots,p-1\}} and \m{t_u \in \ZZ_{\GEQSLANT 0}}. Then \begin{align}
\label{ceqn2}\tag{\emph{c-b}}
\textstyle \TEMPORARYM_u = 1 + \VERIFYIF{u\equiv_{p-1} 0} + \frac {t_u}{s_u}p + \BIGO{t_up^2}.
\end{align}
\item  If \m{n\LEQSLANT 0} then
\begin{align}
\label{ceqn3}\tag{\emph{c-c}}
\textstyle \TEMPORARYM_{u,n} = %
\sum_{i=0}^{-n} (-1)^i \binom{-n}{i} \TEMPORARYM_{u-n-i}.%
\end{align} \item If \m{n\GEQSLANT 0} then
\begin{align}
\label{ceqn4}\tag{\emph{c-d}}\textstyle \TEMPORARYM_{u,n} \equiv (1+\VERIFYIF{u\equiv_{p-1}n\equiv_{p-1} 0})\binom{\MORE{u}}{\LESS{n}} %
 \bmod p.
\end{align}
\item If \m{\FORMX} is a formal variable  then
\begin{align}\label{ceqn8}\tag{\emph{c-e}}
   \textstyle  \binom{\FORMX}{t+l} \binom{t}{w} = \sum_v \binom{-l}{w-v} \binom{\FORMX}{v} \binom{\FORMX-v}{t+l-v}.
\end{align}
Consequently, if \m{b+l\GEQSLANT d+w} then
\begin{align}
\label{ceqn6}\nonumber& \textstyle \sum_i %
 \binom{b-d+l}{i(p-1)+l}\binom{i(p-1)}{w} %
  = \sum_v  \binom{-l}{w-v} \binom{b-d+l}{v} \TEMPORARYM_{b-d+l-v,l-v}.\tag{\emph{c-f}}\end{align}
\end{enumerate}
\end{lemma}

\begin{proof}{Proof. }
\begin{enumerate}
    \item We can rewrite \m{ \TEMPORARYM_{u,n}} by using the equation
\[\textstyle \TEMPORARYM_{u,n} = %
\frac{1}{p-1}\sum_{0\ne \MU\in\FF_p} [\MU]^{-n} (1 + [\MU])^u.\]
This is because 
\[\textstyle \sum_{0\ne \MU\in\FF_p} [\MU]^{\lambda} = \left\{\begin{array}{rl}
    p-1 & \text{if } p-1 \mid \lambda,  \\
    0 & \text{otherwise.}
\end{array}\right.\]
Since \m{1+[\MU] \in \{0\} \cup \ZZ_p^\times} for all \m{\MU \in \FF_p} (as \m{p} is odd),
we have
\[\textstyle{}(1+[\MU])^{u} \equiv (1+[\MU])^v \bmod p^m \] for all \m{\MU\in\FF_p} (because  \m{\textstyle{}|(\ZZ/p^m\ZZ)^\times|} divides \m{u-v} and \m{u,v>0}). %
 Therefore \m{\TEMPORARYM_{u,n} \equiv \TEMPORARYM_{v,n} \bmod p^m}.
\item For \m{\MU \in \FF_p\backslash\{-1,0\}} let us define
\[\textstyle x_\MU = \frac{(1+[\MU])^{p-1}-1}{p} \in \ZZ_p.\]
Then %
 \m{\TEMPORARYM_u} is equal to
\begin{align*}
    \textstyle  & \textstyle 1+ \VERIFYIF{s_u=p-1} + \frac{1}{p-1}\sum_{\MU\in\FF_p\backslash\{-1,0\}} (1+[\MU])^{s_u} \sum_{j=1}^{t_u} \binom{t_u}{j} p^j x_\MU^j .
\end{align*} It is easy to show that \[\textstyle{}\binom{t_u}{j} p^j = \BIGO{t_up^2}\] for \m{j>1}
by writing
\[\textstyle \PADICVALUATION\left(\binom{t_u}{j}p^j\right) \GEQSLANT \PADICVALUATION(t_u) + j - \sum_{i\GEQSLANT 1} \lfloor jp^{-i} \rfloor\] and verifying that the right side is at least \m{ \PADICVALUATION(t_up^2)} for \m{j\in\{2,3\}} and at least \m{\textstyle{}\PADICVALUATION(t_u) + \frac{p-2}{p-1}j \GEQSLANT \PADICVALUATION(t_up^2)}
for \m{j\GEQSLANT 4}.
This implies that if \[\textstyle A_{s_u} = \frac{1}{p-1}\sum_{\MU\in\FF_p\backslash\{-1,0\}} (1+[\MU])^{s_u} x_\MU\]
then
\[\textstyle \TEMPORARYM_u = 1+\VERIFYIF{s_u=p-1} + A_{s_u} t_u p + \BIGO{t_up^2}.\] 
Moreover, the constant \m{A_{s_u}} is independent of \m{t_u}. For \m{t_u=1} we have
\begin{align*}\textstyle \TEMPORARYM_{s_u+p-1} & \textstyle{} = 1 + \VERIFYIF{s_u=p-1} + \binom{s_u+p-1}{p-1}
\\ & \textstyle{} = 1 + \VERIFYIF{s_u=p-1} + \binom{s_u-1}{p-1} + p\binom{s_u-1}{p-1}^\partial + \BIGO{p^2}
\\ & \textstyle{} = 1 + \VERIFYIF{s_u=p-1} + p\binom{s_u-1}{p-1}^\partial + \BIGO{p^2}.\end{align*}
The second equality follows from %
 equation~\eqref{eqTaylor}, and the third equality follows from the fact that \m{s_u-1 \in \{0,\ldots,p-2\}}. Thus, due to lemma~\ref{partHARMONIC},
\[\textstyle{} A_{s_u} = \binom{s_u-1}{p-1}^\partial + \BIGO{p} %
 = \frac {(-1)^{s_u}}{s_u\binom{p-1}{s_u}} + \BIGO{p} = \frac 1{s_u} + \BIGO{p}.\]
\item  This follows from repeated application of Pascal's triangle equation \[\textstyle \TEMPORARYM_{u,v} = \TEMPORARYM_{u-1,v}+\TEMPORARYM_{u-1,v-1},\] which is true for \m{u,v\GEQSLANT1}. We omit the full details.
\item  This follows from the congruence
\[\textstyle %
\TEMPORARYM_{u,n} \equiv -\sum_{0\ne \MU\in\FF_p} [\MU]^{-n}(1+[\MU])^u \equiv -\sum_{0\ne \MU\in\FF_p} [\MU]^{-\LESS{n}}(1+[\MU])^{\MORE{u}}\bmod p.\] 
\item 
We can use
\begin{align*}%
\textstyle \binom{\FORMX}{v} \binom{\FORMX-v}{t+l-v} & \textstyle %
 = \binom{\FORMX}{t+l} \binom{t+l}{v} \end{align*} to rewrite equation~\eqref{ceqn8}  as
\begin{align*}\textstyle \binom{\FORMX}{t+l}  \binom{t}{w} %
 &\textstyle{} = \binom{\FORMX}{t+l}  \sum_v  \binom{-l}{w-v}  \binom{t+l}{v},\end{align*}
 which follows from Vandermonde's convolution formula\footnote{Vandermonde's convolution formula is the fact that, for all \m{A,B,C\in\ZZ}, \[\textstyle{} \sum_v \binom{A}{v} \binom{B}{C-v} = \binom{A+B}{C}.\] %
}. In particular,
if we apply equation~\eqref{ceqn8} to \m{\FORMX=b-d+l} and \m{t=i(p-1)}, we get
\[\textstyle  
 \binom{b-d+l}{i(p-1)+l}\binom{i(p-1)}{w} = \sum_v  \binom{-l}{w-v} \binom{b-d+l}{v} \binom{b-d+l-v}{i(p-1)+l-v}.\]
Then we get equation~\eqref{ceqn6} by summing over all \m{i\in \ZZ}.
\end{enumerate}
\end{proof}

\begin{lemma}\label{part1X5prime}
Suppose throughout this lemma that 
\begin{align*}
    y \in \ZZ, \qquad
    b,d,l,w  \in \ZZ_{\GEQSLANT 0}, \qquad
    u  \in \ZZ_{\GEQSLANT 1}.
\end{align*}
\begin{enumerate}
\item We have
\begin{align}
\label{ceqn13}\tag{\emph{c-g}}
\textstyle \textstyle \sum_j (-1)^{b+j} \binom{j}{y}\binom{b}{j} & \textstyle{} %
 =\VERIFYIF{y=b}.
\end{align}
\item We have \begin{align}
\label{ceqn9}\tag{\emph{c-h}}
\textstyle \sum_{j} (-1)^j \binom{y}{j} \binom{y+l-j}{w-j} = \binom{l}{w}.
\end{align}
\item We have 
\begin{align}
\label{ceqn10}\tag{\emph{c-i}}
\textstyle \sum_{j} \binom{u-1}{j-1} \binom{-l}{j-w} = (-1)^{u-w} \binom{l-w}{u-w}.
\end{align}
\item We have
\begin{align}
\label{ceqn11}\tag{\emph{c-j}}
\textstyle \sum_{j} (-1)^{j}  \binom{l}{j-b}\binom{j}{w} = (-1)^{b+l}\binom{b}{w-l}.
\end{align}
\item If \m{u \GEQSLANT (b+l)d} and \m{l\GEQSLANT w} then
\begin{align}
\label{ceqn5}\tag{\emph{c-k}}\textstyle \sum_j (-1)^{j-b} \binom{l}{j-b}\binom{u-dj}{w} = \VERIFYIF{w=l}d^l.
\end{align}
\end{enumerate}
\end{lemma}

\begin{proof}{Proof. }
\begin{enumerate}
\item We have
\begin{align*}
\textstyle \textstyle \sum_j (-1)^{b+j} \binom{j}{y}\binom{b}{j} & \textstyle{} = \sum_j (-1)^{b+j} \binom{j}{j-y}\binom{b}{b-j} \\ &\textstyle{}
 =  \sum_j (-1)^{b-y} \binom{-y-1}{j-y}\binom{b}{b-j}
 \\ & \textstyle{}
 = (-1)^{b-y} \binom{b-y-1}{b-y} 
 =\VERIFYIF{y=b}.
\end{align*} The third equality follows from Vandermonde's convolution formula.
\item  Since \[\textstyle{}\textstyle (-1)^j \binom{y+l-j}{w-j} = (-1)^{w} \binom{w-y-l-1}{w-j} \text{ and } (-1)^w \binom{w-l-1}{w} = \binom{l}{w},\] this follows from Vandermonde's convolution formula.
\item 
Since \[\textstyle \binom{u-1}{j-1} = \binom{u-1}{u-j}  \text{ and }  \binom{u-l-1}{u-w} = (-1)^{u-w}\binom{l-w}{u-w},\] 
this follows from Vandermonde's convolution formula.
\item We have
\begin{align*}
    \textstyle{} \sum_{j} (-1)^{j}  \binom{l}{j-b}\binom{j}{w} & \textstyle{} = \sum_{j} (-1)^{j} \binom{l}{j-b} \sum_v \binom{j-b}{v} \binom{b}{w-v} 
    \\& \textstyle{} =  \sum_v \binom{b}{w-v} \sum_{j} (-1)^{j} \binom{j-b}{v}  \binom{l}{j-b}
        \\& \textstyle{} =  \sum_v [v=l](-1)^{b+l} \binom{b}{w-v} 
        =  (-1)^{b+l} \binom{b}{w-l}.
\end{align*}
The first equality follows from Vandermonde's convolution formula, the second equality is a simple rearrangement of the sums ``\m{\sum_j}'' and ``\m{\sum_v}'', and the third equality follows from equation~\eqref{ceqn13}.
\item   Because the degree \m{w} polynomial  \m{\binom{u-dX}{w}} is a linear combination of the polynomials  \m{\binom{X}{0},\ldots,\binom{X}{w}}, we have 
\[\textstyle \binom{u-dj}{w} = (-d)^w \binom{j}{w} + h_{w-1} \binom{j}{w-1} + \cdots + h_{0} \binom{j}{0}\]
 for some \m{h_{w-1},\ldots,h_{0} \in \ZZ_p} that depend on \m{u,d,w}. The claim follows from the facts that, due to equation~\eqref{ceqn11}, %
 \[\textstyle{}  \sum_j (-1)^{j-b} \binom{l}{j-b}\binom{j}{w-m} = (-1)^{l} \binom{b}{w-m-l} = 0 \] for all \m{m \in \{1,\ldots,w\}} (since \m{l\GEQSLANT w}), and
 \[\textstyle{}  \sum_j (-1)^{j-b} \binom{l}{j-b}\binom{j}{w} = (-1)^{l} \binom{b}{w-l} = [w=l](-1)^l.  \]
\end{enumerate}
\end{proof}

\begin{lemma}\label{part2X5}
Let \m{\alpha \in \ZZ \cap [0,\ldots,\frac{r}{p+1}]} and let \m{\FAMILY{D_i}_{i\in \ZZ}} be a family of elements of \m{\ZZ_p} such that \m{D_i=0} for \m{i\not\in[0,\frac{r-\alpha}{p-1}]} and
\m{\textstyle \TEMPORARYD_{w} (D_\DARKCIRC) = 0} for all \m{w \in \{0,\ldots,\alpha-1\}}. Then the polynomial \[\textstyle \sum_i D_i x^{i(p-1)+\alpha}y^{r-i(p-1)-\alpha} \in \TILDE\SIGMA_r\] can be written in the form \m{\theta^\alpha h} for some polynomial \m{h} with integer coefficients.
\end{lemma}

\begin{proof}{Proof. }
The \m{\QQ_p}-span of \[\textstyle\left\{\binom{(p-1)\FORMX}{0},\ldots,\binom{(p-1)\FORMX}{\alpha-1}\right\}\] is the same as the \m{\QQ_p}-span of \[\textstyle\left\{\binom{\FORMX}{0},\ldots,\binom{\FORMX}{\alpha-1}\right\}.\] Note that here we do not put a restriction on the size of \m{\alpha}, and in particular \m{\alpha} may be larger than \m{p}: the polynomials \m{\binom{(p-1)\FORMX}{j}} and \m{\binom{\FORMX}{j}} are in \m{\QQ_p[X]} (and not necessarily in \m{\ZZ_p[X]}) and by the ``\m{\QQ_p}-span'' we mean the corresponding \m{\QQ_p}-vector subspace of \m{\QQ_p[X]}. Thus the condition that \m{\textstyle \TEMPORARYD_{w} (D_\DARKCIRC) = 0} for all \m{w \in \{0,\ldots,\alpha-1\}} %
 is equivalent to \[\textstyle \sum_i D_i \binom{i}{w} = 0\] for all \m{w \in \{0,\ldots,\alpha-1\}}. %
  The coefficients of \m{\theta^\alpha h} for any polynomial \m{h} (which has degree \m{r-\alpha(p+1)}) satisfy this set of \m{\alpha} equations. %
 We can find an \m{h} such that
\[\textstyle \sum_i D_i x^{i(p-1)+\alpha}y^{r-i(p-1)-\alpha}  = \theta^\alpha h + \sum_i D_i^{\prime} x^{i(p-1)+\alpha}y^{r-i(p-1)-\alpha}\] and \m{D_i^{\prime}=0} if \m{i\not\in\{0,\ldots,\alpha-2\}}. Since the matrix %
 \[\textstyle\det \left(\binom{i}{w}\right)_{0\LEQSLANT i,w<\alpha}\] 
 is lower triangular with \m{1}'s on the diagonal, it has full rank, and therefore 
  we have a set of \m{\alpha -1} constants \m{D_0^{\prime},\ldots,D_{\alpha-2}^{\prime}} that satisfy \m{\alpha} independent linear equations, so all of them must be zero. Hence \[\textstyle \sum_i D_i x^{i(p-1)+\alpha}y^{r-i(p-1)-\alpha} = \theta^\alpha h.\] Moreover,  \m{\theta^\alpha h} must be integral since the coefficients \m{\FAMILY{D_i}_{i\in \ZZ}} are integers. It is easy to prove that ``\m{\theta g} is integral'' \m{\Longrightarrow} ``\m{g} is integral'', implying by induction that \m{h} is integral as well.
\end{proof}

\begin{lemma}\label{part2X3a}
For \m{u,v,c \in \ZZ} let us define
\[\textstyle F_{u,v,c}(X) = \sum_w (-1)^{w-c} \binom{w}{c} \binom{X}{w}^\partial \binom{X+u-w}{v-w} \in \QQ_p[X].\]
Then
\[\textstyle F_{u,v,c}(X) =  \binom{u}{v-c} \binom{X}{c}^\partial-\binom{u}{v-c}^\partial \binom{X}{c}.\]
\end{lemma}

\begin{proof}{Proof. }
If \m{c<0} or \m{v < c} then the claim is trivial. Let  \m{v\GEQSLANT c \GEQSLANT 0}. It is enough to show that
 \m{\polyH^{\prime}(0)=0} for
\[\textstyle \polyH(\polyr) = \sum_w (-1)^{w-c} \binom{w}{c} \binom{\polyr+X}{w} \binom{X+u-w}{v-w} -  \binom{\polyr+X}{c} \binom{u-\polyr}{v-c}\in \QQ_p[X][\polyr].\]
 In fact, \m{\polyH(\polyr)} is the zero polynomial (over \m{\QQ_p[X]}) as can be seen from
\begin{align*}
\textstyle  \sum_w (-1)^{w-c} \binom{w}{c} \binom{\polyr+X}{w} \binom{X+u-w}{v-w}   & \textstyle  =  (-1)^{v-c} \binom{\polyr+X}{c} \sum_w \binom{\polyr+X-c}{w-c} \binom{v-u-X-1}{v-w}%
\\& \textstyle  = (-1)^{v-c} \binom{\polyr+X}{c} \binom{\polyr+v-u-c-1}{v-c}
\\& \textstyle  = \binom{\polyr+X}{c} \binom{u-\polyr}{v-c}.
\end{align*}
Here the first equality is a simple rewrite, \IE we use the equations
\[\textstyle \binom{w}{c} \binom{\polyr+X}{w}  =  \binom{\polyr+X}{c} \binom{\polyr+X-c}{w-c}  \text{ and }\binom{X+u-w}{v-w} =(-1)^{v-w} \binom{v-u-X-1}{v-w}.\] The second equality follows from Vandermonde's convolution formula, and the third equality is a simple rewrite as well.
\end{proof}

\textbf{Combinatorial identities involving  Stirling numbers.}
In this subsection we prove two lemmas about Stirling numbers. The reason  Stirling numbers appear in the proof is because they encode the change of basis of \m{\QQ_p[X]} from powers \m{\left\{1,X,X^2,\ldots\right\}} to falling factorials \m{\left\{1,X_1,X_2,\ldots\right\}}.

\begin{lemma}\label{part3X5} For \m{\alpha,\mu \in \ZZ_{\GEQSLANT 0}} and \m{\lambda\in\ZZ_{\GEQSLANT1}} let
\m{\textstyle L_\alpha(\lambda,\mu)} be the 
\m{(\alpha+1)\times(\alpha+1)} matrix with  %
 entries indexed by \m{0\LEQSLANT l,j\LEQSLANT\alpha} and defined by
    \[\textstyle (L_\alpha(\lambda,\mu))_{l,j} = \sum_{k=0}^{\alpha} \frac{j!}{l!} \left(\frac\mu\lambda\right)^k  \sfirst (l,k) \ssecond (k,j),\]
where \m{\sfirst (l,k)}  and \m{\ssecond (k,j)} are the Stirling numbers of the first and second kind, respectively (see section~\ref{secList}). Then
\[\textstyle L_\alpha(\lambda,\mu) \left(\binom{\lambda \FORMX}{0},\ldots,\binom{\lambda \FORMX}{\alpha}\right)^T = \left(\binom{\mu \FORMX}{0},\ldots,\binom{\mu \FORMX}{\alpha}\right)^T.\]
\end{lemma}

\begin{proof}{Proof. }
Straight from the definitions of \m{\sfirst(n,k)} and \m{\ssecond(n,k)} we have
\[\textstyle \left(\frac{\mu^j\sfirst(i,j)}{i!}\right)_{0\LEQSLANT i,j\LEQSLANT \alpha} \cdot \left(1,\ldots,\FORMX^\alpha\right)^T = \left(\binom{\mu\FORMX}{0},\ldots,\binom{\mu\FORMX}{\alpha}\right)^T\]
and
\[\textstyle \left(\frac{j! \ssecond(i,j)}{\lambda^i}\right)_{0\LEQSLANT i,j\LEQSLANT \alpha} \cdot \left(\binom{\lambda\FORMX}{0},\ldots,\binom{\lambda\FORMX}{\alpha}\right)^T = \left(1,\ldots,\FORMX^\alpha\right)^T.\] The claim then follows from the fact that
\[\textstyle L_\alpha(\lambda,\mu) = \left(\frac{\mu^j\sfirst(i,j)}{i!}\right)_{0\LEQSLANT i,j\LEQSLANT \alpha} \cdot \left(\frac{j! \ssecond(i,j)}{\lambda^i}\right)_{0\LEQSLANT i,j\LEQSLANT \alpha}.\]
\end{proof}

\begin{lemma}\label{lemmaComb2} 
For \m{\alpha \in \ZZ_{\GEQSLANT 0}} let \m{B_{\alpha}} be the 
\m{(\alpha+1)\times(\alpha+1)} matrix with %
 entries  indexed by \m{0\LEQSLANT i,j\LEQSLANT\alpha} and defined by
    \[\textstyle (B_\alpha)_{i,j} = j! \sum_{k,l=0}^{\alpha} \frac{(-1)^{i+l+k}}{l!} \binom{l}{i} (1-p)^{-k} \sfirst (l,k) \ssecond (k,j),\]
where \m{\sfirst (l,k)}  and \m{\ssecond (k,j)} are the Stirling numbers of the first and second kind, respectively (see section~\ref{secList}). Let \m{\{X_{i,j}\}_{i,j\GEQSLANT 0}} be a family of formal variables.
For \m{\beta \in \{0,\ldots,\alpha\}} %
 let \m{S({\alpha,\beta})} %
  be the \m{(\alpha+1)\times(\alpha+1)} matrix with %
 entries indexed by \m{0\LEQSLANT w,j\LEQSLANT\alpha} and defined by
 \[\textstyle S({\alpha,\beta})_{w,j} = \sum_{i=1}^\beta X_{i,j} \binom{i(p-1)}{w}.\]
 Then the matrix \m{B_\alpha S({\alpha,\beta})} is zero outside the rows indexed by \m{1,\ldots,\beta} and \[\textstyle (B_\alpha S({\alpha,\beta}))_{i,j} = X_{i,j}\]
for \m{i \in \{1,\ldots,\beta\}}.
\end{lemma}

\begin{proof}{Proof. }
Let \m{L=L_\alpha(p-1,1)} be the matrix defined in lemma~\ref{part3X5}. Then it follows from lemma~\ref{part3X5} that \[\textstyle (LS(\alpha,\beta))_{l,j} = \sum_{w=1}^\beta X_{w,j}\binom{w}{l}.\]
Let \m{E=\left(\binom{w}{l}\right)_{0\LEQSLANT l,w\LEQSLANT\alpha}}, so that \m{E^{-1} = ((-1)^{i+l}\binom{l}{i})_{0\LEQSLANT i,l\LEQSLANT\alpha}} and hence \[\textstyle{}B_\alpha = E^{-1}L.\]
Then
\begin{align*}
    \textstyle{} (B_\alpha S({\alpha,\beta}))_{i,j} & \textstyle{} = \sum_{l=0}^\alpha (E^{-1})_{i,l} (LS(\alpha,\beta))_{l,j}
    \\& \textstyle{} = \sum_{l=0}^\alpha (-1)^{i+l}\binom{l}{i} \sum_{w=1}^\beta X_{w,j}\binom{w}{l}
      \\& \textstyle{} = \sum_{w=1}^\beta X_{w,j} \sum_{l=0}^\alpha (-1)^{i+l}\binom{l}{i}\binom{w}{l}
            \\& \textstyle{} =  \sum_{w=1}^\beta [w=i] X_{w,j}  = [i \in \{1,\ldots,\beta\}] X_{i,j}.
\end{align*} The fourth equality follows from lemma~\ref{part1X5prime}, equation~\eqref{ceqn13}.
\end{proof}

\textbf{Two ad hoc combinatorial lemmas.}
 The statements of the two lemmas in this subsection are specifically tailored to be applicable in the main proof. This is why they appear overly specific, and why their proofs are unfortunately necessarily somewhat long.

\begin{lemma}\label{part4X5} Let \m{\alpha\in\ZZ_{\GEQSLANT1}}, let \m{\FORMX} and \m{\FORMY} be formal variables, and  let
\[\textstyle c_j = (-1)^j \alpha! \left( \frac{\FORMX+j+1}{j+1}\binom{\FORMY}{\alpha-j-1} + \binom{\FORMY}{\alpha-j}\right)\in\QQ[\FORMX,\FORMY]\subset \QQ(\FORMX,\FORMY)\]
  for \m{j \in \{1,\ldots,\alpha\}}. Let \m{M} %
 be the \m{(\alpha+1)\times(\alpha+1)} matrix over \m{\QQ(\FORMX,\FORMY)} with entries  indexed by \m{0\LEQSLANT w,j\LEQSLANT\alpha} and defined by
 \begin{align*}
        M_{w,0} & =\textstyle (-1)^w\frac{(\FORMY-\FORMX)\FORMX_w}{\FORMY_{w+1}}, %
        \\ M_{w,j} & = \textstyle \sum_v \binom{-j}{w-v} \binom{\FORMX+j}{v} \left(\binom{\FORMY+j-v}{j-v}-\binom{\FORMX+j-v}{j-v}\right), %
\end{align*}
for \m{w \in\{0,\ldots,\alpha\}} and \m{j \in \{1,\ldots,\alpha\}}. Let \m{c = (\FORMY_{\alpha},c_1,\ldots,c_\alpha)^T}. Then the first \m{\alpha} entries of
\[\textstyle Mc = (d_0,\ldots,d_\alpha)^T \] are zero, and \m{d_\alpha=\frac{(\FORMY-\FORMX)_{\alpha+1}}{\FORMY-\alpha}}. \end{lemma}

\begin{proof}{Proof. }
 First let us compute \m{d_0}. We have
\begin{align*}
        \textstyle \frac{d_0}{\alpha!} & \textstyle = %
        \textstyle \frac{\FORMY-\FORMX}{\FORMY}\binom{\FORMY}{\alpha} +  \sum_{j=1}^\alpha (-1)^j \left( \frac{\FORMX+j+1}{j+1}\binom{\FORMY}{\alpha-j-1} + \binom{\FORMY}{\alpha-j}\right) \left(\binom{\FORMY+j}{j} - \binom{\FORMX+j}{j}\right)%
       \\ & =\textstyle %
       \textstyle \frac{\FORMY-\FORMX}{\FORMY}\binom{\FORMY}{\alpha} + (\FORMX+1)\binom{\FORMY}{\alpha-1} +  \sum_{j=1}^\alpha (-1)^j \left( \frac{\FORMX+j+1}{j+1}\binom{\alpha-1}{j}\binom{\FORMY+j}{\alpha-1} + \binom{\alpha}{j} \binom{\FORMY+j}{\alpha}\right) %
       \\ & =\textstyle %
       -\frac{\FORMX}{\FORMY} \binom{\FORMY}{\alpha} + \sum_{j=0}^\alpha (-1)^j \left( \frac{\FORMX+j+1}{j+1}\binom{\alpha-1}{j}\binom{\FORMY+j}{\alpha-1} + \binom{\alpha}{j} \binom{\FORMY+j}{\alpha}\right) = 0.
\end{align*}
The first equality follows directly from the definition of \m{d_0}. The second equality follows from telescoping the sum
\begin{align*}
   &  \textstyle{}\sum_{j=1}^\alpha (-1)^j \left( \frac{\FORMX+j+1}{j+1}\binom{\FORMY}{\alpha-j-1} + \binom{\FORMY}{\alpha-j}\right) \binom{\FORMX+j}{j}\\ & \null\qquad\textstyle{}= \sum_{j=1}^\alpha (-1)^j \left( \binom{\FORMX+j+1}{j+1}\binom{\FORMY}{\alpha-j-1} + \binom{\FORMX+j}{j} \binom{\FORMY}{\alpha-j}\right)
   \\ & \null\qquad\textstyle{}=  \binom{\FORMX+\alpha+1}{\alpha+1} \binom{\FORMY}{-1} - \binom{\FORMX+1}{1} \binom{\FORMY}{\alpha-1}
\end{align*}
and using the equation \[\textstyle{}\binom{\FORMY}{\beta-j} \binom{\FORMY +j}{j} = \binom{\beta}{j}\binom{\FORMY+j}{\beta},\] which is true for \m{\beta,j\in \ZZ_{\GEQSLANT0}}.
 The third equality follows from the fact that the ``\m{j=0}'' term in the last sum is
 \[\textstyle{} (\FORMX+1) \binom{\FORMY}{\alpha-1}+\binom{\FORMY}{\alpha}.\] Thus we have shown that \m{d_0=0}.
Now let us suppose that \m{w \in \{1,\ldots,\alpha\}}. If \m{j \in \{1,\ldots,\alpha\}} then it is clear that \m{ M_{w,j}} is a polynomial that belongs to the ideal generated by \m{Y-X}, and due to lemma~\ref{part1X5}, equation~\eqref{ceqn8} we have
\[\textstyle M_{w,j} = \textstyle \sum_v  \binom{-j}{w-v} \binom{\FORMX+j}{v} \binom{\FORMY+j-v}{j-v}.\]
We want to show that \[\textstyle (\FORMY-\alpha)(Mc)_w = (\FORMY-\alpha) d_w =  \VERIFYIF{w=\alpha} (\FORMY-\FORMX)_{\alpha+1}.\] %
 Since \m{\FORMY_{\alpha+1} M_{w,0}} is a polynomial and \m{ M_{w,j}} is a polynomial for all \m{j\in\{1,\ldots,\alpha\}}, it follows that \m{(\FORMY-\alpha)(Mc)_w} is also a polynomial.
 The degree of \m{(\FORMY-\alpha)(Mc)_w} %
  is at most \m{\alpha+1}. Let us show that \m{(\FORMY-\alpha)(Mc)_w} belongs to the ideal generated by \m{\FORMY-\FORMX-t} for all \m{t \in \{0,\ldots,\alpha\}}. %
 If \m{t=0} then this is obvious (because \m{\FORMY_{\alpha+1} M_{w,0},M_{w,1},\ldots,M_{w,\alpha}} all belong to the ideal generated by \m{Y-X}), so let us suppose that \m{t \in \{1,\ldots, \alpha\}}. In that case \m{(Mc)_w} is  a polynomial of degree at most \m{\alpha-1} in the quotient ring \[\textstyle \QQ[\FORMX,\FORMY]/(\FORMY-\FORMX-t)\cong \QQ[\FORMX]\cong \QQ[\FORMY].\] (We can write an element of this quotient ring either as a polynomial in \m{X} or as a polynomial in \m{Y}, and it will have the same degree in either variable.) By lemma~\ref{part1X5}, equation~\eqref{ceqn8},
 \begin{align*}
     \textstyle M_{w,j}  & \textstyle{}= \sum_v  \binom{-j}{w-v} \binom{\FORMX+j}{v} \binom{\FORMX+t+j-v}{j-v}
    \\ & \textstyle{}= \sum_v  \binom{-j}{w-v} \binom{\FORMX+j}{v} \sum_l \binom{\FORMX+j-v}{j-v-l} \binom{t}{l}
     \\ & \textstyle{}= \sum_l \binom{t}{l} \sum_v  \binom{-j}{w-v} \binom{\FORMX+j}{v} \binom{\FORMX+j-v}{j-v-l}
          \\ & \textstyle{}= \sum_l \binom{t}{l} \binom{\FORMX+j}{j-l} \sum_v  \binom{-j}{w-v} \binom{j-l}{v}
                    \\ & \textstyle{}= \sum_l \binom{t}{l} \binom{\FORMX+j}{j-l} \binom{-l}{w}
                    \\ & \textstyle{} = (-1)^w \sum_{l} \binom{t}{l} \binom{\FORMX+j}{j-l}\binom{w+l-1}{w}.
 \end{align*} The second and fifth equalities follow from Vandermonde's convolution formula, the third equality is a simple rearrangement of the sums ``\m{\sum_v}'' and ``\m{\sum_l}'', the fourth equality follows from \[\textstyle{}\textstyle{}  \binom{\FORMX+j}{v} \binom{\FORMX+j-v}{j-v-l} =  \binom{\FORMX+j}{j-l} \binom{j-l}{v},\]
 and the sixth equality follows from 
 \m{\textstyle{} \binom{-l}{w} = (-1)^w\binom{w+l-1}{w}}.
 So \[\textstyle \frac{(-1)^{w}}{\alpha!} \sum_{j=1}^\alpha M_{w,j}c_j\] is equal to
\begin{align*}
    & \textstyle \sum_{l} \binom{t}{l} \binom{w+l-1}{w} \sum_{j>0} (-1)^j \binom{\FORMX+j}{j-l} \left(\frac{\FORMX+j+1}{j+1}\binom{\FORMX+t}{\alpha-j-1}+\binom{\FORMX+t}{\alpha-j}\right)
    \\ & \null\qquad \textstyle = 
    \sum_{l} \binom{t}{l} \binom{w+l-1}{w} \sum_{j>0} (-1)^j  \left(\frac{j-l+1}{j+1}\binom{\FORMX+j+1}{j-l+1} \binom{\FORMX+t}{\alpha-j-1}+\binom{\FORMX+j}{j-l}\binom{\FORMX+t}{\alpha-j}\right)
        \\ & \null\qquad \textstyle = 
    \sum_{l} l \binom{t}{l} \binom{w+l-1}{w} \sum_{j>0}  \frac{(-1)^j }{j} \binom{\FORMX+j}{j-l} \binom{\FORMX+t}{\alpha-j}
             \\ & \null\qquad \textstyle = 
     \sum_u \sum_{l} l \binom{t}{l} \binom{w+l-1}{w} \sum_{j>0}  \frac{(-1)^j }{j} \binom{\alpha-j+u}{u} \binom{\alpha-t}{j-l-u} \binom{\FORMX+t}{\alpha-j+u}
                \\ & \null\qquad \textstyle = 
    \sum_u \binom{\FORMX+t}{u} \sum_{l} l \binom{t}{l} \binom{w+l-1}{w}  \binom{\alpha-t}{\alpha-l-u} \sum_{j>0}  \frac{(-1)^j }{j} \binom{u}{u-\alpha+j}.
\end{align*} 
Let us provide some justification for this string of equalities. The first equality is a simple rewrite. For the second equality we replace \m{j+1} with \m{j} in the first term of the inner sum, \IE we use the equation
\[\textstyle \sum_{j>0} (-1)^j \frac{j-l+1}{j+1}\binom{\FORMX+j+1}{j-l+1} \binom{\FORMX+t}{\alpha-j-1} = \sum_{j>1} (-1)^{j-1} \frac{j-l}{j}\binom{\FORMX+j}{j-l} \binom{\FORMX+t}{\alpha-j}\]
to rewrite
\begin{align*} & \textstyle \sum_{j>0} (-1)^j  \left(\frac{j-l+1}{j+1}\binom{\FORMX+j+1}{j-l+1} \binom{\FORMX+t}{\alpha-j-1}+\binom{\FORMX+j}{j-l}\binom{\FORMX+t}{\alpha-j}\right) \\ & \textstyle \null\qquad  \vphantom{\sum_{j>0} (-1)^j  \left(\frac{j-l+1}{j+1}\binom{\FORMX+j+1}{j-l+1} \binom{\FORMX+t}{\alpha-j-1}+\binom{\FORMX+j}{j-l}\binom{\FORMX+t}{\alpha-j}\right)}= -\binom{\FORMX+1}{1-l}\binom{\FORMX+t}{\alpha-1} + \sum_{j>1} (-1)^j  \left(\frac{l-j}{j}+1\right)\binom{\FORMX+j}{j-l}\binom{\FORMX+t}{\alpha-j}
\\& \textstyle \null\qquad \vphantom{\sum_{j>0} (-1)^j  \left(\frac{j-l+1}{j+1}\binom{\FORMX+j+1}{j-l+1} \binom{\FORMX+t}{\alpha-j-1}+\binom{\FORMX+j}{j-l}\binom{\FORMX+t}{\alpha-j}\right)} = (l-1)\binom{\FORMX+1}{1-l}\binom{\FORMX+t}{\alpha-1} + \sum_{j>0} \frac{ (-1)^j l}{j} \binom{\FORMX+j}{j-l}\binom{\FORMX+t}{\alpha-j}\end{align*}
The extra term \m{(l-1)\binom{\FORMX+1}{1-l}\binom{\FORMX+t}{\alpha-1}} does not contribute to anything since
\[\textstyle  \sum_{l} \binom{t}{l} \binom{w+l-1}{w} (l-1)\binom{\FORMX+1}{1-l}\binom{\FORMX+t}{\alpha-1}\]
is zero for \m{l<0} (as then \m{\binom{t}{l} =0}), for \m{l=0} (as then \m{\binom{w+l-1}{w}=0}), for \m{l=1} (as then \m{l-1=0}), and for \m{l>1} (as then \m{\binom{\FORMX+1}{1-l}=0}).
For the third equality we use Vandermonde's convolution formula. 
 To be more precise, we write
\begin{align*}\textstyle    \binom{\FORMX+j}{j-l} \binom{\FORMX+t}{\alpha-j} & \textstyle{}=  \sum_u \binom{\alpha-t}{j-l-u}\binom{\FORMX+t-\alpha+j}{u} \binom{\FORMX+t}{\alpha-j} \\ & \textstyle{}=  \sum_u \binom{\alpha-t}{j-l-u}\binom{\alpha-j+u}{u} \binom{\FORMX+t}{\alpha-j+u}.\end{align*}
Finally, for the fourth equality we simply change the variable \m{u} to \m{u-\alpha+j}. Since the degree of the last polynomial \[\textstyle{}    \sum_u \binom{\FORMX+t}{u} \sum_{l} l \binom{t}{l} \binom{w+l-1}{w}  \binom{\alpha-t}{\alpha-l-u} \sum_{j>0}  \frac{(-1)^j }{j} \binom{u}{u-\alpha+j}\] is at most \m{\alpha-1}, we can replace the sum \m{\sum_u} with \m{\sum_{u=0}^{\alpha-1}}, and for \m{u} in this range we have
    \begin{align*}\textstyle \sum_{j>0}  \frac{(-1)^j }{j} \binom{u}{u-\alpha+j} & \textstyle{}\vphantom{\left.-\sum_{j}  \binom{X}{j+\alpha-u}^\partial \binom{u}{u-j}\right|_{X=0}}= \sum_{j}  \frac{(-1)^{j+\alpha-u}}{j+\alpha-u} \binom{u}{j} \\&\textstyle{}\vphantom{\left.-\sum_{j}  \binom{X}{j+\alpha-u}^\partial \binom{u}{u-j}\right|_{X=0}}= \left.-\sum_{j}  \binom{0}{j+\alpha-u}^\partial \binom{u}{u-j}\right.
    \\&\textstyle{}\vphantom{\left.-\sum_{j}  \binom{X}{j+\alpha-u}^\partial \binom{u}{u-j}\right|_{X=0}}= \left.-\frac{\partial}{\partial X}\sum_{j}  \binom{X}{j+\alpha-u} \binom{u}{u-j}\right|_{X=0}
    \\&\textstyle{}\vphantom{\left.-\sum_{j}  \binom{X}{j+\alpha-u}^\partial \binom{u}{u-j}\right|_{X=0}}= \left.- \binom{X+u}{\alpha}^\partial \right|_{X=0}
     \\&\textstyle{}\vphantom{\left.-\sum_{j}  \binom{X}{j+\alpha-u}^\partial \binom{u}{u-j}\right|_{X=0}}= \left.-\binom{u}{\alpha}^\partial\right.\\&\textstyle{}\vphantom{\left.-\sum_{j}  \binom{X}{j+\alpha-u}^\partial \binom{u}{u-j}\right|_{X=0}} = \left. \frac{(-1)^{\alpha-u}(\alpha-u-1)!u!}{\alpha!}.\right.%
    \end{align*}
    The second and sixth equalities follow from lemma~\ref{partHARMONIC}, and the fourth equality follows from Vandermonde's convolution formula. Therefore
    \[\textstyle{}  \textstyle \frac{(-1)^{w}}{\alpha!} \sum_{j=1}^\alpha M_{w,j}c_j = t  \sum_{u,l} \binom{\FORMX+t}{u} \frac{(-1)^{\alpha-u}(\alpha-u-1)!u!}{\alpha!} \binom{t-1}{l-1} \binom{w+l-1}{w}  \binom{\alpha-t}{\alpha-l-u} .\]
    Thus the equation we want to show, \IE 
    \[\textstyle \frac{(-1)^{w}}{\alpha!} (Mc)_w = \frac{(Y-X)X_w}{Y_{w+1}}\binom{Y}{\alpha} + \frac{(-1)^{w}}{\alpha!} \sum_{j=1}^\alpha M_{w,j}c_j = 0\] (modulo \m{Y-X-t}),
is equivalent to
\begin{align}\label{eqwr34}\textstyle \sum_{u,l} \binom{\FORMY}{u} \frac{(-1)^{\alpha-u-1}(\alpha-u-1)!u!}{\alpha!} \binom{t-1}{l-1} \binom{w+l-1}{w} \binom{\alpha-t}{\alpha-l-u} = \frac{(\FORMY-t)_w}{\FORMY_{w+1}}\binom{\FORMY}{\alpha}.\end{align}
 Here for convenience we look at the relevant elements as polynomials in \m{Y} (recall that we identify \m{\QQ[\FORMX,\FORMY]/(\FORMY-\FORMX-t)\cong \QQ[\FORMX]\cong \QQ[\FORMY]}).
 Let us denote the  left and right side of equation~\eqref{eqwr34} by \m{L(t)} and \m{R(t)}, respectively. To conclude that \m{L(t)=R(t)} for \m{t \in \{1,\ldots, \alpha\}}, it is evidently enough to show that \begin{align}\label{34te545y4}\textstyle \sum_{j=0}^m (-1)^j \binom{m}{j} L(j+1)=\sum_{j=0}^m (-1)^j \binom{m}{j}R(j+1)\end{align} for \m{m \in \{0,\ldots,\alpha-1\}}. 
 We can uniquely write each side of equation~\eqref{34te545y4} %
  in the form
 \[\textstyle h_{\alpha-1}\binom{Y}{\alpha-1} + \cdots + h_0 \binom{Y}{0}.\]
 To show that the polynomials on the two sides of equation~\ref{34te545y4} are equal, it is enough to show that the coefficients of \m{\binom{\FORMY}{u}} (\IE the corresponding \m{h_u}) are the same on both sides, for all \m{u \in \{0,\ldots,\alpha-1\}}. The coefficient of \m{\binom{\FORMY}{u}} on the left side of equation~\eqref{34te545y4} is
\[\textstyle{} \kappa_{\alpha,m,u}^{-1}\sum_{j,l} (-1)^{m+j} \binom{m}{j} \binom{j}{l-1} \binom{w+l-1}{w} \binom{\alpha-j-1}{l+u-j-1},\]
with \m{\kappa_{\alpha,m,u} = (-1)^{\alpha-m-u-1} \alpha \binom{\alpha-1}{u}}. 
 We have %
\begin{align*}\textstyle \sum_{j} (-1)^{j}\binom{m}{j} R(j+1) & \textstyle{}= \sum_{j}(-1)^{j}\binom{m}{j} \frac{(\FORMY-j-1)_w}{\FORMY_{w+1}}\binom{\FORMY}{\alpha}
\\& \textstyle{}= \frac{w!}{\FORMY_{w+1}}\binom{\FORMY}{\alpha} \sum_{j} (-1)^{j}\binom{m}{j} \binom{\FORMY-j-1}{w}
\\& \textstyle{}= \frac{w!}{\FORMY_{w+1}}\binom{\FORMY}{\alpha} \sum_{j} (-1)^{w+j}\binom{m}{j} \binom{w+j-\FORMY}{w}
\\& \textstyle{}= \frac{w!}{\FORMY_{w+1}}\binom{\FORMY}{\alpha} \sum_{j} (-1)^{w+j}\binom{m}{j} \sum_l \binom{j}{l}\binom{w-\FORMY}{w-l}
\\& \textstyle{}= \frac{w!}{\FORMY_{w+1}}\binom{\FORMY}{\alpha}  \sum_l \binom{w-\FORMY}{w-l} \sum_{j} (-1)^{w+j}\binom{m}{j} \binom{j}{l}
\\& \textstyle{}= \frac{w!}{\FORMY_{w+1}}\binom{\FORMY}{\alpha}   (-1)^{w+m} \sum_l [l=m] \binom{w-\FORMY}{w-l}
\\& \textstyle{}= \frac{w!}{\FORMY_{w+1}}\binom{\FORMY}{\alpha}   \binom{\FORMY-m-1}{w-m}
\\& \textstyle{}= 
\frac{w_{m}}{\alpha_{m+1}} \binom{\FORMY-m-1}{\alpha-m-1} \\&\textstyle{}= \sum_u \binom{\FORMY}{u} \frac{ w_{m}}{\alpha_{m+1}} \binom{-m-1}{\alpha-m-u-1}
\\&\textstyle{}= \sum_u \binom{\FORMY}{u} \frac{(-1)^{\alpha-m-u-1} w_{m}}{\alpha_{m+1}} \binom{\alpha-u-1}{m}.\end{align*}
Here the second, third, fifth, seventh, eighth, and tenth equalities are simple rewrites and rearrangements of sums, the fourth and ninth equalities follow from Vandermonde's convolution formula, and the sixth equality follows from  lemma~\ref{part1X5prime}, equation~\eqref{ceqn13}.
Therefore the coefficient of \m{\binom{\FORMY}{u}} on the right side of equation~\eqref{34te545y4} is \[\textstyle{} \kappa_{\alpha,m,u}^{-1} \alpha \binom{\alpha-1}{u} \frac{w_m}{\alpha_{m+1}} \binom{\alpha-u-1}{m} =  \kappa_{\alpha,m,u}^{-1} \binom{\alpha-m-1}{u} \binom{w}{m}.\]
Thus we want to show that
 \begin{align}\label{r34t5h4t3t}\textstyle \sum_{j,l} (-1)^{m+j} \binom{m}{j} \binom{j}{l-1} \binom{w+l-1}{w} \binom{\alpha-j-1}{l+u-j-1} = \binom{\alpha-m-1}{u} \binom{w}{m}. \end{align}  Let us show that equation~\eqref{r34t5h4t3t} is true more generally for all \m{ \alpha,m,u,w\in\ZZ_{\GEQSLANT 0}}. %
 Let \m{L(\alpha,u)} and \m{R(\alpha,u)} denote the left and right side of equation~\eqref{r34t5h4t3t}, respectively. Then it is easy to verify that \[\textstyle\asterisk\,(\alpha,u) = \asterisk\,(\alpha-1,u)+\asterisk\,(\alpha-1,u-1)\] for \m{\asterisk\in\{L,R\}}.  Therefore we only need to show equation~\eqref{r34t5h4t3t} in the boundary cases, \IE when either \m{u=0} or \m{\alpha=0}. Suppose  that \m{u=0}. If \m{l \in \{0,\ldots,j\}} then \m{\textstyle{} \binom{\alpha-j-1}{l-j-1} = 0}, and if \m{l \not\in \{0,\ldots,j+1\}} then \m{\binom{j}{l-1}=0}, so the only terms in \m{L(\alpha,0)} that are non-zero are the ones such that \m{l-1=j}. Thus, the equation is
 \[\textstyle \sum_{j} (-1)^{m+j} \binom{m}{j}  \binom{w+j}{j} =  \binom{w}{m}.\]
 This follows from 
 \[\textstyle \binom{w+j}{j} = (-1)^j\binom{-w-1}{j} \text{ and } \binom{w}{m} = (-1)^m \binom{m-w-1}{m}\]
 and Vandermonde's convolution formula.
 If \m{\alpha=0} then the equation is
 \begin{align*}
     \label{gg45h56g6h55}\textstyle \sum_{j,l} (-1)^{m+l+1} \binom{m}{j} \binom{j}{l-1} \binom{w+l-1}{w} \binom{l+u-1}{j} = \binom{m+u}{u} \binom{w}{m},
 \end{align*}
 which follows from the fact that
 \[\textstyle \textstyle %
 \sum_{w,m,u,j,l\GEQSLANT 0} (-1)^{m+l+1} \binom{m}{j} \binom{j}{l-1} \binom{w+l-1}{w} \binom{l+u-1}{j} Z_1^wZ_2^mZ_3^u \in \QQ(Z_1,Z_2,Z_3)\] 
 is equal to
 \begin{align*}
    & \textstyle{}  \sum_{w,m,u,j,l\GEQSLANT 0} (-1)^{l+1} \binom{m}{j} \binom{j}{l-1} \binom{-l}{w} \binom{l+u-1}{j} (-Z_1)^w(-Z_2)^mZ_3^u 
       \\ & \textstyle{}\null\qquad\vphantom{\frac{Z_3}{1+Z_2-Z_3-Z_2Z_3} \sum_{j,l\GEQSLANT 0} (-1)^{l+1}  \binom{j}{l-1}  \left(\frac{1}{Z_3-Z_1Z_3}\right)^{l}\left(-\frac{Z_2Z_3}{1+Z_2-Z_3-Z_2Z_3}\right)^j}=  \sum_{w,m,u,j,l\GEQSLANT 0} (-1)^{l+1} \binom{m}{j} \binom{j}{l-1} \binom{-l}{w} \binom{l+u-1}{j} (-Z_1)^w(-Z_2)^mZ_3^u 
              \\ & \textstyle{}\null\qquad\vphantom{\frac{Z_3}{1+Z_2-Z_3-Z_2Z_3} \sum_{j,l\GEQSLANT 0} (-1)^{l+1}  \binom{j}{l-1}  \left(\frac{1}{Z_3-Z_1Z_3}\right)^{l}\left(-\frac{Z_2Z_3}{1+Z_2-Z_3-Z_2Z_3}\right)^j}=  - \sum_{m,u,j,l\GEQSLANT 0}  \binom{m}{j} \binom{j}{l-1} \binom{l+u-1}{j} (Z_1-1)^{-l}(-Z_2)^mZ_3^u 
                            \\ & \textstyle{}\null\qquad\vphantom{\frac{Z_3}{1+Z_2-Z_3-Z_2Z_3} \sum_{j,l\GEQSLANT 0} (-1)^{l+1}  \binom{j}{l-1}  \left(\frac{1}{Z_3-Z_1Z_3}\right)^{l}\left(-\frac{Z_2Z_3}{1+Z_2-Z_3-Z_2Z_3}\right)^j}=  -\frac{1}{1+Z_2}\sum_{u,j,l\GEQSLANT 0}   \binom{j}{l-1} \binom{l+u-1}{j} (Z_1-1)^{-l}\left(-\frac{Z_2}{1+Z_2}\right)^j  Z_3^u                    \\ & \textstyle{}\null\qquad\vphantom{\frac{Z_3}{1+Z_2-Z_3-Z_2Z_3} \sum_{j,l\GEQSLANT 0} (-1)^{l+1}  \binom{j}{l-1}  \left(\frac{1}{Z_3-Z_1Z_3}\right)^{l}\left(-\frac{Z_2Z_3}{1+Z_2-Z_3-Z_2Z_3}\right)^j}=  -\frac{Z_3}{1+Z_2-Z_3-Z_2Z_3} \sum_{j,l\GEQSLANT 0}   \binom{j}{l-1}  \left(-\frac{1}{Z_3-Z_1Z_3}\right)^{l}\left(-\frac{Z_2Z_3}{1+Z_2-Z_3-Z_2Z_3}\right)^j          \\ & \textstyle{}\null\qquad\vphantom{\frac{Z_3}{1+Z_2-Z_3-Z_2Z_3} \sum_{j,l\GEQSLANT 0} (-1)^{l+1}  \binom{j}{l-1}  \left(\frac{1}{Z_3-Z_1Z_3}\right)^{l}\left(-\frac{Z_2Z_3}{1+Z_2-Z_3-Z_2Z_3}\right)^j}=  \frac{1}{(1-Z_1)(1+Z_2-Z_3-Z_2Z_3)} \sum_{j\GEQSLANT 0}   \left(\frac{(1-Z_3+Z_1Z_3)Z_2}{(1-Z_1)(1+Z_2-Z_3-Z_2Z_3)}\right)^{j}
                          \\ & \textstyle{}\null\qquad\vphantom{\frac{Z_3}{1+Z_2-Z_3-Z_2Z_3} \sum_{j,l\GEQSLANT 0} (-1)^{l+1}  \binom{j}{l-1}  \left(\frac{1}{Z_3-Z_1Z_3}\right)^{l}\left(-\frac{Z_2Z_3}{1+Z_2-Z_3-Z_2Z_3}\right)^j}=  %
       \frac{1}{1-Z_1-Z_3-Z_1Z_2+Z_1Z_3}\in \QQ(Z_1,Z_2,Z_3)
 \end{align*}
 and the fact that \[\textstyle \sum_{w,m,u\GEQSLANT 0} \binom{m+u}{u} \binom{w}{m} Z_1^wZ_2^mZ_3^u\in \QQ(Z_1,Z_2,Z_3)\] is equal to 
 \begin{align*}
    &\textstyle{} \frac{1}{1-Z_1}\sum_{m,u\GEQSLANT 0} \binom{m+u}{u} \left(\frac{Z_1Z_2}{1-Z_1}\right)^mZ_3^u
     \\ & \textstyle{}\null\qquad \vphantom{\frac{1}{1-Z_1-Z_1Z_2}\sum_{u\GEQSLANT 0}  \left(\frac{Z_3-Z_1Z_3}{1-Z_1-Z_1Z_2}\right)^{u}}= \frac{1}{1-Z_1-Z_1Z_2}\sum_{u\GEQSLANT 0}  \left(\frac{Z_3-Z_1Z_3}{1-Z_1-Z_1Z_2}\right)^{u}
          \\ & \textstyle{}\null\qquad \vphantom{\frac{1}{1-Z_1-Z_1Z_2}\sum_{u\GEQSLANT 0}  \left(\frac{Z_3-Z_1Z_3}{1-Z_1-Z_1Z_2}\right)^{u}}= \frac{1}{1-Z_1-Z_3-Z_1Z_2+Z_1Z_3}\in \QQ(Z_1,Z_2,Z_3).
 \end{align*} To summarize, we have shown that  \m{(\FORMY-\alpha) d_w = (\FORMY-\alpha)(Mc)_w} has degree at most \m{\alpha+1} and belongs to the ideal generated by \m{\FORMY-\FORMX-t} for all \m{t \in \{0,\ldots,\alpha\}}. Therefore it must be a constant multiple of \m{(\FORMY-\FORMX)_{\alpha+1}}. We can deduce that the constant is indeed \m{\VERIFYIF{w=\alpha}} by comparing the coefficients of \m{\FORMY^{\alpha+1}}, thus completing the proof.
\end{proof}

\begin{lemma}\label{lemmaComb4} 
Suppose that \m{s,\alpha,\beta\in\ZZ} are such that
\[\textstyle 1 \LEQSLANT \beta \LEQSLANT \alpha \LEQSLANT \frac{s}{2}-2 \LEQSLANT \frac{p-5}{2}.\]
Let \m{B=B_\alpha} denote the matrix defined in lemma~\ref{lemmaComb2}.
Let \m{M} denote the \m{(\alpha+1)\times(\alpha+1)} matrix with entries in \m{\FF_p} such that if \m{i\in\{1,\ldots,\beta\}} and \m{j\in\{0,\ldots,\alpha\}} then
\[\textstyle M_{i,j} = \binom{\beta}{i} \cdot \left\{\begin{array}{rl}
    \binom{s-\alpha-\beta+i}{i}^{-1} (-1)^{i+1} & \text{if } j = 0, \\
    \binom{s-\alpha-\beta+j}{j-i} & \text{if }j>0,
\end{array}\right.\]
and if \m{i\in\{0,\ldots,\alpha\}\backslash\{1,\ldots,\beta\}} and \m{j\in\{0,\ldots,\alpha\}} then \m{M_{i,j}} is the reduction modulo \m{p} of
 \begin{align*}
       & \textstyle p^{-\VERIFYIF{j=0}}  \sum_{w=0}^\alpha B_{i,w} \sum_v   \binom{-j}{w-v} \binom{s+\beta(p-1)-\alpha+j}{v}^\partial  
        \sum_{u=0}^{\beta} \binom{s+\beta(p-1)-\alpha+j-v}{u(p-1)+j-v}
        \\
        & \textstyle \null\qquad - \VERIFYIF{i=0} %
        \textstyle  \binom{s+\beta(p-1)-\alpha+j} {j}^\partial %
     -\VERIFYIF{j=0} \sum_{w=0}^\alpha B_{i,w} \binom{s+\beta(p-1)-\alpha}{w}  \frac{(-1)^{w}w!}{(s-\alpha)_{w+1}}.
\end{align*}
Then there is a solution of
\[M (z_0,\ldots,z_\alpha)^T = (1,0,\ldots,0)^T\] such that \m{z_0 \ne 0}. 
\end{lemma}

\begin{proof}{Proof. }
Let us simplify  \m{M_{i,j}} for \m{i\in\{0,\ldots,\alpha\}\backslash\{1,\ldots,\beta\}} and \m{j\in\{0,\ldots,\alpha\}}. The matrix \m{\overline{B}} can be defined as
\[\textstyle \overline{B} \left(\binom{-X}{0},\ldots,\binom{-X}{\alpha}\right)^T = \left(\sum_{l=0}^\alpha (-1)^{i+l} \binom{l}{i} \binom{X}{l}\right)_{0\LEQSLANT i \LEQSLANT \alpha}^T.\]
Since \[\textstyle\binom{X}{l} = (-1)^l \binom{-X+l-1}{l} = (-1)^l \sum_w \binom{-X}{w}\binom{l-1}{l-w}\] for \m{l\in\ZZ_{\GEQSLANT1}}, it follows that
\begin{align}\label{eq2525}\textstyle \overline{B}_{i,w} = \VERIFYIF{(i,w)=(0,0)} + \sum_{l=1}^\alpha (-1)^{i} \binom{l}{i} \binom{l-1}{l-w}.\end{align}
In particular, \m{\overline{B}_{0,w} = \binom{\alpha}{w}}, and %
\begin{align*}
    \textstyle \sum_{w=0}^\alpha \overline{B}_{i,w} \binom{-j}{w-v} & \textstyle = \sum_{w=0}^\alpha \left(\VERIFYIF{(i,w)=(0,0)} + \sum_{l=1}^\alpha (-1)^{i} \binom{l}{i} \binom{l-1}{l-w}\right) \binom{-j}{w-v}
    \\ & \textstyle = \vphantom{\sum_{w=0}^\alpha \left(\VERIFYIF{(i,w)=(0,0)} + \sum_{l=1}^\alpha (-1)^{i} \binom{l}{i} \binom{l-1}{l-w}\right) \binom{-j}{w-v}} \VERIFYIF{(i,v)=(0,0)} + \sum_{w,l=1}^\alpha  (-1)^{i} \binom{l}{i} \binom{l-1}{w-1} \binom{-j}{w-v}
    \\ & \textstyle = \vphantom{\sum_{w=0}^\alpha \left(\VERIFYIF{(i,w)=(0,0)} + \sum_{l=1}^\alpha (-1)^{i} \binom{l}{i} \binom{l-1}{l-w}\right) \binom{-j}{w-v}} \VERIFYIF{(i,v)=(0,0)} + (-1)^{i+v} \sum_{l=1}^\alpha  (-1)^{l} \binom{l}{i}  \binom{j-v}{l-v}
        \\ & \textstyle = \vphantom{\sum_{w=0}^\alpha \left(\VERIFYIF{(i,w)=(0,0)} + \sum_{l=1}^\alpha (-1)^{i} \binom{l}{i} \binom{l-1}{l-w}\right) \binom{-j}{w-v}}(-1)^{i+v} \sum_{l=0}^\alpha  (-1)^{l} \binom{l}{i}  \binom{j-v}{l-v}
        .
\end{align*} The first equality follows from equation~\eqref{eq2525}. 
The third equality follows from lemma~\ref{part1X5prime}, equation~\eqref{ceqn10}. When \m{0\LEQSLANT v\LEQSLANT j\LEQSLANT \alpha}, due to lemma~\ref{part1X5prime}, equation~\eqref{ceqn11}, %
 \[\textstyle (-1)^{i+v} \sum_{l=0}^\alpha  (-1)^{l} \binom{l}{i}  \binom{j-v}{l-v} = (-1)^{i+j+v}  \binom{v}{j-i}.\] Hence we can simplify \m{M_{i,j}} for \m{i\not\in \{1,\ldots,\beta\}} as
\[\textstyle \left\{\begin{array}{rl}
     \sum_{w=0}^\alpha  \overline{B}_{i,w} \left(\beta  \binom{s-\alpha-\beta}{w}^\partial-\binom{s-\alpha-\beta}{w}\right) \frac{(-1)^{w}w!}{(s-\alpha)_{w+1}} & \text{if } j = 0, \\
    \sum_{v} (-1)^{i+j+v} \binom{v}{j-i} \binom{s-\alpha-\beta+j}{v}^\partial  \binom{s-\alpha+j-v}{j-v} - \VERIFYIF{i=0} \textstyle  \binom{s-\alpha-\beta+j} {j}^\partial & \text{if } j>0.
\end{array}\right.\]
Here the simplification for \m{j=0} follows from the equation
\[\textstyle{} %
        \sum_{u} \binom{s+\beta(p-1)-\alpha-w}{u(p-1)-w} = \frac{(-1)^w w! \beta p}{(s-\alpha)_{w+1}},\] which follows from lemma~\ref{part1X5}, equation~\eqref{ceqn3}. If \m{i\in \{\beta+1,\ldots,\alpha\}} then, due to lemma~\ref{part2X3a}, \begin{align*}\textstyle \sum_v (-1)^{i+j+v} \binom{v}{i-j} \binom{s-\alpha-\beta+j}{v}^\partial \binom{s-\alpha+j-v}{j-v} 
  & \textstyle = F_{\beta,j,j-i}(s-\alpha-\beta+j)  \\ & \textstyle = \binom{\beta}{i}^\partial \binom{s-\alpha-\beta+j}{j-i},\end{align*}
  so we can further simplify 
\m{M_{i,j}} for \m{i\in \{\beta+1,\ldots,\alpha\}} as
\[\textstyle \left\{\begin{array}{rl}
     \sum_{w=0}^\alpha \overline{B}_{i,w} \left(\beta  \binom{s-\alpha-\beta}{w}^\partial-\binom{s-\alpha-\beta}{w}\right) \frac{ (-1)^{w}w!}{(s-\alpha)_{w+1}} & \text{if } j = 0, \\
   \binom{\beta}{i}^\partial \binom{s-\alpha-\beta+j}{j-i} - \VERIFYIF{i=0} \textstyle  \binom{s-\alpha-\beta+j} {j}^\partial & \text{if } j>0.
\end{array}\right.\]
It follows immediately from these expressions that all entries of \m{M} that are below the diagonal and are not in the zeroth column must be zero, since \[\textstyle \binom{s-\alpha-\beta+j}{j-i}=0\] %
 when \m{j<i}. It also follows that all entries in the zeroth row except for the one that is in the zeroth column must be zero, since the only non-zero term in the relevant sum
\[\textstyle \sum_{v} (-1)^{j+v} \binom{v}{j} \binom{s-\alpha-\beta+j}{v}^\partial  \binom{s-\alpha+j-v}{j-v}\] is the one with \m{v=j}. This is because \m{\textstyle \binom{v}{j}=0} for \m{v<j} and \m{\textstyle \binom{s-\alpha+j-v}{j-v}=0} for \m{v>j}. Therefore
\begingroup
\renewcommand*{\arraystretch}{1.5}
\[\textstyle{}M = \begin{bmatrix} M_{0,0} & 0 \\ v & U
\end{bmatrix}\]
\endgroup
for some vector \m{v} and some upper triangular \m{\alpha\times\alpha} matrix \m{U}.\footnote{For example, when \m{\alpha\GEQSLANT 4} this means that
\begingroup
\renewcommand*{\arraystretch}{1.5}
\[\textstyle{}M = \begin{bmatrix} M_{0,0} & 0 & 0 & 0 & \cdots & 0 \\
M_{1,0} & M_{1,1} & M_{1,2} & M_{1,3} & \cdots & M_{1,\alpha} \\
M_{2,0} & 0 & M_{2,2} & M_{2,3} & \cdots & M_{2,\alpha}
\\
M_{3,0} & 0 & 0 & M_{3,3} & \cdots & M_{3,\alpha}
\\
\vdots & \vdots & \vdots & \vdots & \ddots & \vdots \\
M_{\alpha,0} & 0 & 0 & 0 & \cdots & M_{\alpha,\alpha}
\end{bmatrix}.\]
\endgroup}
Consequently, the equation \[M (z_0,\ldots,z_\alpha)^T = (1,0,\ldots,0)^T\] has a solution such that \m{z_0 \ne 0} as long as the diagonal entries of \m{M} are non-zero.
 Since \m{\overline{B}_{0,w}=\binom{\alpha}{w}}, we have %
\[\textstyle \frac{(s-\alpha)_{\alpha+1}}{\alpha!\beta} M_{0,0} = \textstyle \sum_{w=0}^\alpha (-1)^{w} \left(\binom{s-\alpha-\beta}{w}^\partial-\frac1\beta\binom{s-\alpha-\beta}{w}\right) \binom{s-\alpha-w-1}{\alpha-w}.\] Due to  lemma~\ref{part1X5prime}, equation~\eqref{ceqn9},
\[\textstyle{} -\frac1\beta\sum_{w=0}^\alpha (-1)^{w} \binom{s-\alpha-\beta}{w} \binom{s-\alpha-w-1}{\alpha-w} = -\frac1\beta\binom{\beta-1}{\alpha} = 0.\]
Due to lemma~\ref{part2X3a},
\[\textstyle{} \sum_{w=0}^\alpha (-1)^{w} \binom{s-\alpha-\beta}{w}^\partial \binom{s-\alpha-w-1}{\alpha-w} = F_{\beta-1,\alpha,0}(s-\alpha-\beta) = -\binom{\beta-1}{\alpha}^\partial.\]
Thus, due to lemma~\ref{partHARMONIC}, \[\textstyle{}M_{0,0} = -\frac{\alpha!\beta}{(s-\alpha)_{\alpha+1}}\binom{\beta-1}{\alpha}^\partial = \frac{(-1)^{\alpha+\beta+1}}{(\alpha+1) \binom{\alpha}{\beta}\binom{s-\alpha}{\alpha+1}}\ne 0.\]
The diagonal entries \m{M_{1,1},\ldots,M_{\beta,\beta}} are equal to \m{\binom{\beta}{1},\ldots,\binom{\beta}{\beta}}, respectively, all of which are non-zero. For \m{j\in\{\beta+1,\ldots,\alpha\}} %
 we have, due to lemma~\ref{partHARMONIC},
\[\textstyle M_{j,j} = \binom{\beta}{j}^\partial %
 =\frac{(-1)^{\beta+j+1}}{(\beta+1)\binom{j}{\beta+1}} \ne 0.\] So the diagonal entries of \m{M} are indeed non-zero.
\end{proof}

\section{Computing \texorpdfstring{\m{\overline\Theta_{k,a}}}{\textTheta}}\label{secComputing}

This section comprises three lemmas about properties of the ideal \m{\Img(T-a)}. In sections~\ref{secProof} and~\ref{secProof2} we use these three lemmas to show that
\m{\overline\Theta_{k,a}} is irreducible under the conditions of theorem~\ref{MAINTHM1aa}. As remarked in section~\ref{secNot}, the main strategy we use to compute
\m{\overline\Theta_{k,a}} is
to find elements of \m{\SCR I_a} that represent non-trivial elements of the subquotients \m{\FAMILY{\FILT N_\alpha}_{0\LEQSLANT\alpha<\TEMPORARYI}}, and this relates to \m{\Img(T-a)} because if an element of \m{\COMPACTIND{KZ}{G} \SIGMA_{r}} is the reduction modulo \m{\maxideal} of an integral element of \m{\Img(T-a)} then it is also in \m{\SCR I_a}.

We refer to section~\ref{secList} for the relevant notation. We recall that  \m{\TEMPORARYI\LEQSLANT\frac{p-1}{2}<p} and \m{k>p^{100}} (and therefore \m{r=k-2>p^{99}}). Let us write \m{n=p^2}, so that \m{r>np^2}. If \m{w \in \{0,\ldots,\lfloor \frac{r}{p+1}\rfloor\}} and \m{u \in \{-w,\ldots,r-wp\}} then we write
\[\textstyle{}\theta^wx^uy^{r-w(p+1)-u}=\textstyle{}\sum_{j} (-1)^j \binom{w}{j} x^{u+w+j(p-1)}y^{r-w-j(p-1)-u},\] which is an element of \m{\COMPACTIND{KZ}{G} \SIGMA_{r}}. That is, the notation \m{\theta^wx^uy^{v}} makes sense even if either \m{u} or \m{v} is in \m{\{-w,\ldots,-1\}} since \m{\theta^w} is divisible by \m{(xy)^w}.

\begin{lemma}\label{LEMMAX2}
{
Suppose that \m{\alpha \in \{0,\ldots,\TEMPORARYI-1\}}.
\begin{enumerate}
    \item \label{111}
We have
  \begin{align*}
      \textstyle & \textstyle (T-a) \left(1 \sqbrackets{}{KZ}{\QQp}{ \theta^{n}x^{\alpha-n}y^{r-np-\alpha}}\right)
      \\ & \textstyle \null\qquad 
      = \sum_j (-1)^j \binom{n}{j} p^{j(p-1)+\alpha} \begin{Lmatrix}
      1 & 0 \\ 0 & p
      \end{Lmatrix} \sqbrackets{}{KZ}{\QQp}{x^{j(p-1)+\alpha}y^{r-j(p-1)-\alpha}} 
      \\ & \textstyle \null\qquad\qquad - a \sum_j (-1)^j \binom{n-\alpha}{j} \sqbrackets{}{KZ}{\QQp}{\theta^\alpha x^{j(p-1)}y^{r-j(p-1)-\alpha(p+1)}} + \BIGO{p^n}.
  \end{align*}
\item \label{222} The submodule \m{\Img(T-a)\subset\COMPACTIND{KZ}{G} \TILDE\SIGMA_{r}} contains
  \begin{align*}
      & \textstyle \sum_{i} \left(\sum_{l=\beta-\gamma}^{\beta} C_l \binom{r-\beta+l}{i(p-1)+l}\right)\sqbrackets{}{KZ}{\QQp}{x^{i(p-1)+\beta}y^{r-i(p-1)-\beta}} \\ & \textstyle \null\qquad\qquad\qquad\qquad\qquad\qquad\qquad\qquad\qquad+ \BIGO{ap^{-\beta+v_C}+p^{p-1}}
  \end{align*}
  for all \m{0 \LEQSLANT \beta \LEQSLANT \gamma < \TEMPORARYI} and all families \m{\FAMILY{C_l}_{l\in \ZZ}} of elements of \m{\ZZ_p}, where
  \[\textstyle v_C = \min_{\beta-\gamma\LEQSLANT l \LEQSLANT \beta}(\PADICVALUATION(C_l)+l).\]
  The \m{\BIGO{ap^{-\beta+v_C}+p^{p-1}}} term is equal to  \m{\BIGO{p^{p-1}}} plus
  \[\textstyle -\frac{ap^{-\beta}}{p-1} \sum_{l=\beta-\gamma}^{\beta} C_{l}p^l\sum_{0\ne\MU\in\FF_p} [\MU]^{-l}\begin{Lmatrix}
p & [\MU] \\ 0 & 1
\end{Lmatrix}  \sqbrackets{}{KZ}{\QQp}{\theta^n  x^{\beta-l-n}y^{r-np-\beta+l}}.\]
\end{enumerate}
}
\end{lemma}

\begin{proof}{Proof. }
(\ref{111}) Let us recall the formula
    \[\textstyle \HECKET ({\gamma}\sqbrackets{}{KZ}{\QQp}{v}) = {\sum_{\MU\in\FF_p} \gamma\begin{Lmatrix}
    p &[\MU] \\ 0 & 1
    \end{Lmatrix} }\sqbrackets{}{KZ}{\QQp}{\left(\begin{Lmatrix}
    1 & -[\MU] \\ 0 & p
    \end{Lmatrix} \cdot v\right)} + \gamma\begin{Lmatrix}
    1 & 0 \\ 0 & p
    \end{Lmatrix} \sqbrackets{}{KZ}{\QQp}{\left(\begin{Lmatrix}
    p & 0 \\ 0 & 1
    \end{Lmatrix} \cdot v\right)}\] for \m{T}. Directly from the definition of \m{\theta} we have
    \begin{align*}\textstyle {\gamma}\sqbrackets{}{KZ}{\QQp}{v} & \textstyle = 1 \sqbrackets{}{KZ}{\QQp}{ \theta^{n}x^{\alpha-n}y^{r-np-\alpha}}\\ & \textstyle  = \sum_{j=0}^n (-1)^j \binom{n}{j} \sqbrackets{}{KZ}{\QQp}{ x^{\alpha+j(p-1)}y^{r-j(p-1)-\alpha}}\end{align*}
    If we apply the formula for \m{T} to \m{{\gamma}\sqbrackets{}{KZ}{\QQp}{v}= 1 \sqbrackets{}{KZ}{\QQp}{ \theta^{n}x^{\alpha-n}y^{r-np-\alpha}}},
    then the first part of the formula, the part
    \[\textstyle {\sum_{\MU\in\FF_p} \gamma\begin{Lmatrix}
    p &[\MU] \\ 0 & 1
    \end{Lmatrix} }\sqbrackets{}{KZ}{\QQp}{\left(\begin{Lmatrix}
    1 & -[\MU] \\ 0 & p
    \end{Lmatrix} \cdot v\right)},\]
    can be expanded into
\[\textstyle {\sum_{\MU\in\FF_p} \sum_{j,\TEMPORARYA\GEQSLANT 0}   \begin{Lmatrix}
    p &[\MU] \\ 0 & 1
    \end{Lmatrix} }\sqbrackets{}{KZ}{\QQp}{(-1)^j \binom{n}{j}\binom{r-j(p-1)-\alpha}{\TEMPORARYA}(-[\MU])^{r-j(p-1)-\alpha-\TEMPORARYA}p^\TEMPORARYA x^{r-\TEMPORARYA }y^\TEMPORARYA }.\]
The part of this sum with \m{\TEMPORARYA \GEQSLANT n} is in \m{\BIGO{p^n}}, and moreover \m{(-[\MU])^{r-j(p-1)-\alpha-\TEMPORARYA }} is independent of \m{j} when \m{\TEMPORARYA <n} since \m{r>p^{99}} implies \[\textstyle r-j(p-1)-\alpha-\TEMPORARYA  > r-np-p-n > 0.\] The coefficient of \m{x^{r-\TEMPORARYA }y^\TEMPORARYA } vanishes when \m{\TEMPORARYA <n} since
\[\textstyle \sum_j (-1)^j \binom{n}{j} \binom{r-j(p-1)-\alpha}{\TEMPORARYA } = 0,\]
due to lemma~\ref{part1X5prime}, equation~\eqref{ceqn5} applied to \m{(u,b,d,l,w)=(r-\alpha,0,p-1,n,\TEMPORARYA)}. The second part of the formula, the part
\[\textstyle \gamma\begin{Lmatrix}
    1 & 0 \\ 0 & p
    \end{Lmatrix} \sqbrackets{}{KZ}{\QQp}{\left(\begin{Lmatrix}
    p & 0 \\ 0 & 1
    \end{Lmatrix} \cdot v\right)},\]
is precisely
\[\textstyle \sum_j (-1)^j \binom{n}{j} p^{j(p-1)+\alpha} \begin{Lmatrix}
      1 & 0 \\ 0 & p
      \end{Lmatrix} \sqbrackets{}{KZ}{\QQp}{x^{j(p-1)+\alpha}y^{r-j(p-1)-\alpha}},\]
      and \m{-a \sqbrackets{}{KZ}{\QQp}{ \theta^{n}x^{\alpha-n}y^{r-np-\alpha}}} is precisely
      \[\textstyle - a \sum_j (-1)^j \binom{n-\alpha}{j} \sqbrackets{}{KZ}{\QQp}{\theta^\alpha x^{j(p-1)}y^{r-j(p-1)-\alpha(p+1)}}.\]
      By adding these three terms up we obtain the required expression for \[\textstyle(T-a) \left(1 \sqbrackets{}{KZ}{\QQp}{ \theta^{n}x^{\alpha-n}y^{r-np-\alpha}}\right).\]
     
(\ref{222}) Let us multiply the equation in the statement of part~\ref{111} on the left by
\[\textstyle C_{\beta-\alpha}p^{-\alpha}\sum_{0\ne\MU\in\FF_p} [\MU]^{\alpha-\beta}\begin{Smatrix}
p & [\MU] \\ 0 & 1
\end{Smatrix}\]
(that is, let us act on both sides of the equation in part~\ref{111} by the above element of \m{\QQp[G]}). Since \m{\Img(T-a)} is a \m{G}-module, both sides still end up in \m{\Img(T-a)}. The ``\m{j=0}'' part of the first sum on the right side becomes
\begin{align*}
      \textstyle & \textstyle C_{\beta-\alpha} \sum_{0\ne\MU\in\FF_p} [\MU]^{\alpha-\beta}\begin{Smatrix}
1 & [\MU] \\ 0 & 1
\end{Smatrix}\sqbrackets{}{KZ}{\QQp}{x^{\alpha}y^{r-\alpha}}
      \\ & \textstyle \null\qquad 
      = C_{\beta-\alpha} \sqbrackets{}{KZ}{\QQp}{\sum_{0\ne\MU\in\FF_p} [\MU]^{\alpha-\beta}\begin{Smatrix}
1 & [\MU] \\ 0 & 1
\end{Smatrix} x^{\alpha}y^{r-\alpha}}
 \\ & \textstyle \null\qquad 
      = C_{\beta-\alpha} \sqbrackets{}{KZ}{\QQp}{\sum_{0\ne\MU\in\FF_p} [\MU]^{\alpha-\beta} x^{\alpha}([\MU]x+y)^{r-\alpha}}
      \\ & \textstyle \null\qquad 
      =
      (p-1) \sum_i C_{\beta-\alpha}\binom{r-\alpha}{i(p-1)+\beta-\alpha} \sqbrackets{}{KZ}{\QQp}{x^{i(p-1)+\beta}y^{r-i(p-1)-\beta}}.
  \end{align*}
The third equality follows from the fact that \[\textstyle \sum_{0\ne\MU\in\FF_p} [\MU]^{w} = (p-1)\VERIFYIF{w\equiv_{p-1} 0}\] for \m{w\in\ZZ}.
 The ``\m{j>0}'' part of the first sum on the right becomes \m{\BIGO{p^{p-1}}}. The rest of the right side becomes
\[\textstyle -aC_{\beta-\alpha}p^{-\alpha}\sum_{0\ne\MU\in\FF_p} [\MU]^{\alpha-\beta}\begin{Smatrix}
p & [\MU] \\ 0 & 1
\end{Smatrix} \sqbrackets{}{KZ}{\QQp}{\theta^n  x^{\alpha-n}y^{r-np-\alpha}} + \BIGO{p^{n-\alpha}}.\]
Thus, since \m{n-\alpha > p^2-p>p}, we get that \m{\Img(T-a)} contains
\begin{align*}
    & \textstyle (p-1) \sum_i C_{\beta-\alpha}\binom{r-\alpha}{i(p-1)+\beta-\alpha} \sqbrackets{}{KZ}{\QQp}{x^{i(p-1)+\beta}y^{r-i(p-1)-\beta}} + \BIGO{p^{p-1}}
    \\ & \null\qquad\textstyle -aC_{\beta-\alpha}p^{-\alpha}\sum_{0\ne\MU\in\FF_p} [\MU]^{\alpha-\beta}\begin{Smatrix}
p & [\MU] \\ 0 & 1
\end{Smatrix} \sqbrackets{}{KZ}{\QQp}{\theta^n  x^{\alpha-n}y^{r-np-\alpha}}.
\end{align*}
If we sum this over all \m{\alpha\in\{0,\ldots,\gamma\}} and then divide by \m{p-1}, we get 
  \begin{align*}
      \textstyle & \textstyle \sum_{\alpha=0}^\gamma \sum_i C_{\beta-\alpha}\binom{r-\alpha}{i(p-1)+\beta-\alpha} \sqbrackets{}{KZ}{\QQp}{x^{i(p-1)+\beta}y^{r-i(p-1)-\beta}} + \BIGO{p^{p-1}}
      \\ & \textstyle \null\qquad 
      -\frac{a}{p-1}\sum_{\alpha=0}^\gamma C_{\beta-\alpha}p^{-\alpha}\sum_{0\ne\MU\in\FF_p} [\MU]^{\alpha-\beta}\begin{Smatrix}
p & [\MU] \\ 0 & 1
\end{Smatrix} \sqbrackets{}{KZ}{\QQp}{\theta^n  x^{\alpha-n}y^{r-np-\alpha}},
  \end{align*}
 This element is in \m{\Img(T-a)} and, after changing to the variable \m{l=\beta-\alpha}, it turns into precisely the element we want.
\end{proof}

\begin{lemma}\label{LEMMAX3}
{
Suppose that \m{\alpha \in \ZZ} and \m{v \in \QQ} and the family \m{\FAMILY{D_i}_{i\in \ZZ}} of elements of \m{\ZZ_p} are such that
\begin{align*}\textstyle \alpha & \textstyle \in\{0,\ldots,\TEMPORARYI-1\}, \\D_i&\textstyle{}=0 \text{ for } i\not\in[0,\frac{r-\alpha}{p-1}],\\ \textstyle v & \textstyle \LEQSLANT \PADICVALUATION(\TEMPORARYD_{\alpha}(D_\DARKCIRC)),
\\\min\{\PADICVALUATION(a)-\alpha,v\} & \textstyle  \LEQSLANT \PADICVALUATION(\TEMPORARYD_{w}(D_\DARKCIRC)) \text{ for } \alpha < w < 2\TEMPORARYI-\alpha,
\\ v^{\prime} & \textstyle <\PADICVALUATION(\TEMPORARYD_{w}(D_\DARKCIRC)) \text{ for } 0\LEQSLANT w<\alpha.
\end{align*} Define \m{v^{\prime} \defeq \min\{\PADICVALUATION(a)-\alpha,v\}}. Note that \m{2\TEMPORARYI-\alpha\LEQSLANT 2\cdot\frac{p-1}{2}<p}. %
 If, for \m{j \in \ZZ},
\[\textstyle\DELTAA_j \defeq (-1)^{j-1} (1-p)^{-\alpha} \binom{\alpha}{j-1} \TEMPORARYD_{\alpha}(D_\DARKCIRC),\]
so that \m{\FAMILY{\DELTAA_j}_{j\in \ZZ}} is supported on the set of indices \m{\{1,\ldots,\alpha+1\}} and therefore \m{\TEMPORARYD_{w}(\DELTAA_\DARKCIRC)} is properly defined for \m{0\LEQSLANT w < \alpha},
then \m{v \LEQSLANT\PADICVALUATION(\TEMPORARYD_{\alpha}(\DELTAA_\DARKCIRC))\LEQSLANT \PADICVALUATION(\DELTAA_j)} for all \m{j \in \ZZ}, and
\begin{align*}
      \textstyle & \textstyle\vphantom{\sum_{w=0}^\alpha \left(\VERIFYIF{(i,w)=(0,0)} + \sum_{l=1}^\alpha (-1)^{i} \binom{l}{i} \binom{l-1}{l-w}\right) \binom{-j}{w-v}} \sum_i (\DELTAA_i - D_i) \sqbrackets{}{KZ}{\QQp}{x^{i(p-1)+\alpha}y^{r-i(p-1)-\alpha}}
      \\ & \textstyle \null\vphantom{\sum_{w=0}^\alpha \left(\VERIFYIF{(i,w)=(0,0)} + \sum_{l=1}^\alpha (-1)^{i} \binom{l}{i} \binom{l-1}{l-w}\right) \binom{-j}{w-v}}\qquad 
      =  \VERIFYIF{\alpha\LEQSLANT s}(-1)^{n+1} D_{\frac{r-s}{p-1}}
 \sqbrackets{}{KZ}{\QQp}{\theta^n x^{r-np-s+\alpha}y^{s-\alpha-n}} 
            \\ & \textstyle \null\vphantom{\sum_{w=0}^\alpha \left(\VERIFYIF{(i,w)=(0,0)} + \sum_{l=1}^\alpha (-1)^{i} \binom{l}{i} \binom{l-1}{l-w}\right) \binom{-j}{w-v}}\qquad\qquad  - D_0
 \sqbrackets{}{KZ}{\QQp}{\theta^n x^{\alpha-n}y^{r-np-\alpha}}  %
  \\ & \textstyle \null\vphantom{\sum_{w=0}^\alpha \left(\VERIFYIF{(i,w)=(0,0)} + \sum_{l=1}^\alpha (-1)^{i} \binom{l}{i} \binom{l-1}{l-w}\right) \binom{-j}{w-v}}\qquad\qquad   + E \sqbrackets{}{KZ}{\QQp}{\theta^{\alpha+1}h} + %
  F \sqbrackets{}{KZ}{\QQp}{h^{\prime}} + \mathrm{ERR}_1 + \mathrm{ERR}_2,
  \end{align*}
  for some integral \m{\mathrm{ERR}_1\in \COMPACTIND{KZ}{G}\TILDE\SIGMA_{r}} and \m{\mathrm{ERR}_2\in\COMPACTIND{KZ}{G}\TILDE\SIGMA_{r}} such that
  \[\textstyle \mathrm{ERR}_1 \in\Img(T-a) \text{ and }\mathrm{ERR}_2 = \BIGO{p^{\TEMPORARYI-\PADICVALUATION(a)+v}+p^{\TEMPORARYI-\alpha}},\]
 some polynomials %
 \m{h} and \m{h^{\prime}}, and some %
  \m{E,F\in\QQp} such that %
  \m{\PADICVALUATION(E)\GEQSLANT v^{\prime}} and \m{\PADICVALUATION(F)>v^{\prime}}.
}
\end{lemma}

\begin{proof}{Proof. }
By using the equation
\begin{align}
    \label{eqeriufr3}
\textstyle g = a^{-1}T g - (T-a) (a^{-1}g)\end{align}
we can rewrite
 \[\textstyle g = \sum_i (\DELTAA_i - D_i) \sqbrackets{}{KZ}{\QQp}{x^{i(p-1)+\alpha}y^{r-i(p-1)-\alpha}}\] as
\begin{align}
       & \textstyle 
      \vphantom{\sum_{w=0}^\alpha \left(\VERIFYIF{(i,w)=(0,0)} + \sum_{l=1}^\alpha (-1)^{i} \binom{l}{i} \binom{l-1}{l-w}\right) \binom{-j}{w-v}}- \sum_{\TEMPORARYA =0}^{2\TEMPORARYI-1-\alpha} a^{-1} X_\TEMPORARYA  p^\TEMPORARYA \sum_{0 \ne \LAMBDA\in \FF_p} [-\LAMBDA]^{r-\alpha-\TEMPORARYA } \begin{Lmatrix}
      p & [\LAMBDA] \\ 0 & 1
      \end{Lmatrix} \sqbrackets{}{KZ}{\QQp}{x^{r-\TEMPORARYA }y^\TEMPORARYA }\nonumber
             \\ & \textstyle \null\vphantom{\sum_{w=0}^\alpha \left(\VERIFYIF{(i,w)=(0,0)} + \sum_{l=1}^\alpha (-1)^{i} \binom{l}{i} \binom{l-1}{l-w}\right) \binom{-j}{w-v}}\qquad - \VERIFYIF{\alpha\LEQSLANT s}(-1)^{n+1} D_{\frac{r-s}{p-1}}
 \sqbrackets{}{KZ}{\QQp}{\theta^n x^{r-np-s+\alpha}y^{s-\alpha-n}} \nonumber
            \\ & \textstyle \null\vphantom{\sum_{w=0}^\alpha \left(\VERIFYIF{(i,w)=(0,0)} + \sum_{l=1}^\alpha (-1)^{i} \binom{l}{i} \binom{l-1}{l-w}\right) \binom{-j}{w-v}}\qquad  + D_0
 \sqbrackets{}{KZ}{\QQp}{\theta^n x^{\alpha-n}y^{r-np-\alpha}}+ \BIGO{p^{\TEMPORARYI-\alpha}}\nonumber \\ & \textstyle \null\vphantom{\sum_{w=0}^\alpha \left(\VERIFYIF{(i,w)=(0,0)} + \sum_{l=1}^\alpha (-1)^{i} \binom{l}{i} \binom{l-1}{l-w}\right) \binom{-j}{w-v}}\qquad   + \mathrm{ERR_3},     \label{equation367}
  \end{align}
where \m{\mathrm{ERR_3}\in\COMPACTIND{KZ}{G}\TILDE\SIGMA_{r}}
is such that \m{\mathrm{ERR_3}  \in\Img(T-a)}, and \[\textstyle X_\TEMPORARYA  = \sum_i (\DELTAA_i-D_i)\binom{r-i(p-1)-\alpha}{\TEMPORARYA }.\]  The first three lines of equation~\eqref{equation367} come from the ``\m{a^{-1}T g}'' part of equation~\eqref{eqeriufr3}, and to obtain the second and third lines we use  part~(\ref{111}) of lemma~\ref{LEMMAX2} to replace terms of type \m{a^{-1}\begin{Smatrix}
1 & 0 \\ 0 & p
\end{Smatrix}  \sqbrackets{}{KZ}{\QQp}{h_1}} with terms of type \m{ 1 \sqbrackets{}{KZ}{\QQp}{h_2}} plus some error terms that are in \m{\Img(T-a)} and other error terms that are in \m{\BIGO{p^{n-\alpha}}}---we omit the full details.
 Note that \m{X_\TEMPORARYA} is well-defined since \m{\FAMILY{D_i}_{i\in\ZZ}} and \m{\FAMILY{\DELTAA_i}_{i\in\ZZ}} are both supported on the finite  set of indices \m{i\in\{0,\ldots,\lfloor\frac{r-\alpha}{p-1}\rfloor\}}.
The constants \m{\FAMILY{\DELTAA_i}_{i\in\ZZ}} are precisely designed in a way that \m{\TEMPORARYD_{w}(\DELTAA_\DARKCIRC)=0} for all \m{w\in\{0,\ldots,\alpha-1\}} and \m{\TEMPORARYD_{\alpha}(\DELTAA_\DARKCIRC) = \TEMPORARYD_{\alpha}(D_\DARKCIRC)}, and this follows immediately from lemma~\ref{part1X5prime}, equation~\eqref{ceqn5} applied to \m{(b,l,u)=(0,\alpha,(p-1)(\alpha+1))}. By writing \m{\binom{r-X-\alpha}{\TEMPORARYA} \in \QQ_p[X]} as a \m{\ZZ_p}-linear combination of  \m{\binom{X}{0},\ldots,\binom{X}{\TEMPORARYA}\in\QQ_p[X]} we get that
\[\textstyle{} \sum_i (\DELTAA_i-D_i)\binom{r-i(p-1)-\alpha}{\TEMPORARYA } \] is a \m{\ZZ_p}-linear combination of \m{\TEMPORARYD_{0}(\DELTAA_\DARKCIRC)-\TEMPORARYD_{0}(D_\DARKCIRC),\ldots,\TEMPORARYD_{\TEMPORARYA}(\DELTAA_\DARKCIRC)-\TEMPORARYD_{\TEMPORARYA}(D_\DARKCIRC)}. Thus we can conclude that
\[\textstyle{}\PADICVALUATION(X_\TEMPORARYA)>v^{\prime} \text{ for all }\TEMPORARYA\in\{0,\ldots,\alpha\}.\] We also 
have \[\textstyle{}\PADICVALUATION(\Delta_i) \GEQSLANT \PADICVALUATION(\TEMPORARYD_{\alpha}(D_\DARKCIRC)) \GEQSLANT v \GEQSLANT v^{\prime}\] by assumption, so we can conclude that \[\textstyle{}\textstyle \PADICVALUATION(X_\TEMPORARYA ) \GEQSLANT v^{\prime}\text{ for all }\TEMPORARYA \in \{\alpha+1,\ldots,2\TEMPORARYI-1-\alpha\}.\] Thus the part of the sum
``\m{\sum_{\TEMPORARYA =0}^{2\TEMPORARYI-1-\alpha}}'' in the first line of~\eqref{equation367} consisting of the indices  \m{\TEMPORARYA \in \{\TEMPORARYI,\ldots,2\TEMPORARYI-1-\alpha\}} is in \[\textstyle\BIGO{p^{\TEMPORARYI-\PADICVALUATION(a)+v^{\prime}}} = \BIGO{p^{\TEMPORARYI-\PADICVALUATION(a)+v}+p^{\TEMPORARYI-\alpha}}.\] Due to part~(\ref{111}) of lemma~\ref{LEMMAX2}, the part of the sum consisting of the indices \m{ \TEMPORARYA \in \{\alpha+1,\ldots,\TEMPORARYI-1\}} is
\begin{align*}
      \textstyle & \textstyle a^{-1}\sum_{\TEMPORARYA =\alpha+1}^{\TEMPORARYI-1} X_\TEMPORARYA  p^\TEMPORARYA \sum_{0 \ne \LAMBDA\in \FF_p} [-\LAMBDA]^{r-\alpha-\TEMPORARYA } \begin{Lmatrix}
      p & [\LAMBDA] \\ 0 & 1
      \end{Lmatrix} \sqbrackets{}{KZ}{\QQp}{x^{r-\TEMPORARYA }y^\TEMPORARYA }
      \\ & \textstyle \null\qquad 
    = a^{-1}\sum_{\TEMPORARYA =\alpha+1}^%
      {\TEMPORARYI-1} X_\TEMPORARYA  a \sum_{0 \ne \LAMBDA\in \FF_p} [-\LAMBDA]^{r-\alpha-\TEMPORARYA } \begin{Lmatrix}
      1 & [\LAMBDA] \\ 0 & 1
      \end{Lmatrix} \sqbrackets{}{KZ}{\QQp}{\theta^\TEMPORARYA  h_\TEMPORARYA }+ \mathrm{ERR}_4
       \\ & \textstyle \null\qquad 
      = %
      E %
      \sqbrackets{}{KZ}{\QQp}{\theta^{\alpha+1}h%
      }  + \mathrm{ERR}_4,
  \end{align*}
  where \m{\mathrm{ERR_4}\in\COMPACTIND{KZ}{G}\TILDE\SIGMA_{r}}
is such that  \[\textstyle \mathrm{ERR}_4 + \BIGO{p^{\TEMPORARYI-\PADICVALUATION(a)+v} + p^{\TEMPORARYI-\alpha}} \in \Img(T-a)\]
 and \m{h} and \m{E} are such that \m{\PADICVALUATION(E) \GEQSLANT v^{\prime}} (the constant \m{E} can be chosen in this way since \m{\PADICVALUATION(X_{\TEMPORARYA }) \GEQSLANT v^{\prime}} for all \m{ \TEMPORARYA \in \{\alpha+1,\ldots,\TEMPORARYI-1\}}). Similarly, the part of the sum consisting of the indices \m{ \TEMPORARYA \in \{0,\ldots,\alpha\}} is \[\textstyle \sum_{\TEMPORARYA =0}^{\alpha} F_\TEMPORARYA  \sqbrackets{}{KZ}{\QQp}{\theta^{\TEMPORARYA }h_\TEMPORARYA }+ \mathrm{ERR}_5, \]
 where \m{\mathrm{ERR}_5\in\COMPACTIND{KZ}{G}\TILDE\SIGMA_{r}}
is such that  \[\textstyle \mathrm{ERR}_5 + \BIGO{p^{\TEMPORARYI-\PADICVALUATION(a)+v} + p^{\TEMPORARYI-\alpha}} \in \Img(T-a)\] and \m{h_\TEMPORARYA } and \m{F_\TEMPORARYA } are such that \m{\PADICVALUATION(F_\TEMPORARYA ) \GEQSLANT \PADICVALUATION(X_\TEMPORARYA ) > v^{\prime}}. We get the identity we want after writing \[\textstyle Fh^{\prime} = \sum_{\TEMPORARYA =0}^{\alpha}F_\TEMPORARYA \theta^{\TEMPORARYA }h_\TEMPORARYA \]
 for suitable \m{h^{\prime}} and \m{F}.
\end{proof}

\begin{lemma}\label{LEMMAX4}
{
Let  \m{\FAMILY{C_l}_{l\in \ZZ}} be any family of elements of \m{\ZZ_p}. Suppose that 
\m{\alpha \in \{0,\ldots, \TEMPORARYI-1\}} and \m{v \in \QQ}, and suppose that the constants
 \[\textstyle D_i \defeq %
 \VERIFYIF{i=0} C_{-1} + \VERIFYIF{0<i(p-1)<r-2\alpha}\sum_{l=0}^\alpha C_{l} \binom{r-\alpha+l}{i(p-1)+l}\]
satisfy the conditions of lemma~\ref{LEMMAX3}, \IE 
\begin{align*}\textstyle  \textstyle v & \textstyle \LEQSLANT \PADICVALUATION(\TEMPORARYD_{\alpha}(D_\DARKCIRC)),\\v^{\prime} \defeq \min\{\PADICVALUATION(a)-\alpha,v\} & \textstyle  \LEQSLANT \PADICVALUATION(\TEMPORARYD_{w}(D_\DARKCIRC)) \text{ for } \alpha < w < 2\TEMPORARYI-\alpha,
\\ v^{\prime} & \textstyle <\PADICVALUATION(\TEMPORARYD_{w}(D_\DARKCIRC)) \text{ for } 0\LEQSLANT w<\alpha.
\end{align*} Note that \m{\FAMILY{D_i}_{i\in \ZZ}} is supported on the finite set of indices \m{\{0,\ldots,\lfloor\frac{r-2\alpha}{p-1}\rfloor\}} by definition.
Moreover, suppose that \m{C_0} is a unit. Let
\[\textstyle\TEMPORARYD^{\prime} \defeq (1-p)^{-\alpha}\TEMPORARYD_{\alpha}(D_\DARKCIRC) - C_{-1}   .\]
Suppose that \m{\PADICVALUATION(C_{-1})\GEQSLANT \PADICVALUATION(\TEMPORARYD^{\prime})}.
\begin{enumerate}
    \item \label{7111} If \m{\PADICVALUATION(\TEMPORARYD^{\prime}) \LEQSLANT v^{\prime}}
then there is some element  \m{\mathrm{gen}_1 \in \SCR I_a} that
represents a generator of \m{\FILT N_\alpha}.  \\ \item \label{7222} If \m{\PADICVALUATION(a)-\alpha < v} then there is some element \m{\mathrm{gen}_2 \in \SCR I_a}  that represents a generator of a finite-codimensional submodule of 
\[\textstyle \HECKEQ{
\alpha}\left(\COMPACTIND{KZ}{G}\QUOTIENTI(\alpha)\right) =\HECKEQ{
\alpha}\left(\FILT N_\alpha/ \COMPACTIND{KZ}{G}\SUBMODULE(\alpha)\right),\]
 where \m{\HECKEQ{
\alpha}} denotes the  endomorphism  of \m{\COMPACTIND{KZ}{G}\QUOTIENTI(\alpha)} corresponding to the double coset of \m{\begin{Smatrix}
p & 0 \\ 0 & 1
\end{Smatrix}}.
\end{enumerate}
}
\end{lemma}

\begin{proof}{Proof. }
Before proceeding to the proofs of (\ref{7111}) and (\ref{7222}), let us make some initial remarks. Due to part~(\ref{222}) of lemma~\ref{LEMMAX2}, we know that \m{\Img(T-a)} contains
  \begin{align}\label{3094rf309f}
      \textstyle \sum_{i} D_i \sqbrackets{}{KZ}{\QQp}{x^{i(p-1)+\alpha}y^{r-i(p-1)-\alpha}} + \BIGO{ap^{-\alpha}}.
\end{align} 
  The conditions imposed on the constants \m{D_i} are  designed in a way that
 \[\textstyle \sum_i D_i x^{i(p-1)+\alpha}y^{r-i(p-1)-\alpha} = \theta^\alpha h + %
 h^{\prime}\] for some polynomials \m{h,h^{\prime}} such that %
  \m{\PADICVALUATION(h^{\prime})>v^{\prime}} (that is, such that the valuation of each of the coefficients of \m{h^{\prime}} is strictly greater than \m{v^{\prime}}). In the very special case when \m{\TEMPORARYD_{w}(D_\DARKCIRC)=0} for \m{0\LEQSLANT w<\alpha}, this is because these conditions are precisely the equations needed to imply that  \begin{align*}\textstyle \sum_i D_i x^{i(p-1)+\alpha}y^{r-i(p-1)-\alpha}\end{align*} can be factored as \m{\theta^\alpha h}, due to lemma~\ref{part2X5}. The general case when \m{\TEMPORARYD_{w }(D_\DARKCIRC)>v^{\prime}} for \m{0\LEQSLANT w<\alpha} can be reduced to this special case by adding a term \m{h^{\prime}} such that \m{\PADICVALUATION(h^{\prime})>v^{\prime}}.
     Then \m{\overline\theta^\alpha \overline h} is an element of \m{N_\alpha}, and the number \m{\TEMPORARYD^{\prime}} is specifically designed in a way that \m{\overline\theta^\alpha \overline g} is precisely \m{\overline{\TEMPORARYD ^{\prime}}} times a generator of \m{N_\alpha}. If \m{\TEMPORARYD ^{\prime}} is a unit then this immediately gives us an integral element of \m{\Img(T-a)} whose reduction modulo \m{\maxideal} represents a generator of \m{\FILT N_\alpha}. Note that in general the valuation of \m{\TEMPORARYD ^{\prime}} is an integer. If that integer is strictly positive then we would like to multiply \m{\theta^\alpha h} by \m{p^{-\PADICVALUATION(\TEMPORARYD ^{\prime})}} prior to reducing modulo \m{\maxideal}. We cannot do this directly, since  almost certainly there exist some \m{D_i} whose valuation is strictly smaller than the valuation of \m{\TEMPORARYD ^{\prime}}, so after multiplying \m{\theta^\alpha h} by \m{p^{-\PADICVALUATION(\TEMPORARYD ^{\prime})}} we might not get an integral element. Thus we would like to ``smoothen'' the \m{D_i} to better constants \m{\DELTAA_i} which have much of the same qualities as the \m{D_i} (\IE satisfy the same conditions) but whose valuations are all at least as large as the valuation of \m{\TEMPORARYD ^{\prime}}. We use the constants \m{\DELTAA_i} from lemma~\ref{LEMMAX3}, and we replace \m{D_i}
      with \m{\DELTAA_i} by adding
 \begin{align*}%
 \textstyle \sum_{i} (\DELTAA_i-D_i) \sqbrackets{}{KZ}{\QQp}{x^{i(p-1)+\alpha}y^{r-i(p-1)-\alpha}}+D_0
 \sqbrackets{}{KZ}{\QQp}{\theta^n x^{\alpha-n}y^{r-np-\alpha}}\end{align*} to~\eqref{3094rf309f}. We know, directly from the definition of the constants \m{\DELTAA_i}, that 
\begin{align*}
 & \textstyle \sum_{i} \DELTAA_i \sqbrackets{}{KZ}{\QQp}{x^{i(p-1)+\alpha}y^{r-i(p-1)-\alpha}} +D_0
 \sqbrackets{}{KZ}{\QQp}{\theta^n x^{\alpha-n}y^{r-np-\alpha}}\\
 &\textstyle\null\qquad = (1-p)^{-\alpha} \TEMPORARYD_{\alpha}(D_\DARKCIRC) \sqbrackets{}{KZ}{\QQp}{\theta^\alpha x^{p-1}y^{r-\alpha(p+1)-p+1}} \\ & \textstyle{}\null\qquad\qquad+ C_{-1}  \sqbrackets{}{KZ}{\QQp}{\theta^n x^{\alpha-n}y^{r-np-\alpha}},
\end{align*}
and, for any \m{A, B\in \ZZp}, the reduction modulo \m{\maxideal} of %
\[\textstyle \textstyle A \sqbrackets{}{KZ}{\QQp}{\theta^\alpha x^{p-1}y^{r-\alpha(p+1)-p+1}} + B  \sqbrackets{}{KZ}{\QQp}{\theta^n x^{\alpha-n}y^{r-np-\alpha}}\] is \m{\overline{A-B}} times a generator of \m{\FILT N_\alpha}.
We use lemma~\ref{LEMMAX3} to deduce that if we add to \begin{align*}%
 \textstyle \sum_{i} (\DELTAA_i-D_i) \sqbrackets{}{KZ}{\QQp}{x^{i(p-1)+\alpha}y^{r-i(p-1)-\alpha}}+D_0
 \sqbrackets{}{KZ}{\QQp}{\theta^n x^{\alpha-n}y^{r-np-\alpha}}\end{align*}%
  a certain error term that looks like 
\[\textstyle{} -E \sqbrackets{}{KZ}{\QQp}{\theta^{\alpha+1}h}  %
  -F \sqbrackets{}{KZ}{\QQp}{h^{\prime}} + \mathrm{ERR}_1 + \mathrm{ERR}_2\] 
 then we get an element of \m{\Img(T-a)}, so that \m{\Img(T-a)} contains
\begin{align}& \textstyle (\TEMPORARYD^{\prime}+C_{-1}) \sqbrackets{}{KZ}{\QQp}{\theta^\alpha x^{p-1}y^{r-\alpha(p+1)-p+1}}  + C_{-1}  \sqbrackets{}{KZ}{\QQp}{\theta^n x^{\alpha-n}y^{r-np-\alpha}}\nonumber
\\ & \textstyle\qquad + %
 E\sqbrackets{}{KZ}{\QQp}{\theta^{\alpha+1}h%
 }  +  F \sqbrackets{}{KZ}{\QQp}{h^{\prime}}
+ \BIGO{p^{\TEMPORARYI-\PADICVALUATION(a)+v}+p^{\PADICVALUATION(a)-\alpha}}\label{EQIMP}%
,\end{align}
for some \m{h%
,h^{\prime},E%
,F} with \m{\PADICVALUATION(E%
) \GEQSLANT v^{\prime}} and \m{\PADICVALUATION(F)>v^{\prime}}.
Now let us proceed to the proofs of (\ref{7111}) and (\ref{7222}). %

(\ref{7111}) Suppose that \m{\PADICVALUATION(\TEMPORARYD^{\prime}) \LEQSLANT v^{\prime}}. Let us note that \m{\PADICVALUATION(C_{-1})\GEQSLANT \PADICVALUATION(\TEMPORARYD^{\prime})}.
The terms in the first line of equation~\eqref{EQIMP} have valuations at least \m{\PADICVALUATION(\TEMPORARYD^{\prime})}, and  the terms in the second line of equation~\eqref{EQIMP} have valuations at least \m{v^{\prime}\GEQSLANT \PADICVALUATION(\TEMPORARYD^{\prime})}. Thus if we  we multiply~\eqref{EQIMP} by \m{|\TEMPORARYD^{\prime}|=p^{-\PADICVALUATION(\TEMPORARYD^{\prime})}} then we get an integral element \m{\widetilde{\mathrm{gen}}_1} that is still an element of \m{\Img(T-a)}. As \m{\PADICVALUATION(E)\GEQSLANT v^{\prime}} and \m{\theta^{\alpha+1}h} is trivial in \m{N_\alpha}, the
reduction modulo \m{\maxideal} of \[\textstyle{}|\TEMPORARYD^{\prime}|E\sqbrackets{}{KZ}{\QQp}{\theta^{\alpha+1}h%
 }\] is trivial in \m{\FILT N_\alpha}.
As \m{\PADICVALUATION(F)>v^{\prime}}, the reduction modulo \m{\maxideal} of \[\textstyle{}|\TEMPORARYD^{\prime}|F \sqbrackets{}{KZ}{\QQp}{h^{\prime}}\] is trivial. And, as \m{\PADICVALUATION(a)} is not an integer, the inequality \m{\PADICVALUATION(a)-\alpha \GEQSLANT \PADICVALUATION(\TEMPORARYD^{\prime})} which follows from the assumption that \m{v^{\prime}\GEQSLANT\PADICVALUATION(\TEMPORARYD^{\prime})} must in fact be a strict inequality. Therefore  the reduction modulo \m{\maxideal} of the last term \[\textstyle{}|\TEMPORARYD^{\prime}|\BIGO{p^{\TEMPORARYI-\PADICVALUATION(a)+v}+p^{\PADICVALUATION(a)-\alpha}}\] is trivial. Thus the reduction modulo \m{\maxideal} of \m{\widetilde{\mathrm{gen}}_1} is the reduction of
\[\textstyle{} |\TEMPORARYD^{\prime}|(\TEMPORARYD^{\prime}+C_{-1}) \sqbrackets{}{KZ}{\QQp}{\theta^\alpha x^{p-1}y^{r-\alpha(p+1)-p+1}}  + |\TEMPORARYD^{\prime}|C_{-1}  \sqbrackets{}{KZ}{\QQp}{\theta^n x^{\alpha-n}y^{r-np-\alpha}},\]
 which, as mentioned, is \m{\overline{|\TEMPORARYD^{\prime}|(\TEMPORARYD^{\prime}+C_{-1})-|\TEMPORARYD^{\prime}|C_{-1} }\in \FFptimes} times a generator of \m{\FILT N_\alpha}. %
  We define  the desired element as the reduction modulo \m{\maxideal} of \m{\widetilde{\mathrm{gen}}_1}.

(\ref{7222}) Suppose that \m{\PADICVALUATION(a)-\alpha<v}.  We may without loss of generality assume that \m{\PADICVALUATION(\TEMPORARYD^{\prime}) > v^{\prime} = \PADICVALUATION(a)-\alpha} since otherwise we could apply the first part of this lemma and reach an even stronger conclusion. Then %
 the dominant term in~\eqref{EQIMP} %
  comes from %
 the error term \m{\BIGO{ap^{-\alpha}}}. This is because the terms in the first line of equation~\eqref{EQIMP}
  have valuation at least \m{\PADICVALUATION(\TEMPORARYD^{\prime}) > \PADICVALUATION(a)-\alpha}, the ``\m{F}'' and ``\m{\BIGO{p^{\TEMPORARYI-\PADICVALUATION(a)+v}}}'' terms in the second line of equation~\eqref{EQIMP}
  have valuation bigger than \m{v^{\prime}=\PADICVALUATION(a)-\alpha}, and  \m{\PADICVALUATION(E)\GEQSLANT v^{\prime}=\PADICVALUATION(a)-\alpha} and \m{\theta^{\alpha+1}h} is trivial in \m{N_\alpha} so the reduction modulo  \m{\maxideal} of 
  \[\textstyle{}p^{\alpha-\PADICVALUATION(a)}E\sqbrackets{}{KZ}{\QQp}{\theta^{\alpha+1}h}\]
 is trivial in \m{\FILT N_\alpha}.
The term \m{\BIGO{ap^{-\alpha}}} %
  is described in the last displayed equation in the statement of  lemma~\ref{LEMMAX2}, and it is
    \[\textstyle -\frac{ap^{-\alpha}C_0}{p-1}  \sum_{0\ne\LAMBDA\in\FF_p} \begin{Lmatrix}
p & [\LAMBDA] \\ 0 & 1
\end{Lmatrix}  \sqbrackets{}{KZ}{\QQp}{\theta^n  x^{\alpha-n}y^{r-np-\alpha}}.\]
Thus after we multiply the element given by equation~\eqref{EQIMP} with \m{\frac{p-1}{ap^{-\alpha}C_0}} and reduce the resulting element modulo \m{\maxideal}, we get a representative of
\[\textstyle \left(\sum_{\LAMBDA\in\FF_p}\begin{Lmatrix}
p & [\LAMBDA] \\ 0 & 1
\end{Lmatrix} + A \begin{Lmatrix}
p & 0 \\0 & 1
\end{Lmatrix} +\VERIFYIF{r\equiv_{p-1}2\alpha} B \begin{Lmatrix}
1 & 0 \\0 & p
\end{Lmatrix}\right)\sqbrackets{}{KZ}{\FFp}{X^{\LESS{2\alpha-r}}},\] for some \m{A,B \in \FFp}. Let \m{\gamma} denote this element, which belongs to \[\COMPACTIND{KZ}{G}\QUOTIENTI(\alpha) = \FILT N_\alpha/\COMPACTIND{KZ}{G}\SUBMODULE(\alpha) \cong \COMPACTIND{KZ}{G}(\textstyle \sigma_{\LESS{2\alpha-r}}(r-\alpha))\] and is represented by the reduction modulo \m{\maxideal} of an integral element of \m{\Img(T-a)} (and therefore by an element of \m{\SCR I_a}). 
The classification given in theorem~\ref{BERGERX1} implies that either \m{\HECKEQ{
\alpha}-\lambda} acts trivially on \m{\COMPACTIND{KZ}{G}\QUOTIENTI(\alpha)}
modulo \m{\SCR I_{a}} for some \m{\lambda \in \FFp}, or \m{s\in\{1,3,\ldots,2\TEMPORARYI-1\}} and 
\m{(\HECKEQ{
\alpha}-\lambda)(\HECKEQ{
\alpha}-\lambda^{-1})} acts trivially on \m{\COMPACTIND{KZ}{G}\QUOTIENTI(\alpha)}
modulo \m{\SCR I_{a}} for some  \m{\lambda \in \overline {\FF}_p^\times}. %
\begin{itemize}[leftmargin=*]
    \item Let us first consider the case when \m{\HECKEQ{
\alpha}-\lambda} acts trivially.
Then
\begin{align}\label{eq4}%
\textstyle \left(A\begin{Lmatrix}
p & 0 \\ 0 & 1
\end{Lmatrix} + \lambda - \VERIFYIF{r\equiv_{p-1}2\alpha}(1-B)\begin{Lmatrix}
1 & 0 \\ 0 & p
\end{Lmatrix} \right) \sqbrackets{}{KZ}{\FFp}{X^{\LESS{2\alpha-r}}}\end{align}
 in \m{\COMPACTIND{KZ}{G}\QUOTIENTI(\alpha)} is represented by an element of \m{\SCR I_a}.
First suppose that \m{\LESS{2\alpha-r}>0}. If \m{A\ne 0} then either \m{\lambda = 0} in which case \m{\COMPACTIND{KZ}{G}\QUOTIENTI(\alpha)} is trivial modulo \m{\SCR I_a}, or \m{\lambda \ne 0} in which case we can multiply~\eqref{eq4} on the left by \m{[\MU]^{p-2}\begin{Smatrix}
1 & 0 \\ [\MU] & 1 
\end{Smatrix}} and sum over all \m{\MU \in\FF_p} to obtain that \m{\lambda \sqbrackets{}{KZ}{\FFp}{X^{\LESS{2\alpha-r}-1}Y}} is represented by an element of  \m{\SCR I_a}, so we can take \[\textstyle \mathrm{gen}_2 = \lambda \sqbrackets{}{KZ}{\FFp}{X^{\LESS{2\alpha-r}-1}Y}.\]
 If \m{A = 0} then either \m{\lambda\ne 0} in which case \m{\COMPACTIND{KZ}{G}\QUOTIENTI(\alpha)} is trivial modulo \m{\SCR I_a}, or \m{\lambda=0} in which case %
  we can take \[\textstyle\mathrm{gen}_2 = \HECKEQ{
\alpha}\left(1\sqbrackets{}{KZ}{\FFp}{X^{\LESS{2\alpha-r}}}\right).\] Now suppose that \m{\LESS{2\alpha-r}=0}. Then we can use the decomposition of \m{G} into cosets of \m{KZ} given in section~2.1.2 of~\cite{Bre03b} together with the fact that the element given by equation~\eqref{eq4} is trivial modulo \m{\SCR I_a} to conclude that any element of \m{\COMPACTIND{KZ}{G}\QUOTIENTI(\alpha)}
  can be written in the form
  \[\textstyle \left(\mu_1\begin{Lmatrix}
p & 0 \\ 0 & 1
\end{Lmatrix} + \mu_2 \right) \sqbrackets{}{KZ}{\FFp}{1} + h^{\prime\prime}\] for some \m{h^{\prime\prime}} which is represented by an element of \m{\SCR I_a}, and thus we can find a suitable \m{\mathrm{gen}_2 \in \SCR I_a} that represents a generator of a submodule of 
\m{\HECKEQ{
\alpha}(\COMPACTIND{KZ}{G}\QUOTIENTI(\alpha))} which has codimension at most two.
\item  Now let us consider the case when \m{(\HECKEQ{
\alpha}-\lambda)(\HECKEQ{
\alpha}-\lambda^{-1})} acts trivially. Thus \m{\lambda \ne 0} and \m{s \in \{1,3,\ldots,2\TEMPORARYI-1\}}, and in particular \m{\LESS{2\alpha-r}>0}. As \begin{align*}%
\textstyle \left(\sum_{\LAMBDA\in\FF_p}\begin{Lmatrix}
p & [\LAMBDA] \\ 0 & 1
\end{Lmatrix} + A\begin{Lmatrix}
p & 0 \\ 0 & 1
\end{Lmatrix}  \right) \sqbrackets{}{KZ}{\FFp}{X^{\LESS{2\alpha-r}}}\end{align*}
is trivial in \m{\COMPACTIND{KZ}{G}\QUOTIENTI(\alpha)} modulo \m{\SCR I_a}, it follows that so is
\begin{align*}%
\textstyle \left(A^2\begin{Lmatrix}
p^2 & 0 \\ 0 & 1
\end{Lmatrix} + (\lambda+\lambda^{-1})A\begin{Lmatrix}
p & 0 \\ 0 & 1
\end{Lmatrix} + 1 \right) \sqbrackets{}{KZ}{\FFp}{X^{\LESS{2\alpha-r}}}.\end{align*}
Again  after multiplying on the left by \m{[\MU]^{p-2}\begin{Smatrix}
1 & 0 \\ [\MU] & 1 
\end{Smatrix}} and summing over all \m{\MU \in\FF_p} we conclude that \m{1 \sqbrackets{}{KZ}{\FFp}{X^{\LESS{2\alpha-r}-1}Y}} is trivial (since \m{\LESS{2\alpha-r}>0}) and therefore that \m{\COMPACTIND{KZ}{G}\QUOTIENTI(\alpha)} is  trivial modulo \m{\SCR I_a}. %
\end{itemize}
 Consequently, either \m{\SCR I_a} contains a representative of a generator of a finite-codimensional submodule of \m{\HECKEQ{
\alpha}(\COMPACTIND{KZ}{G}\QUOTIENTI(\alpha))}, or it contains a representative of a generator of \m{\COMPACTIND{KZ}{G}\QUOTIENTI(\alpha)} (the latter being a stronger statement), and we can take \m{\mathrm{gen}_2} to be that element of \m{\SCR I_a}.
\end{proof}

\section{\texorpdfstring{Proof of theorem~\ref{MAINTHM1}}{Proof of theorem~1}}
\label{secProof}

The proof is based on the approach outlined in~\cite{BG09},  and  roughly consists of finding enough elements in \m{\SCR I_{a}}, consequently eliminating enough subquotients of \m{\COMPACTIND{KZ}{G}\SIGMA_r}, and using that information to find enough data about \m{\overline\Theta_{k,a}}.

Throughout  sections~\ref{secProof} and~\ref{secProof2} we assume that \[\textstyle r = s + \beta(p-1) + u_0p^t + \BIGO{p^{t+1}}\] for some \m{\beta \in \{0,\ldots,p-1\}}  and \m{u_0 \in \ZZ_p^\times} and \m{t \in \ZZ_{\GEQSLANT1}}. Let us write \m{\TEMPORARYJ = u_0p^t}. Recall also that we assume 
\m{\TEMPORARYI-1<\PADICVALUATION(a)<\TEMPORARYI} for some  \m{\TEMPORARYI\in\{1,\ldots,\frac{p-1}{2}\}}, that \m{s\in\{2\TEMPORARYI,\ldots,p-1\}}, and that \m{k>p^{100}} (and consequently \m{r > p^{99}}).

Recall that theorem~\ref{MAINTHM1} is equivalent to  theorem~\ref{MAINTHM1aa}.
Let us show that theorem~\ref{MAINTHM1aa} follows from the following proposition.

\begin{proposition}\label{mainPROP}
For each \m{\alpha \in \{0,\ldots,\TEMPORARYI-2\}},
there is some element \m{\mathrm{gen} \in \SCR I_a}  that represents a generator of a finite-codimensional submodule of 
\[\textstyle \HECKEQ{
\alpha}\left(\COMPACTIND{KZ}{G}\QUOTIENTI(\alpha)\right) =\HECKEQ{
\alpha}\left(\FILT N_\alpha/ \COMPACTIND{KZ}{G}\SUBMODULE(\alpha)\right).\]
\end{proposition}

\begin{proof}{Proof that proposition~\ref{mainPROP} implies theorem~\ref{MAINTHM1aa}. }
The classifications given by theorem~\ref{BERGERX1} and theorem~2.7.1 in~\cite{Bre03a} imply that if \m{\overline \Theta_{k,a}} is reducible then it must have exactly two infinite-dimensional factors. There may be an additional one-dimensional factor, a twist of the Steinberg representation. Thus suppose that \m{\overline \Theta_{k,a}} is reducible and its infinite-dimensional factors are \m{R_1} and \m{R_2}. Recall from section~\ref{secNot} that \m{\overline{\Theta}_{k,a}} is a subquotient of \[\textstyle \COMPACTIND{KZ}{G}(\Sigma_r/\langle y^r,\ldots,x^{\TEMPORARYI-1}y^{r-\TEMPORARYI+1},\overline{\theta}^{\TEMPORARYI}\Sigma_{r-\TEMPORARYI(p+1)}\rangle),\]
 a module which has a series whose factors are subquotients of \m{\FILT N_0,\ldots,\FILT N_{\TEMPORARYI-1}}. Thus
each of \m{R_1,R_2} is a quotient of one of the representations
\[\textstyle{} \COMPACTIND{KZ}{G}\SUBMODULE(0),\ldots,\COMPACTIND{KZ}{G}\SUBMODULE(\TEMPORARYI-1),\COMPACTIND{KZ}{G}\QUOTIENTI(0), \ldots,\COMPACTIND{KZ}{G}\QUOTIENTI(\TEMPORARYI-1).\]
If \m{R_j} happened to be a quotient of \m{\QUOTIENTI(\alpha)} for some \m{\alpha \in \{0,\ldots,\TEMPORARYI-2\}} and some \m{j \in \{1,2\}}, then proposition~\ref{mainPROP} would imply that \m{R_j} must be a quotient of a representation that contains \m{\COMPACTIND{KZ}{G}\QUOTIENTI(\alpha)/\HECKEQ{
\alpha}} as a finite-codimensional submodule. Since \m{\COMPACTIND{KZ}{G}\QUOTIENTI(\alpha)/\HECKEQ{
\alpha}} is irreducible (and  infinite-dimensional) and \m{R_j} is infinite-dimensional, it would follow in light of the classification given by theorem~\ref{BERGERX1} that \m{R_j \cong \COMPACTIND{KZ}{G}\QUOTIENTI(\alpha)/\HECKEQ{
\alpha}}, resulting in a contradiction. Thus neither \m{R_1} nor \m{R_2} can happen to be a quotient of \m{\QUOTIENTI(\alpha)} for some \m{\alpha \in \{0,\ldots,\TEMPORARYI-2\}}. Moreover, neither \m{R_1} nor \m{R_2} can happen to be a quotient of \m{\COMPACTIND{KZ}{G}\SUBMODULE(0)}, since \m{\SCR I_a} contains \m{1\sqbrackets{}{KZ}{\QQp}{y^{r}}}, which generates \m{\COMPACTIND{KZ}{G}\SUBMODULE(0)}. Hence each of them is a quotient of one of the representations
\[\textstyle{} \COMPACTIND{KZ}{G}\SUBMODULE(1),\ldots,\COMPACTIND{KZ}{G}\SUBMODULE(\TEMPORARYI-1),\COMPACTIND{KZ}{G}\QUOTIENTI(\TEMPORARYI-1).\]
Suppose that \m{R_1} and \m{R_2}  are quotients of \m{\COMPACTIND{KZ}{G} (\sigma_b(b^{\prime}))} and \m{\COMPACTIND{KZ}{G} (\sigma_d(d^{\prime}))}, respectively. 
By theorem~\ref{BERGERX1} we must have
\begin{align*}\textstyle b^{\prime}-d^{\prime}& \textstyle{}\equiv d+1 \bmod{p-1} , \\  \textstyle{} b+d&\textstyle{}\equiv-2 \bmod{p-1}.\end{align*}
In particular, \m{\COMPACTIND{KZ}{G} (\sigma_b(b^{\prime}))} and \m{\COMPACTIND{KZ}{G} (\sigma_d(d^{\prime}))} cannot happen to be the representations \m{\COMPACTIND{KZ}{G}\SUBMODULE(\alpha)} and \m{\COMPACTIND{KZ}{G}\SUBMODULE(\alpha^{\prime})}, as that would imply that \[\textstyle \{2,\ldots,2\TEMPORARYI-2\} \ni \alpha+\alpha^{\prime}\equiv_{p-1}r+1\equiv_{p-1}s+1 \not\in \{2,\ldots,2\TEMPORARYI\},\]
resulting in a contradiction.
Similarly,  \m{\COMPACTIND{KZ}{G} (\sigma_b(b^{\prime}))} and \m{\COMPACTIND{KZ}{G} (\sigma_d(d^{\prime}))} cannot happen to be the representations \m{\COMPACTIND{KZ}{G}\SUBMODULE(\alpha)} and \m{\COMPACTIND{KZ}{G}\QUOTIENTI(\TEMPORARYI-1)}, as that would imply that \[\textstyle{} \alpha \equiv \TEMPORARYI\bmod{p-1},\] resulting in a contradiction. And finally, \m{\COMPACTIND{KZ}{G} (\sigma_b(b^{\prime}))} and \m{\COMPACTIND{KZ}{G} (\sigma_d(d^{\prime}))} cannot happen to both be the representation \m{\COMPACTIND{KZ}{G}\QUOTIENTI(\TEMPORARYI-1)},  as that would imply that \[\textstyle{} s \equiv 2\TEMPORARYI-1\bmod{p-1},\] resulting in a contradiction. Therefore we reach a contradiction in all cases, implying that \m{\overline \Theta_{k,a}} must in fact be irreducible.
\end{proof}

The next section is dedicated to proving proposition~\ref{mainPROP}.

\section{\texorpdfstring{Proof of proposition~\ref{mainPROP}}{Proof of proposition~15}}
\label{secProof2}

Recall from section~\ref{secList} that \[\textstyle{}\FORMX_n = \FORMX(\FORMX-1)\cdots(\FORMX-n+1) \in \ZZ_p[X]\] is the falling factorial. Let \m{\alpha \in \{0,\ldots,\TEMPORARYI-2\}}. The goal is to apply lemma~\ref{LEMMAX4} with a certain  family  \m{\FAMILY{C_l}_{l\in \ZZ}} of elements of \m{\ZZ_p}  and certain \m{v \in \QQ}. For any such \m{\FAMILY{C_l}_{l\in \ZZ}} and \m{v \in \QQ} let us define  \[\textstyle D_i \defeq 
 \VERIFYIF{i=0} C_{-1} + \VERIFYIF{0<i(p-1)<r-2\alpha}\sum_{l=0}^\alpha C_{l} \binom{r-\alpha+l}{i(p-1)+l}.\] Thus throughout the proof the constants \m{\FAMILY{D_i}_{i\in\ZZ}} depend on and are defined entirely by the constants \m{(C_{-1},C_0,\ldots,C_{\alpha})}. Let us note from the definition of \m{\TEMPORARYD_j(D_\DARKCIRC)} for \m{j \in \{0,\ldots,p-1\}} that it is a \m{\ZZ_p}-linear combination of the constants \m{(C_{-1},C_0,\ldots,C_{\alpha})}. We always pick \m{C_j=0} for \m{j\not\in\{-1,0,\ldots,\alpha\}}. In the proof we make  choices for \m{(C_{-1},C_0,\ldots,C_{\alpha})} and \m{v} so that the constants \m{\FAMILY{D_i}_{i\in \ZZ}} satisfy the conditions of lemma~\ref{LEMMAX4}, \IE 
\begin{align}\textstyle  \textstyle v & \textstyle \LEQSLANT \PADICVALUATION(\TEMPORARYD_{\alpha}(D_\DARKCIRC)),\label{cond1}\tag{\taga}\\v^{\prime} \defeq \min\{\PADICVALUATION(a)-\alpha,v\} & \textstyle  \LEQSLANT \PADICVALUATION(\TEMPORARYD_{w}(D_\DARKCIRC)) \text{ for } \alpha < w < 2\TEMPORARYI-\alpha\label{cond2}\tag{\tagb},
\\ v^{\prime} & \textstyle <\PADICVALUATION(\TEMPORARYD_{w}(D_\DARKCIRC)) \text{ for } 0\LEQSLANT w<\alpha,\label{cond3}\tag{\tagc}
\\ \textstyle C_0 & \in \ZZ_p^\times,\label{cond4}\tag{\tagd}
\\ \textstyle \PADICVALUATION(C_{-1}) & \textstyle\GEQSLANT \PADICVALUATION(\TEMPORARYD^{\prime}).\label{cond5}\tag{\tage}
\end{align} We call \m{(C_{-1},C_0,\ldots,C_{\alpha})} and \m{v \in \QQ} ``suitable'' if the five conditions~\eqref{cond1}, \eqref{cond2}, \eqref{cond3}, \eqref{cond4}, \eqref{cond5} are satisfied and if moreover either \m{\PADICVALUATION(\TEMPORARYD^{\prime})\LEQSLANT v^{\prime}} or \m{\PADICVALUATION(a)-\alpha<v}.
Suppose that we are indeed  able to find suitable \m{(C_{-1},C_0,\ldots,C_{\alpha})} and \m{v}.
If \m{\PADICVALUATION(a)-\alpha<v} then we can conclude from part~(\ref{7222}) of lemma~\ref{LEMMAX4} that there is some element \m{\mathrm{gen} \in \SCR I_a}  that represents a generator of a finite-codimensional submodule of
\[\textstyle \HECKEQ{
\alpha}\left(\COMPACTIND{KZ}{G}\QUOTIENTI(\alpha)\right) =\HECKEQ{
\alpha}\left(\FILT N_\alpha/ \COMPACTIND{KZ}{G}\SUBMODULE(\alpha)\right),\] and if \m{\PADICVALUATION(\TEMPORARYD^{\prime})\LEQSLANT v^{\prime}} we  can use part~(\ref{7111}) of lemma~\ref{LEMMAX4} to reach an even stronger conclusion, that  there is some element  \m{\mathrm{gen} \in \SCR I_a} that
represents a generator of \m{\FILT N_\alpha} (and therefore \m{\mathrm{gen}} also represents a generator of 
\m{\textstyle \HECKEQ{
\alpha}(\COMPACTIND{KZ}{G}\QUOTIENTI(\alpha))}). %
 Thus all we need to do to prove proposition~\ref{mainPROP} is to find  suitable \m{(C_{-1},C_0,\ldots,C_\alpha)} and \m{v}, for  \m{\alpha \in \{0,\ldots,\TEMPORARYI-2\}}.

Let \m{\UUU} denote the \m{(2\TEMPORARYI-\alpha)\times (\alpha+1)} matrix with entries in \m{\QQ_p} indexed by \m{0\LEQSLANT w<2\TEMPORARYI-\alpha} and \m{0\LEQSLANT j\LEQSLANT\alpha}, defined by \[\textstyle \UUU (C_0,C_1,\ldots,C_\alpha)^T = (\TEMPORARYD_0(D_\DARKCIRC)-C_{-1},\TEMPORARYD_1(D_\DARKCIRC),\ldots,\TEMPORARYD_{2\TEMPORARYI-\alpha-1}(D_\DARKCIRC))^T.\] Note that only the constants \m{C_0,\ldots,C_\alpha} appear in the linear combinations on the right side, so this equation does indeed define a linear map.
Let \m{\UUU^{\subm}} denote the \m{(\alpha+1) \times (\alpha+1)} submatrix consisting of the top \m{\alpha+1} rows of \m{\UUU}, \IE of those entries indexed by  \m{0\LEQSLANT w,j\LEQSLANT\alpha}.\footnote{Note that  \m{\alpha \in \{0,\ldots,\TEMPORARYI-2\}} implies that \m{\alpha+1<2\TEMPORARYI-\alpha}.} Due to lemma~\ref{part1X5}, equation~\eqref{ceqn6},
\begin{align*}
     \UUU_{w,j} & = \textstyle \sum_v  \binom{-j}{w-v} \binom{r-\alpha+j}{v} \left(\binom{s-\alpha+j-v}{j-v}-\binom{r-\alpha+j-v}{j-v}\right) +\BIGO{p}
\end{align*}
for \m{0\LEQSLANT w < 2\TEMPORARYI-\alpha} and \m{0\LEQSLANT j\LEQSLANT\alpha}, so the entries of \m{\UUU} are in fact integers.

Let us consider four different cases: \begin{itemize}
    \item \m{\alpha=0},
    \\[-11pt] \item \m{\alpha>0} and \m{\beta \not\in \{0,\ldots,\alpha\}},
    \\[-11pt] \item \m{\alpha>0} and \m{\beta=0}, 
    \\[-11pt] \item \m{\alpha>0} and \m{\beta\in\{1,\ldots,\alpha\}}.
\end{itemize}
In order to prove proposition~\ref{mainPROP}, we must find constants \m{(C_{-1},C_0,\ldots,C_\alpha)} and \m{v\in\QQ_p} that are ``suitable'' (\IE such that the five conditions~\eqref{cond1}, \eqref{cond2}, \eqref{cond3}, \eqref{cond4}, \eqref{cond5} are satisfied and either \m{\PADICVALUATION(\TEMPORARYD^{\prime})\LEQSLANT v^{\prime}} or \m{\PADICVALUATION(a)-\alpha<v}) in each case. %

\begin{itemize}[leftmargin=*]
    \item \emph{Suppose that \m{{\alpha=0}}.} In this case we choose \m{(C_{-1},C_0) = (0,1)}. Then the two conditions~\eqref{cond1}, \eqref{cond2} are equivalent to
    \begin{align*}\textstyle  \textstyle v & \textstyle \LEQSLANT \PADICVALUATION(\TEMPORARYD_{0}(D_\DARKCIRC)),\\v^{\prime} \defeq \min\{\PADICVALUATION(a),v\} & \textstyle  \LEQSLANT \PADICVALUATION(\TEMPORARYD_{w}(D_\DARKCIRC)) \text{ for } 0 < w < 2\TEMPORARYI,
\end{align*} and the three conditions  \eqref{cond3}, \eqref{cond4}, \eqref{cond5} are vacuously true. Moreover, we set \m{v = \PADICVALUATION(\TEMPORARYD_{0}(D_\DARKCIRC))}, so that condition~\eqref{cond1} is also satisfied and we trivially have either \m{\PADICVALUATION(\TEMPORARYD^{\prime}) = v \LEQSLANT v^{\prime}} or \m{\PADICVALUATION(a)<v}. Thus all we need to check is that \m{v^{\prime} \LEQSLANT \PADICVALUATION(\TEMPORARYD_{w}(D_\DARKCIRC)) } for \m{0 < w < 2\TEMPORARYI}. %
 Due to lemma~\ref{part1X5}, equation~\eqref{ceqn2},
\[\textstyle \TEMPORARYD_0(D_\DARKCIRC) = \frac{s-r}{s}p + \BIGO{(s-r)p^2}.\]
Moreover, due to lemma~\ref{part1X5}, equation~\eqref{ceqn2} and lemma~\ref{part1X5}, equation~\eqref{ceqn3},
    \begin{align}
        \TEMPORARYD_w(D_\DARKCIRC) & =\textstyle \sum_{0<i(p-1)<r} \binom{r}{i(p-1)} \binom{i(p-1)}{w} \nonumber\\ &\textstyle{} = \sum_{0<i(p-1)<r} \binom{r-w}{i(p-1)-w} \binom{r}{w}\nonumber \\
        & = \textstyle \binom{r}{w} \sum_{0<i(p-1)<r} \sum_{u=0}^w (-1)^u \binom{w}{u} \binom{r-u}{i(p-1)} \nonumber\\
        & = \textstyle \binom{r}{w} (s-r)p \sum_{u=0}^w (-1)^u \binom{w}{u} \frac{1}{s-u} + \BIGO{(s-r)p^2}\nonumber
        \\ & = \textstyle  \frac{(-1)^w(s-r)r_w}{s_{w+1}}p + \BIGO{(s-r)p^2} = \BIGO{(s-r)p}\label{f398j98}
\end{align}
for \m{0<w<2\TEMPORARYI}. The first equality follows from the definition of \m{ \TEMPORARYD_w(D_\DARKCIRC)}, the second equality is a simple rewrite, the third equality follows from lemma~\ref{part1X5}, equation~\eqref{ceqn3}, the fourth equality follows from lemma~\ref{part1X5}, equation~\eqref{ceqn2}, and the fifth equality follows from the equation \[\textstyle{} \sum_{u=0}^w (-1)^u \binom{w}{u} \frac{1}{\FORMX} = \frac{(-1)^w w!}{\FORMX_{w+1}} \in \QQ(\FORMX)\]
and the fact that \m{s_{w+1} \ne 0} (because \m{s\GEQSLANT 2\TEMPORARYI}).
If \m{0=\alpha <\max\{\TEMPORARYI-\jj-1,\beta\}}, then either \m{\beta>0} in which case \m{s-r \in \ZZ_p^\times} and \m{\TEMPORARYD_0(D_\DARKCIRC) \in \ZZ_p^\times} and therefore the constants \m{(C_{-1},C_0)} are suitable for \m{v=0}, or \m{\beta = 0} and \m{\jj \LEQSLANT \TEMPORARYI-2} in which case \m{\PADICVALUATION(\TEMPORARYD_0(D_\DARKCIRC) ) = \jj+1} and \m{\PADICVALUATION(\TEMPORARYD_w(D_\DARKCIRC) ) \GEQSLANT \jj+1} for   \m{0<w<2\TEMPORARYI} and therefore the constants \m{(C_{-1},C_0)} are suitable for \m{v=v^{\prime}=\jj+1}. If, on the other hand, \m{\max\{\TEMPORARYI-\jj-1,\beta\}=0}, then \m{\beta=0} and \m{\jj\GEQSLANT\TEMPORARYI-1} and \m{v = \jj+1} and \m{v^{\prime}=\PADICVALUATION(a)<v\LEQSLANT \PADICVALUATION(\TEMPORARYD_{w}(D_\DARKCIRC))} for \m{0 < w < 2\TEMPORARYI}, so again the constants \m{(C_{-1},C_0)} are suitable for \m{v=\jj+1}.

\item \emph{Suppose that \m{\alpha>0} and \m{{\beta \not\in \{0,\ldots,\alpha\}}}.}  The second condition implies that \m{{\PADICVALUATION((s-r)_{\alpha+1})=0}}. Let \m{M} and \m{c} be the matrix and vector that are defined in lemma~\ref{part4X5}. The entries of \m{M} and \m{c} are in \m{\QQ(X,Y)\subset \QQ_p(X,Y)}. Let us make the substitutions \m{X=r-\alpha} and \m{Y=s-\alpha}, so that we end up with the matrix \m{\EXTEND{M}} and vector \m{\EXTEND{v}} with entries in \m{\QQ_p}.
In this case we choose
 \m{(C_{-1},C_0)=(0,1)} and 
\[\textstyle C_j = \frac{\EXTEND{c}_jp}{(s-\alpha)_{\alpha}}\] for \m{0<j\LEQSLANT\alpha}. Recall that
    \[\textstyle \UUU (1,C_1,\ldots,C_\alpha)^T = (\TEMPORARYD_0(D_\DARKCIRC),\ldots,\TEMPORARYD_{2\TEMPORARYI-\alpha-1}(D_\DARKCIRC))^T.\]
Due to lemma~\ref{part1X5}, equation~\eqref{ceqn2} and lemma~\ref{part1X5}, equation~\eqref{ceqn6},
     \begin{align*}
        \UUU_{w,0} & =\textstyle \frac{(-1)^w(s-r)(r-\alpha)_w}{(s-\alpha)_{w+1}}p + \BIGO{p^2}, 
        \\ \UUU_{w,j} & = \textstyle \sum_v  \binom{-j}{w-v} \binom{r-\alpha+j}{v} \left(\binom{s-\alpha+j-v}{j-v}-\binom{r-\alpha+j-v}{j-v}\right) +\BIGO{p},
\end{align*}
for \m{0\LEQSLANT w<2\TEMPORARYI-\alpha} and \m{0<j\LEQSLANT\alpha}. The first line follows as in equation~\eqref{f398j98}.
In particular,
    \[\textstyle \UUU_{w,j} = p^{\VERIFYIF{j=0}}(\EXTEND{M}_{w,j}+\BIGO{p})\] for \m{0\LEQSLANT w<2\TEMPORARYI-\alpha} and \m{0<j\LEQSLANT\alpha}.
 Therefore, due to lemma~\ref{part4X5}, the first \m{\alpha} entries of the resulting column vector 
 \m{(\TEMPORARYD_0(D_\DARKCIRC),\ldots,\TEMPORARYD_{2\TEMPORARYI-\alpha-1}(D_\DARKCIRC))^T}
 are \m{\BIGO{p^2}}, and the entry indexed \m{\alpha} is \[\textstyle \TEMPORARYD_{\alpha}(D_\DARKCIRC) = \frac{(s-r)_{\alpha+1}}{(s-\alpha)_{\alpha+1}}p+\BIGO{p^2} \in p\ZZ_p^\times.\] Since all subsequent entries are evidently \m{\BIGO{p}}, we can conclude that the constants \m{(C_{-1},C_0,\ldots,C_{\alpha})} are suitable for \m{v=1} (we omit the fine details of checking that all of the conditions are satisfied).
 \item \emph{Suppose that \m{\alpha>0} and \m{\beta =0}.} Let us  define \m{\EXTEND {M}} and \m{\EXTEND{c}} as in the previous case, and let
 \m{(C_{-1},C_0)=(0,1)} and 
\[\textstyle C_j = \frac{\EXTEND{c}_jp}{(s-\alpha)_{\alpha}}\] for \m{0<j\LEQSLANT\alpha}. We can similarly compute
  \begin{align*}
        \UUU_{w,0} & =\textstyle \frac{(-1)^w(s-r)(r-\alpha)_w}{(s-\alpha)_{w+1}}p + \BIGO{p^{t+2}}, 
        \\ \UUU_{w,j} & = \textstyle \sum_v  \binom{-j}{w-v} \binom{r-\alpha+j}{v} \left(\binom{s-\alpha+j-v}{j-v}-\binom{r-\alpha+j-v}{j-v}\right) +\BIGO{p^{t+1}},
\end{align*}
for \m{0\LEQSLANT w<2\TEMPORARYI-\alpha} and \m{0<j\LEQSLANT\alpha}. Again,
\[\textstyle{} \UUU (1,C_1,\ldots,C_\alpha)^T =  (\TEMPORARYD_0(D_\DARKCIRC),\ldots,\TEMPORARYD_{2\TEMPORARYI-\alpha-1}(D_\DARKCIRC))^T.\]
The first \m{\alpha} entries of this column vector are in \m{\BIGO{p^{t+2}}}, and the entry indexed \m{\alpha} is \[\textstyle \frac{(s-r)_{\alpha+1}}{(s-\alpha)_{\alpha+1}}p+\BIGO{p^{t+2}}\in p^{t+1}\ZZ_p^\times.\]
    Moreover, as long as \m{w<2\TEMPORARYI-\alpha}, we have %
    \[\textstyle \UUU_{w,0}=\frac{(-1)^w(s-r)(r-\alpha)_w}{(s-\alpha)_{w+1}}p + \BIGO{p^{t+2}} = \BIGO{p^{t+1}} ,\]
and, since \m{r = s + \BIGO{p^t}}, we have
\[\textstyle{}\textstyle{} \binom{r-\alpha+l}{l} = \binom{s-\alpha+l}{l} + \BIGO{p^t}\]
when \m{l <p}, implying directly from the formula for \m{\UUU_{w,j}} that \m{\UUU_{w,j}=\BIGO{p^t}}
 for \m{0<w<2\TEMPORARYI-\alpha} and \m{j>0}. Since \m{\PADICVALUATION(C_j)\GEQSLANT 1} for \m{1\LEQSLANT j\LEQSLANT\alpha}, it follows that the last \m{2\TEMPORARYI-2\alpha-1} entries of \m{(\TEMPORARYD_0(D_\DARKCIRC),\ldots,\TEMPORARYD_{2\TEMPORARYI-\alpha-1}(D_\DARKCIRC))^T} are in \m{\BIGO{p^{t+1}}}. Therefore the constants \m{(C_{-1},C_0,\ldots,C_{\alpha})} are suitable for \m{v=t+1} (again we omit the fine details of checking that all of the conditions are satisfied).

\item \emph{Suppose that \m{\alpha>0} and \m{\beta \in \{1,\ldots,\alpha\}}.} This is the final and most involved case. In this case we are not able to give an explicit formula for the constants \m{(C_{-1},C_0,\ldots,C_\alpha)}, instead we can only show that such constants exist. This is also the first case in which we must choose a non-zero \m{C_{-1}}. 
Let \m{A} be the  \m{(2\TEMPORARYI-\alpha)\times(\alpha+1)} matrix with entries  indexed by \m{0\LEQSLANT w<2\TEMPORARYI-\alpha} and \m{0\LEQSLANT j\LEQSLANT\alpha} and defined by %
\[\textstyle A_{w,j} = \textstyle p^{-\VERIFYIF{j=0}} \sum_{i>0} \binom{r-\alpha+j}{i(p-1)+j} \binom{i(p-1)}{w},\] for all \m{0\LEQSLANT w<2\TEMPORARYI-\alpha} and \m{0\LEQSLANT j\LEQSLANT\alpha}. %
 The matrix \m{A} is precisely the matrix obtained from \m{\UUU} by dividing the first column by \m{p}, and we have
\begin{align*} & \textstyle{}A (C_0,C_1/p,\ldots,C_\alpha/p)^T \\ & \textstyle{}\null\qquad=((\TEMPORARYD_0(D_\DARKCIRC)-C_{-1})/p,\TEMPORARYD_1(D_\DARKCIRC)/p,\ldots,\TEMPORARYD_{2\TEMPORARYI-\alpha-1}(D_\DARKCIRC)/p)^T.\end{align*} Let \m{A^{\subm}} denote the \m{(\alpha+1) \times (\alpha+1)} submatrix consisting of the top \m{\alpha+1} rows of \m{A}, \IE of those entries indexed by  \m{0\LEQSLANT w,j\LEQSLANT\alpha}. We have
 \begin{align*}
        A_{w,0} & =\textstyle \frac{(-1)^w(s-r)(r-\alpha)_w}{(s-\alpha)_{w+1}} + \BIGO{p}, 
        \\ A_{w,j} & = \textstyle \sum_v  \binom{-j}{w-v} \binom{r-\alpha+j}{v} \left(\binom{s-\alpha+j-v}{j-v}-\binom{r-\alpha+j-v}{j-v}\right) +\BIGO{p},
\end{align*}
for all \m{0\LEQSLANT w<2\TEMPORARYI-\alpha} and \m{0\LEQSLANT j\LEQSLANT\alpha},
so in particular \m{A} has integer entries. However, these expressions for the entries of \m{A} are not useful as we need to compute \m{A} up to precision \m{\BIGO{p^{t+1}}}.
Recall the definition of \[\textstyle \binom{\FORMX}{n}^\partial = \frac{\partial}{\partial \FORMX}\binom{\FORMX}{n} \in \QQ_p[\FORMX]\] from section~\ref{secAss3}. Let us also consider the \m{(2\TEMPORARYI-\alpha)\times(\alpha+1)} matrices \m{S} and \m{N} with integer entries indexed by \m{0\LEQSLANT w<2\TEMPORARYI-\alpha} and \m{0\LEQSLANT j\LEQSLANT\alpha} and defined by
    \begin{align*}
        S_{w,j} & = \textstyle p^{-\VERIFYIF{j=0}} \sum_{i=1}^{\beta} \binom{s+\beta(p-1)-\alpha+j}{i(p-1)+j} \binom{i(p-1)}{w}, \\
       N_{w,j} & = \textstyle p^{-\VERIFYIF{j=0}}   \sum_v   \binom{-j}{w-v} \binom{s+\beta(p-1)-\alpha+j}{v}^\partial  \sum_{i=0}^{\beta} \binom{s+\beta(p-1)-\alpha+j-v}{i(p-1)+j-v} \\
        & \textstyle \null\qquad - \VERIFYIF{w=0}%
        \textstyle  \binom{s+\beta(p-1)-\alpha+j} {j}^\partial
        -\VERIFYIF{j=0}\binom{s+\beta(p-1)-\alpha}{w}  \frac{(-1)^{w}w!}{(s-\alpha)_{w+1}},
\end{align*} for all \m{0\LEQSLANT w<2\TEMPORARYI-\alpha} and \m{0\LEQSLANT j\LEQSLANT\alpha}. The fact that \m{S} has integer entries can be shown in the same way as the fact that \m{A} has integer entries, and \m{N} has integer entries because %
\begin{align*} & \textstyle{}\sum_v   \binom{0}{w-v} \binom{s+\beta(p-1)-\alpha}{v}^\partial  \sum_{i=0}^{\beta} \binom{s+\beta(p-1)-\alpha-v}{i(p-1)-v} \\ & \textstyle{}\null\qquad  =  \binom{s+\beta(p-1)-\alpha}{w}^\partial  \TEMPORARYM_{s+\beta(p-1)-\alpha-w,-w} 
\\ & \textstyle{}\null\qquad  =  \binom{s+\beta(p-1)-\alpha}{w}^\partial  \sum_{i=0}^{w} (-1)^i \binom{w}{i} \TEMPORARYM_{s+\beta(p-1)-\alpha-i} 
\\ & \textstyle{}\null\qquad  =  \binom{s+\beta(p-1)-\alpha}{w}^\partial  \sum_{i=0}^{w} (-1)^i \binom{w}{i} (1 + [s-\alpha-i=0] + \BIGO{p}) 
\\ & \textstyle{}\null\qquad  =  \binom{s+\beta(p-1)-\alpha}{w}^\partial  \sum_{i=0}^{w} (-1)^i \binom{w}{i} (1 + \BIGO{p}) 
\\ & \textstyle{}\null\qquad  = \VERIFYIF{w=0}\binom{s+\beta(p-1)-\alpha}{w}^\partial +\BIGO{p}  = \BIGO{p}.\end{align*} The second equality follows from lemma~\ref{part1X5}, equation~\eqref{ceqn3},  the third  equality follows from lemma~\ref{part1X5}, equation~\eqref{ceqn2}, and the fourth equality follows from  \m{s\GEQSLANT 2\TEMPORARYI>\alpha+w}. %
 Let also \m{S^{\subm}} and \m{N^{\subm}} denote the \m{(\alpha+1) \times (\alpha+1)} submatrices  consisting of the top \m{\alpha+1} rows of  \m{S} and \m{N}, respectively.
Recall that \m{\TEMPORARYJ = u_0p^t}. Let us prove the following claim.

\emph{Approximation claim.}\label{approx}
\[\textstyle A = S + \TEMPORARYJ N + \BIGO{\TEMPORARYJ p}.\]

\emph{Proof of the approximation claim. } First let us look at the entries with \m{j>0}. Due to lemma~\ref{part1X5}, equation~\eqref{ceqn6}, \begin{align*}
    \textstyle A_{w,j} & = \textstyle \sum_v   \binom{-j}{w-v} \binom{r-\alpha+j}{v} \sum_{i>0} \binom{r-\alpha+j-v}{i(p-1)+j-v}
    \\  & = \textstyle \sum_v  \binom{-j}{w-v} \binom{r-\alpha+j}{v} \sum_{i} \binom{r-\alpha+j-v}{i(p-1)+j-v} %
    -\binom{r-\alpha+j}{j}\binom{0}{w}.
\end{align*}
The second equality here simply amounts to extracting the ``\m{i=0}'' term from the sum. We now use the equation
\begin{align}\label{wfru984uf38]uh}\textstyle \binom{r-\alpha+j}{j} = \binom{s+\beta(p-1)-\alpha+j}{j} + \TEMPORARYJ \binom{s+\beta(p-1)-\alpha+j}{j}^\partial + \BIGO{\TEMPORARYJ p},\end{align} which is true because \m{j<p} %
 (see %
 section~\ref{secNot}). Due to lemma~\ref{part1X5}, equation~\eqref{ceqn6},
\[\textstyle S_{w,j} = \textstyle \sum_v \binom{-j}{w-v} \binom{s+\beta(p-1)-\alpha+j}{v} \sum_{i>0} \binom{s+\beta(p-1)-\alpha+j-v}{i(p-1)+j-v}.\]
By combining these facts we get that \m{A_{w,j}-S_{w,j}-\TEMPORARYJ N_{w,j}} is
\begin{align*}
    & \textstyle \sum_v \binom{-j}{w-v} \binom{s+\beta(p-1)-\alpha+j}{v} \sum_{i} \left( \binom{r-\alpha+j-v}{i(p-1)+j-v}-\binom{s+\beta(p-1)-\alpha+j-v}{i(p-1)+j-v}\right) \\  &\textstyle  \null\qquad 
    + \TEMPORARYJ \sum_v \binom{-j}{w-v} \binom{s+\beta(p-1)-\alpha+j}{v}^\partial
   \sum_{i} \left( \binom{r-\alpha+j-v}{i(p-1)+j-v}-\binom{s+\beta(p-1)-\alpha+j-v}{i(p-1)+j-v}\right)
    \\  &\textstyle  \null\qquad
    -\binom{r-\alpha+j}{j}\binom{0}{w}+\binom{s+\beta(p-1)-\alpha+j}{j}\binom{0}{w}+\TEMPORARYJ\binom{s+\beta(p-1)-\alpha+j}{j}^\partial\binom{0}{w}.
\end{align*}
The first two lines are in \m{\BIGO{\TEMPORARYJ p}} since, due to lemma~\ref{part1X5}, equation~\eqref{ceqn1},
\[\textstyle  \sum_{i}   \binom{r-\alpha+j-v}{i(p-1)+j-v} \equiv \sum_{i}   \binom{s+\beta(p-1)-\alpha+j-v}{i(p-1)+j-v} \bmod{\TEMPORARYJ p}.\]
The third line is in \m{\BIGO{\TEMPORARYJ p}} because of equation~\eqref{wfru984uf38]uh}. Thus
\[A_{w,j}=S_{w,j}+\TEMPORARYJ N_{w,j}+\BIGO{\TEMPORARYJ p}\]
for all \m{0\LEQSLANT w<2\TEMPORARYI-\alpha} and \m{0< j\LEQSLANT\alpha}. Now let us look at the entries with \m{j=0}. In this case
\begin{align*}
    pA_{w,0} & \textstyle \vphantom{{\binom{s+\beta(p-1)-\alpha}{w}^\partial}
    \sum_{i>0} \binom{s+\beta(p-1)-\alpha-w}{i(p-1)-w}%
    - \binom{s+\beta(p-1)-\alpha}{w} \frac{(-1)^{w}w!}{(s-\alpha)_{w+1}}p}= \sum_{i>0} \binom{r-\alpha}{i(p-1)} \binom{i(p-1)}{w} %
    = \binom{r-\alpha}{w} \sum_{i>0}  \binom{r-\alpha-w}{i(p-1)-w},
    \\
    pS_{w,0} & \textstyle \vphantom{{\binom{s+\beta(p-1)-\alpha}{w}^\partial}
    \sum_{i>0} \binom{s+\beta(p-1)-\alpha-w}{i(p-1)-w}%
    - \binom{s+\beta(p-1)-\alpha}{w} \frac{(-1)^{w}w!}{(s-\alpha)_{w+1}}p}= \sum_{i>0} \binom{s+\beta(p-1)-\alpha}{i(p-1)} \binom{i(p-1)}{w} %
    = \binom{s+\beta(p-1)-\alpha}{w} \sum_{i>0}  \binom{s+\beta(p-1)-\alpha-w}{i(p-1)-w},
    \\
    p N_{w,0} & \textstyle \vphantom{{\binom{s+\beta(p-1)-\alpha}{w}^\partial}
    \sum_{i>0} \binom{s+\beta(p-1)-\alpha-w}{i(p-1)-w}%
    - \binom{s+\beta(p-1)-\alpha}{w} \frac{(-1)^{w}w!}{(s-\alpha)_{w+1}}p}= {\binom{s+\beta(p-1)-\alpha}{w}^\partial}
    \sum_{i>0} \binom{s+\beta(p-1)-\alpha-w}{i(p-1)-w}%
    - \binom{s+\beta(p-1)-\alpha}{w} \frac{(-1)^{w}w!}{(s-\alpha)_{w+1}}p.
\end{align*}
So \m{p(A_{w,0}-S_{w,0}-\TEMPORARYJ N_{w,0})} is equal to
\begin{align*}
    & \textstyle \binom{r-\alpha}{w} \sum_{i>0}  \left(\binom{r-\alpha-w}{i(p-1)-w}-\binom{s+\beta(p-1)-\alpha-w}{i(p-1)-w}\right)
    \\  &\textstyle  \null\qquad
    +\left(\binom{r-\alpha}{w}-\binom{s+\beta(p-1)-\alpha}{w}-\TEMPORARYJ\binom{s+\beta(p-1)-\alpha}{w}^\partial\right)\sum_{i>0} \binom{s+\beta(p-1)-\alpha-w}{i(p-1)-w}
    \\  &\textstyle  \null\qquad +
    \binom{s+\beta(p-1)-\alpha}{w} \frac{ (-1)^{w}w!}{(s-\alpha)_{w+1}}\TEMPORARYJ p.
\end{align*} 
Due to lemma~\ref{part1X5}, equation~\eqref{ceqn2} in combination with lemma~\ref{part1X5}, equation~\eqref{ceqn3}, we have
 \[\textstyle\sum_{i>0} \binom{s+\beta(p-1)-\alpha-w}{i(p-1)-w} = \BIGO{p},\] so the second line is in \m{\BIGO{\TEMPORARYJ p^2}}. 
 Since \[\textstyle{}s_{r-\alpha-j} = s_{s+\beta(p-1)-\alpha-j}=s-\alpha-j\] (in the notation of lemma~\ref{part1X5}, equation~\eqref{ceqn2}), and since  
\begin{align*}
   \TEMPORARYM_u \textstyle  & =\textstyle 1+ \VERIFYIF{s_u=p-1} + \frac{1}{p-1}\sum_{\MU\in\FF_p\backslash\{-1,0\}} (1+[\MU])^{s_u} \sum_{j=1}^{t_u} \binom{t_u}{j} p^j x_\MU^j
\end{align*} (this equation is shown in the proof of lemma~\ref{part1X5}, equation~\eqref{ceqn2}),
we can write \m{\TEMPORARYM_{r-\alpha-j}-\TEMPORARYM_{s+\beta(p-1)-\alpha-j}} as
 \begin{align*}
    \textstyle  & \textstyle  \frac{1}{p-1}\sum_{\MU\in\FF_p\backslash\{-1,0\}} (1+[\MU])^{s-\alpha-j} \sum_{j=1}^{\infty} \left[\binom{t_{r-\alpha-j}}{j}-\binom{t_{s+\beta(p-1)-\alpha-j}}{j}\right] p^j x_\MU^j.
\end{align*}
 Let \m{b = t_{s+\beta(p-1)-\alpha-j}}, so that \[\textstyle{}t_{r-\alpha-j} = b + \frac{r-s-\beta(p-1)}{p-1} = b - \epsilon(1+h_0p)\] for some \m{h_0\in \ZZ_p}. Then
 \begin{align*}\textstyle{} \binom{t_{r-\alpha-j}}{j}-\binom{t_{s+\beta(p-1)-\alpha-j}}{j} & \textstyle{} = \binom{b - \epsilon(1+h_0p)}{j}-\binom{b}{j}  \\ & \textstyle{} = -\epsilon (1+h_0p) \binom{b}{j}^\partial + \BIGO{\epsilon^2p^{-\PADICVALUATION(j!)}}\end{align*}
 for all \m{j\in\ZZ_{\GEQSLANT1}} (because \m{p^{\PADICVALUATION(j!)}\binom{\FORMX}{j}\in \ZZ_p[\FORMX]}). 
 Since  \m{j-\PADICVALUATION(j!)\GEQSLANT 2} for  \m{j\in\ZZ_{\GEQSLANT2}} (as shown in the proof of lemma~\ref{part1X5}, equation~\eqref{ceqn2}), we have 
 \begin{align*}
    \textstyle \TEMPORARYM_{r-\alpha-j}-\TEMPORARYM_{s+\beta(p-1)-\alpha-j}  & \textstyle =  -A_{s-\alpha-j} \binom{b}{1}^\partial \epsilon p  + \BIGO{\epsilon p^2} =   \frac{-\epsilon p}{s-\alpha-j} + \BIGO{\epsilon p^2}
\end{align*}
(where \m{A_{s-\alpha-j}} is as in the proof of lemma~\ref{part1X5}, equation~\eqref{ceqn2}). Thus 
\begin{align*}
   & \textstyle \sum_{i>0}  \left(\binom{r-\alpha-w}{i(p-1)-w}-\binom{s+\beta(p-1)-\alpha-w}{i(p-1)-w}\right) \\& \textstyle   \null\qquad  \vphantom{\sum_{j=0}^w (-1)^j \binom{w}{j} \left(\frac{t_{r-\alpha-j}}{s_{r-\alpha-j}}p-\frac{t_{s+\beta(p-1)-\alpha-j}}{s_{s+\beta(p-1)-\alpha-j}}p + \BIGO{\TEMPORARYJ p^2}\right)} = \sum_{j=0}^w (-1)^j \binom{w}{j} (\TEMPORARYM_{r-\alpha-j}-\TEMPORARYM_{s+\beta(p-1)-\alpha-j})
    \\  &\textstyle   \null\qquad
    \vphantom{\sum_{j=0}^w (-1)^j \binom{w}{j} \left(\frac{t_{r-\alpha-j}}{s_{r-\alpha-j}}p-\frac{t_{s+\beta(p-1)-\alpha-j}}{s_{s+\beta(p-1)-\alpha-j}}p + \BIGO{\TEMPORARYJ p^2}\right)}= \sum_{j=0}^w (-1)^j \binom{w}{j} \left(\frac{-\TEMPORARYJ p}{s-\alpha-j} + \BIGO{\TEMPORARYJ p^2}\right)
        \\  &\textstyle  \null\qquad
    \vphantom{\sum_{j=0}^w (-1)^j \binom{w}{j} \left(\frac{t_{r-\alpha-j}}{s_{r-\alpha-j}}p-\frac{t_{s+\beta(p-1)-\alpha-j}}{s_{s+\beta(p-1)-\alpha-j}}p + \BIGO{\TEMPORARYJ p^2}\right)}= -  \frac{(-1)^{w} w!}{(s-\alpha)_{w+1}}\TEMPORARYJ p + \BIGO{\TEMPORARYJ p^2}.
\end{align*}
 The third equality follows from the equation
\[\textstyle{} \sum_{j=0}^w (-1)^j \binom{w}{j} \frac{1}{s-\alpha-j} =  \frac{(-1)^{w} w!}{(s-\alpha)_{w+1}}.\] Thus \m{p(A_{w,0}-S_{w,0}-\TEMPORARYJ N_{w,0})} is equal to \begin{align*} & \textstyle 
\textstyle  \frac{(-1)^{w} w!}{(s-\alpha)_{w+1}}\left(\binom{s+\beta(p-1)-\alpha}{w}-\binom{r-\alpha}{w}\right)\TEMPORARYJ p + \BIGO{\TEMPORARYJ p^2}=\BIGO{\TEMPORARYJ p^2}.\end{align*} Therefore we have
\[A_{w,j}=S_{w,j}+\TEMPORARYJ N_{w,j}+\BIGO{\TEMPORARYJ p}\]
for all \m{0\LEQSLANT w<2\TEMPORARYI-\alpha} and \m{0\LEQSLANT j\LEQSLANT\alpha}, implying that
\[\textstyle A = S + \TEMPORARYJ N + \BIGO{\TEMPORARYJ p}.\] \hspace{\stretch{1}} \m{\square}

Let
\m{B = B_\alpha} be the 
\m{(\alpha+1)\times(\alpha+1)} matrix defined in lemma~\ref{lemmaComb2}. That lemma implies that %
 \[\textstyle BS^{\subm} = \left([w\in\{1,\ldots,\beta\}]\, p^{-\VERIFYIF{j=0}}\binom{s+\beta(p-1)-\alpha+j}{w(p-1)+j}\right)_{0\LEQSLANT w,j\LEQSLANT\alpha}.\]
 Let \m{\overline{Q}} be the matrix that is obtained from \m{\overline{BN^{\subm}}} by replacing the rows indexed \m{1,\ldots,\beta} with the corresponding rows of \m{\overline{BS^{\subm}}}. If \m{i \in \{1,\ldots,\beta\}} then we can write \m{\overline{Q}_{i,j}} as the reduction modulo \m{p} of
\[\textstyle p^{-\VERIFYIF{j=0}}\binom{s+\beta(p-1)-\alpha+j}{i(p-1)+j},\]
which can be simplified as
\[\textstyle \binom{\beta}{i} \cdot \left\{\begin{array}{rl}
    \binom{s-\alpha-\beta+i}{i}^{-1} (-1)^{i+1} & \text{if } j = 0, \\
    \binom{s-\alpha-\beta+j}{j-i} & \text{if }j>0.
\end{array}\right.\]
If \m{i\not\in\{1,\ldots,\beta\}} then we can write \m{\overline{Q}_{i,j}} as the reduction modulo \m{p} of
    \begin{align*}
       & \textstyle p^{-\VERIFYIF{j=0}}  \sum_{w=0}^\alpha B_{i,w} \sum_v \binom{-j}{w-v} \binom{s+\beta(p-1)-\alpha+j}{v}^\partial  
        \sum_{u=0}^{\beta} \binom{s+\beta(p-1)-\alpha+j-v}{u(p-1)+j-v}
        \\
        & \textstyle \null\qquad - \VERIFYIF{i=0} %
        \textstyle  \binom{s+\beta(p-1)-\alpha+j} {j}^\partial %
    -\VERIFYIF{j=0} \sum_{w=0}^\alpha B_{i,w} \binom{s+\beta(p-1)-\alpha}{w}  \frac{(-1)^{w}w!}{(s-\alpha)_{w+1}}.
\end{align*}
Thus \m{\overline{Q}} is the matrix \m{M} from lemma~\ref{lemmaComb4}, and that lemma implies that there is a solution of
\[\overline{Q} (z_0,\ldots,z_\alpha)^T = (1,0,\ldots,0)^T\] %
 such that \m{z_0 \ne 0}. 
Let \m{\overline u = (z_0,\ldots,z_\alpha)^T}. Then \m{\overline{u}} is %
  in \m{\Ker \overline{BS^{\subm}} = \Ker \overline {S^{\subm}}}.
  Let \m{\widetilde{u}} be any lift of \m{\overline{u}} with entries in \m{\ZZ_p}, so that \m{{BS^{\subm}}\widetilde{u}} is zero outside the entries indexed \m{1,\ldots,\beta}, and any non-zero entries of \m{BS^{\subm}\widetilde{u}} are \m{\BIGO{p}}. Because \m{{BS^{\subm}}} is zero outside the rows indexed \m{1,\ldots,\beta} and the submatrix of \m{{BS^{\subm}}} consisting of the rows indexed \m{1,\ldots,\beta} has full rank \m{\beta}, there is a solution of the equation
  \[\textstyle{} {BS^{\subm}}y = \frac1p BS^{\subm}\widetilde{u},\] for some vector \m{y} with entries in \m{\ZZ_p}. Then \m{u=\widetilde{u}-py} is a lift of \m{\overline{u}} such that \m{u\in \Ker {BS^{\subm}}=\Ker {S^{\subm}}}. Moreover, because the submatrix of \m{\overline{BS^{\subm}}} consisting of the rows indexed \m{1,\ldots,\beta} has full rank \m{\beta}, we have
 \[\textstyle \overline{BN^{\subm}} \overline u + \overline{BS^{\subm}} \overline v = (1,0,\ldots,0)^T\] for some \m{\overline v} with entries in \m{\FF_p}.
 Let \m{v} be any lift of \m{\overline{v}} with entries in \m{\ZZ_p}. Then \[\textstyle B (S^{\subm}+\TEMPORARYJ N^{\subm}) (u+\TEMPORARYJ v) = (\TEMPORARYJ,0,\ldots,0)^T + \BIGO{\TEMPORARYJ p}.\]
 Since \m{A = S+\TEMPORARYJ N + \BIGO{\TEMPORARYJ p}}, if we write \[\textstyle (C_0,C_1/p,\ldots,C_\alpha/p)^T = u+\TEMPORARYJ v,\] then \m{C_0} must be a unit since \m{z_0\ne 0}, and \m{C_j = \BIGO{p}} for \m{j>0}, and
\[\textstyle B A^{\subm} (C_0,C_1/p,\ldots,C_\alpha/p)^T = (\TEMPORARYJ,0,\ldots,0)^T + \BIGO{\TEMPORARYJ p}.\]
 Since the zeroth column of \m{B} is \m{(1,0,\ldots,0)^T},  the zeroth column of \m{B^{-1}} is also \m{(1,0,\ldots,0)^T}, and consequently we also have \begin{align*}& \textstyle A^{\subm} (C_0,C_1/p,\ldots,C_\alpha/p)^T \\ & \null\qquad \textstyle = (\TEMPORARYJ,0,\ldots,0)^T + \BIGO{\TEMPORARYJ p} \\ & \null\qquad \textstyle = ((\TEMPORARYD_0(D_\DARKCIRC)-C_{-1})/p,\TEMPORARYD_1(D_\DARKCIRC)/p,\ldots,\TEMPORARYD_\alpha(D_\DARKCIRC)/p)^T.\end{align*}
 Let us choose \m{C_{-1}=-\TEMPORARYJ p}, so that
 \[\textstyle{} (\TEMPORARYD_0(D_\DARKCIRC)/p,\TEMPORARYD_1(D_\DARKCIRC)/p,\ldots,\TEMPORARYD_\alpha(D_\DARKCIRC)/p)^T = \BIGO{\TEMPORARYJ p}.\] Then \m{\TEMPORARYD_\alpha(D_\DARKCIRC) = \BIGO{p^{t+2}}}, which makes \[\textstyle \textstyle \PADICVALUATION(\TEMPORARYD^{\prime})=\PADICVALUATION(C_{-1})=t+1.\] We moreover have
 \[\textstyle\TEMPORARYD_w(D_\DARKCIRC) = \BIGO{p^{t+2}}\] for all \m{ w \in \{0,\ldots, \alpha\}}.
 Since \[\textstyle{} A \equiv  S \bmod\TEMPORARYJ\] (because the entries of \m{ N} are integers), and since all rows of \m{ S} are linear combinations of the rows of \m{BS} indexed \m{1,\ldots,\beta}, it follows that every entry of \m{ A(C_0,C_1/p,\ldots,C_\alpha/p)^T} is \m{\BIGO{\TEMPORARYJ}}.\footnote{\label{foot1}Note that it is crucial here that the entries of \m{N} be integers. In this case they are since all equations in the proof of the approximation claim hold true for \m{w < s-\alpha}. If \m{w \GEQSLANT s-\alpha} then the equation for \m{N_{w,j}} is different---as, for instance,
 \[\textstyle{} \sum_{j=0}^w (-1)^j \binom{w}{j} \frac{1}{s-\alpha-j} \ne  \frac{(-1)^{w} w!}{(s-\alpha)_{w+1}}.\]%
   This is one of the places where the proof breaks down if \m{s < 2\TEMPORARYI}.}
 Therefore we have
\[\textstyle \TEMPORARYD_w(D_\DARKCIRC) = \BIGO{p^{t+1}}\] for all \m{w \in \{\alpha+1,\ldots, 2\TEMPORARYI-\alpha-1\}}. %
   So the constants \m{(C_{-1},C_0,\ldots,C_\alpha)} are suitable for \m{v=t+1} (we omit the fine details of checking that all of the conditions are satisfied).
\end{itemize}
Therefore we can always find  suitable constants \m{(C_{-1},C_0,\ldots,C_\alpha)} and suitable \m{v\in\QQ_p}, and that concludes the proof of proposition~\ref{mainPROP}.
\hspace{\stretch{1}}     {\BLACKSQUARE}

\section{Final remarks}
\label{secFinal}

The proof in section~\ref{secProof2} gives a partial classification of \m{\overline{V}_{k,a}} which we %
 summarize in this section. %

       \begin{theorem}\label{Ars1THM}
 Suppose that \m{\PADICVALUATION(a) \not\in \ZZ} and \m{\TEMPORARYI=\lceil\PADICVALUATION(a)\rceil} is such that
 \[\textstyle k \not\equiv 3,4,\ldots,2\TEMPORARYI,2\TEMPORARYI+1 \bmod p-1.\]
 Then \m{\overline{V}_{k,a} \cong \ind(\omega_2^{\MORE{k-2\TEMPORARYI+1}}) \otimes \omega^{\TEMPORARYI-1}} if 
 \[\textstyle{}\PADICVALUATION(k-s-\beta(p-1)-2)=0 \text{ for  all } \beta \in \{0,\ldots,\TEMPORARYI-2\},\]
 and \m{\overline{V}_{k,a}} is isomorphic to one of the representations in the set
 \[\textstyle{} \left\{\ind(\omega_2^{\MORE{k-2\TEMPORARYI+2l+1}}) \otimes \omega^{\TEMPORARYI-l-1} \,\big|\, l \in \{0,\ldots,\min\{j,\TEMPORARYI-\beta-1\}\}\right\}\]
 if
 \[\textstyle{} \PADICVALUATION(k-s-\beta(p-1)-2)=j>0 \text{ for  some } \beta \in \{0,\ldots,\TEMPORARYI-2\}.\] 
\end{theorem}

\begin{proof}{Proof. }
We know by proposition~\ref{mainPROP} that \m{\overline \Theta_{k,a}} is irreducible, and we know from the proof of proposition~\ref{mainPROP} that it must be a factor of one of the representations \begin{align}\label{listeroij}\textstyle{} \COMPACTIND{KZ}{G}\SUBMODULE(1),\ldots,\COMPACTIND{KZ}{G}\SUBMODULE(\TEMPORARYI-1),\COMPACTIND{KZ}{G}\QUOTIENTI(\TEMPORARYI-1).\end{align} Let us recall that the proof of proposition~\ref{mainPROP} is based on finding  constants \m{(C_{-1},C_0,\ldots,C_\alpha)} and  \m{v\in\QQ_p} that satisfy the five conditions~\eqref{cond1}, \eqref{cond2}, \eqref{cond3}, \eqref{cond4}, \eqref{cond5} and such that either \m{\PADICVALUATION(\TEMPORARYD^{\prime})\LEQSLANT v^{\prime}} or \m{\PADICVALUATION(a)-\alpha<v}. The key point is that if \m{\PADICVALUATION(\TEMPORARYD^{\prime})\LEQSLANT v^{\prime}} then we can reach a stronger conclusion than the one in the statement of proposition~\ref{mainPROP}: that  there is some element  \m{\mathrm{gen} \in \SCR I_a} that
represents a generator of \m{\FILT N_\alpha}.  In particular, if for \m{\alpha\in \{1,\ldots,\TEMPORARYI-2\}} the constants we choose in the proof of proposition~\ref{mainPROP} satisfy \m{\PADICVALUATION(\TEMPORARYD^{\prime})\LEQSLANT v^{\prime}}, then we can eliminate \m{\SUBMODULE(\alpha)} from the list~\eqref{listeroij}. Suppose first that \[\textstyle{}\PADICVALUATION(k-s-\beta(p-1)-2)=0 \text{ for  all } \beta \in \{0,\ldots,\TEMPORARYI-2\}.\] Then, for all \m{\alpha\in \{1,\ldots,\TEMPORARYI-2\}}, the constants we choose come from the second case in the proof of proposition~\ref{mainPROP}. So  \m{v=1} and \m{v^{\prime}=\min\{\PADICVALUATION(a)-\alpha,1\}=1}, and \m{\PADICVALUATION(\TEMPORARYD^{\prime})=\PADICVALUATION(\TEMPORARYD_\alpha(D_\DARKCIRC))=1\LEQSLANT v^{\prime}}. Thus we can eliminate \[\textstyle{}\SUBMODULE(1),\ldots,\SUBMODULE(\TEMPORARYI-2)\] from the list~\eqref{listeroij}, implying that \m{\overline \Theta_{k,a}} is   a factor of one of the representations \begin{align*}\textstyle{} \COMPACTIND{KZ}{G}\SUBMODULE(\TEMPORARYI-1),\COMPACTIND{KZ}{G}\QUOTIENTI(\TEMPORARYI-1)\end{align*} and hence, due to the classification given in theorem~\ref{BERGERX1}, that \[\textstyle{}\overline{V}_{k,a} \cong \ind(\omega_2^{\MORE{k-2\TEMPORARYI+1}}) \otimes \omega^{\TEMPORARYI-1}.\] Now suppose that  \[\textstyle{} \PADICVALUATION(k-s-\beta(p-1)-2)=j>0 \text{ for  some } \beta \in \{0,\ldots,\TEMPORARYI-2\}.\] If \m{\alpha \in \{1,\ldots,\beta-1\}} then the constants we choose come from the second case in the proof of proposition~\ref{mainPROP}, so we can similarly eliminate \m{\SUBMODULE(\alpha)} from the list~\eqref{listeroij}. If \m{\beta=0} then the constants we choose come from the third case in the proof of proposition~\ref{mainPROP}, so
\m{\textstyle{}\PADICVALUATION(\TEMPORARYD^{\prime})=j+1} and  \m{v^{\prime}=\min\{\PADICVALUATION(a)-\alpha,j+1\}}. Thus if \m{\alpha < \TEMPORARYI-j-1} and \m{\beta=0} then  we can eliminate \m{\SUBMODULE(\alpha)} from the list~\eqref{listeroij}. Finally, if \m{\beta>0} and \m{\alpha \in \{\beta,\ldots,\TEMPORARYI-2\}} then the constants we choose come from the fourth case in the proof of proposition~\ref{mainPROP}, so \m{\textstyle{}\PADICVALUATION(\TEMPORARYD^{\prime})=j+1} and  \m{v^{\prime}=\min\{\PADICVALUATION(a)-\alpha,j+1\}}, implying that if \m{\alpha < \TEMPORARYI-j-1}  then  we can eliminate \m{\SUBMODULE(\alpha)} from the list~\eqref{listeroij}. To summarize, we can eliminate \m{\SUBMODULE(\alpha)} from the list~\eqref{listeroij} for all \[\textstyle{}\alpha \in \{1,\ldots,\beta-1\} \cup \{1,\ldots,\TEMPORARYI-j-2\}.\] So   \m{\overline \Theta_{k,a}} is   a factor of one of the representations \begin{align*}\textstyle{} \COMPACTIND{KZ}{G}\SUBMODULE(\max\{\TEMPORARYI-j-1,\beta\}),\ldots,\COMPACTIND{KZ}{G}\SUBMODULE(\TEMPORARYI-1),\COMPACTIND{KZ}{G}\QUOTIENTI(\TEMPORARYI-1)\end{align*} and hence, due to the classification given in theorem~\ref{BERGERX1}, \m{\overline{V}_{k,a}} is isomorphic to one of the representations in the set
 \[\textstyle{} \left\{\ind(\omega_2^{\MORE{k-2\TEMPORARYI+2l+1}}) \otimes \omega^{\TEMPORARYI-l-1} \,\big|\, l \in \{0,\ldots,\min\{j,\TEMPORARYI-\beta-1\}\}\right\}.\] 
\end{proof}

\end{document}